\documentclass[12]{amsart}
\usepackage{amsmath,amsthm}
\usepackage{latexsym}
\usepackage{amssymb,amscd}
\usepackage{xspace}
\usepackage{amscd}
\usepackage{enumerate}
\newtheorem{theo}{Theorem}[section]
\newtheorem{coro}[theo]{Corollary}

\newtheorem{prop}[theo]{Proposition}

\numberwithin{equation}{section}
\newtheorem{lemm}[theo]{Lemma}  
\newcommand{\te}{\mathcal{T}}
\newcommand{\Sym}{\operatorname{Sym}}
\newcommand{\Aut}{\operatorname{Aut}}
\newcommand{\e}{\mathrm{e}}
\newcommand{\vol}{\operatorname{Vol}}

\newcommand{\ML}{\mathcal{ML}}
\newcommand{\Mod}{\operatorname{Mod}}
\newcommand{\CC}{\beta}

\newcommand{\ch}{\S}

\newcommand{\Li}{\mathcal{L}}
\newcommand{\tw}{\operatorname{tw}}
\newcommand{\arccosh}{\operatorname{Arccosh}}
\newcommand{\arcsinh}{\operatorname{Arcsinh}}
\newcommand{\g}{{\bf g}}
\newcommand{\eLL}{L}
\newcommand{\LL}{\mathcal{L}}
\newcommand{\Stab}{\operatorname{Stab}}
\newcommand{\tr}{\operatorname{Tr}}
\begin{document}
\title{Counting Mapping Class group orbits on hyperbolic surfaces}
\author{Maryam Mirzakhani}
\maketitle

\begin{section}{Introduction}

Let $S_{g,n}$ be a surface of genus $g $ with $n$ marked points. Let $X$ be a complete hyperbolic metric on $S_{g,n}$ with $n$ cusps. Every isotopy class $[\gamma]$ of a closed curve $\gamma \in \pi_{1}(S_{g,n})$ contains a unique closed geodesic on $X$. 
 Let  $\ell_{\gamma}(X)$ denote the hyperbolic length of the geodesic 
representative of $\gamma$ on  $X$. In this paper, we study the asymptotic growth of the lengths of closed curves
of a fixed topological type on $S_{g,n}.$ More precisely, 
we consider the mapping class group $\Mod_{g,n}$ of self-diffeomorphisms of $S_{g,n}$ up to homotopy (fixing the marked points). 
There is a natural action of $\Mod_{g,n}$ on the set of isotopy classes of closed curves on $S_{g,n}$. 

By the definition, two 
closed curves $\gamma_{1}$ and $\gamma_{2}$ are of the same topological type if and only if there exists
$\g\in \Mod_{g,n}$ such that $\g \cdot [\gamma_{1}] = [\gamma_{2}].$
In this paper, we study the asymptotics of the counting function
$$
s_{X}(L,\gamma)= \#\{\beta \in \Mod_{g,n} \cdot
\gamma\;|\;\ell_{\beta}(X) \leq L\}
$$
as $L \rightarrow \infty.$ As an application, one can obtain the asymptotics of the growth of 
$s^{k}_{X}(L)$, the number of closed curves of length $\leq L$ on $X$ with at most $k$ self-intersections.
We also discuss properties of random pants decomposition of large length on $X$.
Both these results are based on ergodic properties of the earthquake flow on a natural bundle over the moduli space $\mathcal{M}_{g,n}$ of hyperbolic surfaces of genus $g$ with $n$ cusps. 
\subsection{Growth of lengths of closed geodesics} Let $c_{X}(L)$ be the number of primitive closed geodesics on $X$ of length 
$\leq L$. By work of Delsart, Huber, Selberg and Margulis, we have 
$$c_{X}(L) \sim \e^{L}/L$$ as $L \rightarrow \infty$ \cite{margulis}. 
However, the number $s_{X}(L)$ of simple closed geodesics of length $\leq L$ grows only polynomially in $L$ \cite{BS}. This growth depends on both topology and geometry of a hyperbolic surface 
\cite{Re}, \cite{R1}, \cite{MR}. By \cite{M:ANN}, if $\gamma$ is a {\it simple} closed curve we have
\begin{equation}\label{SC}
\lim\limits_{L \rightarrow \infty} \frac{s_{X}(L, \gamma)}{L^{6g-6+2n}}= n_{\gamma} \; \frac{B(X)}{b_{g,n}}.
\end{equation}
Here $n_{\gamma} \in {\Bbb Q}$ and $B(X)$ is a smooth proper
function of $X$; in other words, $B(X)$ goes to infinity as $X$ develops short closed geodesics. 

We remark that in general it is hard to discern the self-intersection number of a closed curve in $\pi_{1}(S_{g,n})$. 
Roughly speaking, a {\it random} closed geodesic of length $L$ has around $L^2$ self-intersections. There are different ways to make this statement precise 
 \cite{L:1}, \cite{L:2} and \cite{CL}.  Rivin obtained the asymptotics of $s^{1}_{X}(L)$
 \cite{R2}. Later, Sapir obtained bounds on $s^{k}_{X}(L)$ in terms of $k$ and $L$ even in case $k$ is a function of $L$ \cite{S2}.

\subsection{Main new results.}\label{N:R}
Here we briefly discuss the main results of this paper:\\
 
\noindent
{\bf I. Counting closed geodesics of a given combinatorial type.}

We investigate the orbit of an arbitrary closed curve $\gamma$ under the action of an arbitrary finite index subgroup $\Gamma \subset \Mod_{g,n}$ and prove the following result:

\begin{theo}\label{morege}
Let $X$ be a complete hyperbolic metric on $S_{g,n}$ and let $\gamma=\gamma_{1}+ \ldots+ \gamma_{m}$ where $\gamma_{i} \in \pi_{1}(S_{g,n})$.  
Then there exists $n_{\gamma} \in {\Bbb Q},$ and $b_{g,n} \in {\Bbb R}_+$ such that 

\begin{equation}\label{main:theo}
  \frac{ \#\{\alpha \in \Gamma \cdot \gamma\;|\;\ell_{\alpha}(X) \leq L\}}{ L^{6g-6+2n}} \sim [\Mod_{g,n}: \Gamma] \times n_{\gamma} \frac{B(X)}{b_{g,n}},
\end{equation}
as $L \rightarrow \infty.$\\

 \end{theo}

\noindent
{\bf Remarks.} \\

\noindent
{\bf 1.} In the statement of Theorem \ref{morege} for $\beta \in \Mod_{g,n} \cdot \gamma,$ we can write $\beta= \beta_1+\ldots+ \beta_m$ and we have $\ell_{\beta}(X) =\sum_{i=1}^{m} \ell_{\beta_i}(X).$ However, Theorem \ref{morege} holds for other notions of length: a useful example is the case when $m>1$ and in $(\ref{main:theo})$ we use $\widetilde{\ell}_{\beta}(X)= \max \{\ell_{\beta_{1}}(X),\ldots, \ell_{\beta_{m}}(X)\}.$ See \ch \ref{QR} ({\bf 8}). \\
  
\noindent
{\bf 2.} 
In this theorem $\gamma_{i}'$s might not be simple and they can intersect. For the proof in the case of simple closed curves and curves with one self intersection see \cite{M:ANN}
and \cite{R2}. We remark that we do not rely on the proof of $(\ref{SC})$ and this paper gives a different proof of the results in \cite{M:ANN}.

As in \cite{M:ANN}, if $\gamma$ does not fill the surface the number $n_{\gamma}$ is closely related to intersection pairings of tautological line bundles over the moduli spaces of hyperbolic surfaces with punctures. Otherwise, the trace identities are the reason $n_{\gamma} \in {\Bbb Q}$. See \ch \ref{idea-proof}  ({\bf 2}).\\

\noindent
{\bf 3.} The function $B(X)$ is a smooth proper 
function on the moduli space $\mathcal{M}_{g,n}$ of hyperbolic surfaces of genus $g$ with $n$ cusps.
Moreover $B(X)$ is integrable with respect to the Weil-Petersson volume form on $\mathcal{M}_{g,n}$ 
\cite{M:ANN} and 
$$b_{g,n}= \int_{\mathcal{M}_{g,n}} \;B(X) \;dX.$$
See \ch \ref{QR} ({\bf 1}), and \ch \ref{EPEF} for a geometric interpretation of $B(X).$

\noindent
{\bf 4.}  There are only finitely many isotopy classes of closed 
curves on $S_{g,n}$ with $k$ self intersections up to the action of the mapping class group. Therefore, 
summing $s_{X}(L,\gamma)$ over representatives of these orbits gives $s^{k}_{X}(L)$, and the
asymptotics
of the $s_{X}(L,\gamma)$'s determines the asymptotics of $s^{k}_{X}(L)$. However our proof does not give any information about the growth of the number of 
different topological types of closed curves with $k$ self-intersections (as a function of $k$). For results in this direction see \cite{S1}. \\

\noindent
{\bf 5.} A different version of this result for the case of $g=n=1$ was obtained recently by V. Erlandsson  and J. Souto \cite{ES} using different techniques. We remark that their result is stronger than what we prove in this case. \\  

As it was pointed to us by Souto, in view of Theorem $\ref{morege}$ the results obtained in \cite{ES} can be used to obtain the asymptotics of the geodesic flow invariant measures supported on geodesic representatives of the elements of the set $\{ \g \cdot \gamma\;  |\;  \g \in \Mod_{g,n} \}$ on $T^1(X)$.  This answers a question raised by Sarnak in case of $g=n=1$ \cite{S} (see also \ch \ref{C-V}). 
 For more details see $({\bf 1}),$ $({\bf 2})$ and $({\bf 3})$ in \ch \ref{QR}. \\
 
We will consider the case of $g=n=1,$ and $\gamma= \gamma_{1}+\gamma_2+ \gamma_3,$ where $\gamma_1$ and $\gamma_2$ are simple closed curves on $S_{1,1}$ with $i(\gamma_1, \gamma_2)=1,$ and $\gamma_{3}$ is the Dehn twist of $\gamma_1$ around $\gamma_2.$
  For the modular once punctured torus $X_0$, we can measure the length of $\beta \in \Mod_{1,1} \cdot \gamma$ by $\widetilde{\ell}_{\beta}(X_0)= \max \{\ell_{\beta_{1}}(X_0),\ell_{\beta_2}(X_0), \ell_{\beta_{3}}(X_0)\}.$ In this case, our counting problem is the same as the counting of the integral solutions of the Markoff equation $x^2+y^2+z^2=3xyz,$ using 
  the norm $|(x,y,z)|= \max\{x,y,z\}.$ This problem was first studied by Zagier \cite{Z:M}. See $\ch \ref{C-V}$ for details.
  
 For counting results on related equations see \cite{BK}, \cite{B} and \cite{HSZ}.
  For more on counting and strong approximation results for Markoff type affine cubic surfaces see \cite{BGS}. \\
  
 \noindent
{\bf II. Statistics of pants decompositions on a fixed hyperbolic surface}.
We also study the distribution of lengths and twists of curves in a random pants decomposition of $X$. 
Defining the twist parameter is more technical, so here we  
discuss the weaker results on the length distributions. 

Let $\mathcal{P}=\{\alpha_1,\ldots, \alpha_{3g-3+n}\}$ be a pants decomposition of $S_{g,n}$. We study the asymptotic distribution of
$$\{ (\ell_{\g\cdot \alpha_1} (X),\ldots, \ell_{\g\cdot \alpha_{3g-3+n}} (X)) \}_{\g \in \Mod_{g,n}}  \subset {\Bbb R}^ {3g-3+n} . $$

Given ${\bf g} \in \operatorname{Mod}_{g,n}$, $\g \cdot \mathcal{P}$ defines a set of numbers 
$$\{ (\ell_{\g\cdot \alpha_1} (X),\ldots, \ell_{{\bf g}\cdot \alpha_{3g-3+n}} (X))\}.$$
Let $$\Delta= \{(x_1,\ldots, x_k) | \sum x_{i}=1\} \subset {\Bbb R}_{+}^{3g-3+n}.$$
Then each element in $s_{X}(L, \mathcal{P})$ defines a point in $\Delta.$ Given $L$, let $\nu_{L}$ be the discrete measure on $\Delta$ corresponding to the elements of 
$s_{X}(L, \mathcal{P})$.

\begin{theo}\label{distribution:pants}
Let $X$ be a compact hyperbolic surface. 
Then 
$$\lim_{L \rightarrow \infty} \nu_{L}= \nu,$$
where for any $A \subset \Delta$
$$\nu(A)= \int _{\operatorname{Cone}(A)} x_1 \ldots x_{k} d_{x_{1}}\cdots d x_{k}.$$
\end{theo}

\noindent
{\bf Remarks.} 
The proof of this result rely on the ideas used in \cite{M:AB}. The same method could be used to study the $\Gamma$ orbit of a pair of pants, where $[\Mod_{g,n}: \Gamma] < \infty.$\\
For an equivalent form of this result in terms of counting pairs of pants on a surface see Theorem \ref{eq:theo2}.
Note that the limiting measure in Theorem $\ref{distribution:pants}$ does not depend on the topology of $\mathcal{P}.$
Theorem \ref{distribution:pants} is not used in the proof of Theorem $\ref{morege}.$ Even though here we use the ergodicity of the earthquake flow
\ch \ref{proofs}, both statements are corollaries of the ergodicity of the action of the mapping class group on the space of measured laminations. 
 See also \ch \ref{G-B} and the discussion in the beginning of \ch \ref{proof:main:idea}.

\subsection{Asymptotic properties of length functions and distribution of balls.}

Let $\mathcal{M}_{g,n}$ be the moduli space of complete hyperbolic
Riemann surfaces of genus $g$ with $n$ cusps.  
Given a closed curve $\gamma \in \pi_{1}(S_{g,n})$ the geodesic length function of $\gamma$ defined a function over the universal cover $\te_{g,n}$ of $\mathcal{M}_{g,n}.$
Instead of proving statements about the mapping class group orbit $\Mod_{g,n} \cdot \gamma$ of $\gamma$ on $X \in \te_{g,n},$ we study the discrete set 
$$\Gamma \cdot X \subset \te_{g,n},$$
where $[\Mod_{g,n}: \Gamma] < \infty.$

To simplify the notation, for now we assume that $\Gamma=\Mod_{g,n}$. We can define a cover
 $$\pi^{\gamma}: \mathcal{M}_{g,n}^{\gamma}= \te_{g,n}/ \operatorname{Stab}(\gamma) \rightarrow \mathcal{M}_{g,n}.$$ 

Given $L>0,$ the length function 
 $$\ell_{\gamma}: \mathcal{M}_{g,n}^{\gamma} \rightarrow {\Bbb R}$$
$$ X \rightarrow \ell_{\gamma}(X),$$
defines  
$${\bf B}_{\gamma}(L)=\{ Z \; | \ell_{\gamma}(Z)\leq L \} \subset \mathcal{M}_{g,n}^{\gamma}.$$

We study the growth of 
$\# (\g \in \Mod_{g,n}\; |\; \g \cdot X \cap  {\bf B}_{\gamma}(L))$
as $L \rightarrow \infty$ by investigating the asymptotic shape of ${\bf B}_{\gamma}(L)$. This allows us to obtain results on the distribution of $\pi^{\gamma}( {\bf B}_{\gamma}(L)) \subset \mathcal{M}_{g,n}$ as 
$L\rightarrow \infty.$ \\

We say a closed curve $\gamma$ if {\it filling} if there is no essential closed curve $\beta$ disjoint from $\gamma$ on $S_{g,n}.$ 
We remark that if $\gamma$ is filling then $\operatorname{Stab}(\gamma)$  is a finite group. 
In this case, ${\bf B}_{\gamma}(L)$ is a compact subset of  $\mathcal{M}_{g,n}^{\gamma}$ and $\te_{g,n}$ \cite{Bon}. In general ${\bf B}_{\gamma}(L)$ has finite volume with respect to the Weil-Petersson volume form. 

\subsection{Remark on the case of $\Gamma \not = \Mod_{g,n}$.} 
In general,  we will consider the cover $\mathcal{M}_{g,n}[\Gamma]$ of $\mathcal{M}_{g,n}$ corresponding to $\Gamma \subset \Mod_{g,n}.$
 Our result generalizes to this case: this is mainly because the ergodicity of the earthquake flow also holds for the bundle  $\mathcal{P}^{1}\mathcal{M}_{g,n}[\Gamma]$. 
 See Theorem $\ref{er}$.
However, to simplify the notation, we consider the case that $\Gamma= \Mod_{g,n}$ and $\gamma$ is filling in the introduction. 

\subsection{Idea of proof and additional results.}\label{idea-proof} 

Let  $\mathcal{P}$ be a pants decomposition of $S_{g,n}.$ For a marked hyperbolic surface 
$X \in \te_{g,n} $, the {\it Fenchel-Nielsen coordinates}  
associated with 
$\mathcal{P}$, 
\begin{equation}\label{FNC}
 \te_{g,n} \rightarrow  {\Bbb R}_{+}^{3g-3+n} \times {\Bbb R}^{3g-3+n}
\end{equation}

$$
X \rightarrow ((\ell_{\alpha_{i}}(X))_{i=1}^{3g-3+n}, (\tau_{\alpha_{i}}(X))_{i=1}^{3g-3+n}),
$$
consists of the 
set of lengths of 
all geodesics used in the decomposition and the set of the 
{\it twisting} parameters used to glue the pieces. See $\ch 2$ for more details. 
This way, one obtains a natural isomorphism between $\te_{g,n}$ and $ {\Bbb R}_{+}^{3g-3+n} \times {\Bbb R}^{3g-3+n}$.

We prove that the length functions 
$$\ell_{\gamma}: \te_{g,n} \rightarrow {\Bbb R}_{+},$$
$$ X \rightarrow \ell_{\gamma}(X)$$
behave asymptotically like a piecewise linear function in terms of the Fenchel-Nielsen coordinates.  Note that the linear structure induced by the Fenchel-Nielsen coordinates does depend on the choice of $\mathcal{P}$. However, if we choose a different pants decomposition $\mathcal{P}'$ as we will see the change of coordinates is asymptotically piecewise linear (this happens only away from finitely many hyperplanes).   
Recall that the Weil-Petersson volume form $\mu_{wp}$ defined by $\bigwedge_{i=1}^{3g-3+n}  d\ell_{i} \wedge d \tau_{\alpha_i}$ is independent of the choice of $\mathcal{P}.$ See \ch $\ref{Background}$.\\

\noindent
{\bf Steps of proof of Theorem \ref{morege}.} \\

\noindent
{\bf 1.}  We show that for any $\gamma$, the function $\ell_{\gamma}$ {\it asymptotically} behaves like a piecewise linear function in terms of the Fenchel-Nielsen coordinates on $\te_{g,n}.$ In fact one can write down 
a formula for $\ell_{\gamma}$ in terms of the Fenchel-Nielsen coordinates. For an explicit calculation, see \cite{O:explicit}. We don't need this general formula here. Instead we use work of \cite{HT:moves} to get the result by only doing the calculations in four holed spheres and one holed tori.

 Recall that $X \in \te(S_{g,n})$ defines a representation $\rho: \pi_{1}(S_{g,n}) \rightarrow \operatorname{P SL}(2, {\Bbb R})$ (sending the loops around marked points to elliptic elements). 

Given $\gamma \in \pi_1(S_{g,n}),$ the trace of the matrix $\rho(\gamma)$ is related to the length of the corresponding closed geodesWeil-Peterssonic on $X$ by
$$2 \cosh(\frac{\ell_{\gamma}(X)}{2})= |\tr(\rho(\gamma))|.$$Weil-Petersson
 Let $\gamma_{1},\ldots, \gamma_{N}$ is a set of generators for 
$\pi_{1}(S_{g,n})$. For any $\gamma \in \pi_{1}(S_{g,n}) $ the trace of $\rho (\gamma)$ is a polynomial with integral coefficients in terms of 
$\{\operatorname{Tr}(\rho(\gamma_{i})\cdot \rho(\gamma_j)\cdot \rho(\gamma_{k}) \},$
where $1 \leq i <i <j<k \leq N.$ 
 
 One can prove the asymptotically linear behavior of $\ell_{\gamma}$ on $\te_{g,n}$ using these trace identities in $\operatorname{SL}(2, {\Bbb R})$.
See \ch \ref{C-V} for the discussion in case of once punctured torus. Also, see \cite{B:book}, \cite{G3} references within for more details.\\

\noindent
{\bf 2.} We show that apart from finitely many hyperplanes, as $L \rightarrow \infty$, ${\bf B}_{\gamma}(L)$ gets close to being a polyhedral shape with finitely many sides \ch \ref{PLL}. As  a result, we prove that 
$$  \int_{\mathcal{M}_{g,n}} s_{X}(L,\gamma) dX = \operatorname{Vol}_{wp}({\bf B}_{\gamma}(L)) \sim n_{\gamma} L^{6g-6+2n},$$
$L \rightarrow \infty$.  Here the volume is taken with respect to the Weil-Petersson volume form \ch \ref{vol-est}. This generalizes the result obtained in \cite{M:In}.\\
 Moreover, the coefficients of the linear functions estimating $\ell_{\gamma}$ are all in ${\Bbb Q}$ which implies that $n_{\gamma} \in {\Bbb Q}$ if $\gamma$ is filling. In general $n_{\gamma}$ can be written in terms of the coefficients of the linear functions approximating $\ell_{\gamma}$ and the leading coefficients of Weil-Petersson volumes of moduli spaces of hyperbolic surfaces with geodesic boundary components.  See Theorem $\ref{volm}$ and \ch \ref{vol-est}.\\

\noindent
{\bf 3.} Let $\mathcal{C}$ be a cone (with respect to the Fenchel-Nielsen coordinates) and $H \subset \mathcal{C} \subset \te_{g,n}$ be an open set on a hyperplane defined by a linear function.  
Given ${\bf t}>0$ and $X \in \te_{g,n}$ let ${\bf t} \cdot X \in \te_{g,n}$ be the surface obtained by multiplying all the coordinates of $X$ by ${\bf t}.$  
Consider the projection  $\pi: \te_{g,n} \rightarrow \mathcal{M}_{g,n}.$ We use the results motivated by \cite{M:AB} to obtain equidistribution results for  $\pi ({\bf t} \cdot H) \subset \mathcal{M}_{g,n}.$ 

Our method works only if there exists  a constant $M(\mathcal{C})$ such that 
\begin{equation}\tag{*}
 \frac{|\tau_{\alpha_i}(Y)|}{\ell_{\alpha_i}(Y)} \leq M(\mathcal{C})
 \end{equation}
for each $\alpha_i \in \mathcal{P},$ and $Y \in \mathcal{C}.$ In this case, we prove that as ${\bf t} \rightarrow \infty$, $\pi ({\bf t} \cdot H)$ becomes equidistributed with respect to the measure $B(X) \times \mu_{wp}$ as ${\bf t} \rightarrow \infty.$ Here $\mu_{wp}$ is the Weil-Petersson volume form on $\mathcal{M}_{g,n}.$

 Assume that $\mathcal{C} \subset \te_{g,n}$ satisfies $(*)$ and $\Li: \mathcal{C}  \rightarrow {\Bbb R}_{+}$ 
 is a linear function in terms of the Fenchel-Nielsen coordinates corresponding to $\mathcal{P}$. We prove that if $V_{\mathcal{C}}(\Li)= \operatorname{Vol}_{wp}(\{Z \in \mathcal{C}, \Li(Z) \leq 1\}) < \infty,$
then for any 
$X \in \te_{g,n}$, we have

\begin{equation}\label{countl}
\#\{ \g \in \Mod_{g,n} \; | \;\; \g \cdot X \in \mathcal{C}, \;   \Li(\g \cdot X) \leq L\} \sim V_{\mathcal{C}}(\Li)\cdot  L^{6g-6+2n} \cdot B(X), 
\end{equation}
as $L \rightarrow \infty.$

 In order to prove this result, first we use the ergodicity of the earthquake flow
\cite{M:EE} to prove (\ref{countl}) for $\Li_{0}=\sum_{i=1}^{3g-3+n} \ell_{\alpha_{i}}.$
 
The connection between equidistribution of  $\pi ({\bf t} \cdot H)$ and our counting result $(\ref{countl})$ is discussed in \ch \ref{from-eq-to-count}. The main idea goes back to Margulis' thesis \cite{margulis}. See also \cite{EM:mixing}.

 One important point in the argument is controlling the shape of balls of small radius centered at 
 $\g \cdot X$ in Fenchel-Nielsen coordinates \ch \ref{local:calculation}. 
 Here one could use any metric which is invariant under the mapping class group, but we use the Lipschitz distance function $d_{Th}$ defined by Thurston \cite{Thurston:minimal}. Roughly speaking, we need an upper bound on $$\frac{|\tau_{\alpha_i}(\g \cdot X)-\tau_{\alpha_i}(Y)|}{\epsilon \cdot \ell_{\alpha_i}(\g \cdot X)}$$ when $d_{Th}(\g \cdot X,Y) < \epsilon.$  These estimates do not hold everywhere on $\te_{g,n}.$ However, they work around the points in $\mathcal{B}_{\mathcal{P}} (M) \subset \te_{g,n}$ where for all  $1 \leq i,j \leq 3g-3+n\;,$ we have $\frac{\ell_{\alpha_{i}}(X)}{\ell_{\alpha_{j}}(X)} < M.$ See Lemma $\ref{local:main}$ and Corollary \ref{coro:est}.
 
However, this is not a problem for us as most points in $\Mod_{g,n} \cdot X \cap \mathcal{C}_{L}$ are in $\mathcal{B}_{\mathcal{P}} (M).$ We remark that our estimates get worse as $M \rightarrow \infty.$ 
See \ch \ref{basic:estimate:1} for more details.\\

\noindent 
{\bf 4.} We emphasize that the length function behaves like a linear function only away from finitely many hyperplanes. So we still have to show that the set of points 
$\Mod_{g,n} \cdot X$ do not accumulate around hyperplanes in $\te_{g,n}.$ Also, we need to consider the asymptotic behavior in infinitely many cones 
$ m \ell_{\alpha_{i}} \leq  \tau_{\alpha_i} \leq (m+1) \ell_{\alpha_{i}}.$ In order to prove our main result, we prove several estimates on the properties of random pants 
decompositions on a hyperbolic surface $X$.
For example, we show that as $L \rightarrow \infty,$ we have 
 \begin{itemize}
 \item 
 $$\frac{|\{\g \cdot \mathcal{P} \in  s_{X}(L, \mathcal{P}) \; | \min\{\ell_{\g \cdot \alpha_{i}}(X)\} \leq \sqrt{L} \}|}{L^{6g-6+2n}} \longrightarrow 0,$$
 and
 \item $$\frac{|\{\g \cdot \mathcal{P} \in s_{X}(L, \mathcal{P}) \; | \tau_{i}(X, \g \cdot \mathcal{P}) \leq \sqrt{L} \}|}{L^{6g-6+2n}} \longrightarrow 0. $$
\end{itemize}

More generally, if $\mathcal{R}$ is a linear function with positive coefficients, in \ch \ref{estimates} we obtain an upper bound on 
$$|\{\g \cdot \mathcal{P} \in  s_{X}(L, \mathcal{P}) \; | \mathcal{R}(\ell_{\alpha_{1}}(X), \ldots, \ell_{\alpha_{3g-3+n}}(X))\leq L'\}| . $$ 
These bounds are used in \ch \ref{from-eq-to-count}.
If $\mathcal{R}$ has both negative and positive coefficients, we can use $(\ref{countl})$, and the results of \ch \ref{from-eq-to-count} and \ch \ref{limiting} to get a weaker bound.
For example, we show:
$$\frac{|\{\g \cdot \mathcal{P} \in  s_{X}(L, \mathcal{P}) \; | |\ell_{\g \cdot \alpha_{i}}(X)+ \ell_{\g \cdot \alpha_{j}}(X)-\ell_{\g \cdot \alpha_{l}}(X)| \leq \sqrt{L} \}|}{L^{6g-6+2n}} \longrightarrow 0. $$
See Lemma $\ref{cone:count:upper}.$  These estimates might be of independent interest \ch \ref{estimates}. 

\subsection{Questions and Remarks}\label{QR}
\noindent
{\bf 1.} First, we briefly recall the results obtained in \cite{M:ANN}. For any Riemann surface $Y$ with bounded negative curvature and any multi-curve $\gamma$ we have  
 $$ s_{Y}(L,\gamma) \sim
 \frac{B(Y)}{b_{g,n}}\,
n_{\gamma}\; L^{6g-6+2n},$$
as $L \rightarrow \infty$.

A key role in the approach is played by the space $\ML_{g,n}$ of compactly supported measured laminations on $S_{g,n}$ (See \ch \ref{sec:ML}). The space $\ML_{g,n}$ is the completion of the set of rational multi-curves on $S_{g,n}.$ Moreover, $\ML_{g,n}$ carries a mapping class group invariant volume form $\mu_{Th}$.

Let $B_{X}\subset \ML_{g,n}$ be the unit ball in the space of measured geodesic 
laminations with respect to the length function at $X$, and 
$B(X)=\vol_{Th}(B_{X})$. 
Let $\mu_{\Gamma}^{\gamma}$ be the discrete measure on $\ML_{g,n}$ supported on the orbit 
$\gamma$, that is 
$$
\mu_{\Gamma}^{\gamma}= \sum\limits_{g \cdot \Gamma}
\delta_{g\cdot \gamma}.
$$
 Note that $\ML_{g,n}$ has a natural action of 
${\Bbb R}_{+} $ by dilation. 

For $T \in {\Bbb R}_{+}$, Let $T^{*}(\mu^{\gamma})$
denote the rescaling of $\mu^{\gamma}$ by factor $T$.
Although the action of 
$\operatorname{Mod}_{g,n}$ on 
$\ML_{g,n}$ is not linear, it is homogeneous.  
We define the measure $\mu_{T,\gamma, \Gamma}$ by
 \begin{equation}\label{measure}
 \mu_{T,\gamma, \Gamma}=\frac{T^{*}(\mu_{\Gamma}^{\gamma})}{T^{6g-6+2n}}.
 \end{equation} 
  So given $U \subset \ML_{g,n}$ 
 $\mu_{T,\gamma, \Gamma}(U)=\mu_{\Gamma}^{\gamma}(T\cdot U)/T^{6g-6+2n}$.

Then, for any $T>0$, the measure $\mu_{T,\gamma}$ is also invariant under the action of  $\operatorname{Mod}_{g,n}$
on $\ML_{g,n},$ and 
$$
\mu_{T,\gamma}(B_{X})=\frac{s_{X}(T,\gamma)}{T^{6g-6+2n}}.
$$

Hence understanding the asymptotic behavior of $s_{X}(T,\gamma)$ is
closely related to the asymptotic behavior of the sequence 
$\{\mu_{T,\gamma}\}_{T}.$ For any simple closed curve (or multi-curve), 
the method used in \cite{M:ANN} implies that for any finite index subgroup $\Gamma$ of 
$\Mod_{g,n}$ we have:
\begin{equation}\label{MANN}
\mu_{T,\gamma, \Gamma}\rightarrow  \frac{n_{\gamma}}{b_{g,n}}  [\Mod_{g,n}: \Gamma] \cdot \mu_{Th}.
\end{equation}
as $T\rightarrow \infty.$ 
We remark that  $(\ref{MANN})$ is a topological statement: the hyperbolic length can be replaced by any way of measuring the length of simple closed curves, as long as it extends continuously to 
$\ML_{g,n},$ and it is homogeneous with respect to the action of ${\Bbb R}_+$ on $\ML_{g,n}.$  \\
\noindent
{\bf 2.}  Note that similar to $(\ref{measure})$, for any closed curve $\gamma$ and $X$ one can consider
\begin{equation}\tag{**}
\mu^{X}_{T,\gamma, \Gamma}=\frac{ \sum\limits_{\ell_{\g \cdot \gamma}(X) \leq T}
\delta_{\g\cdot \gamma}}{T^{6g-6+2n}}
\end{equation}
 as a point in the space $\mathcal{G}(S_{g,n})$ of geodesic currents on $S_{g,n}$ \cite{Bon}. 
  On the other hand, every $X \in \te_{g,n}$ determines a measure $\nu_{X}$ on $\mathbb{P}\ML_{g,n}$ such that 
for $U \subset \mathbb{P}\ML$, we have                                                                        
\begin{equation}\label{pa}
\nu_{X}(U)=\nu(\{\eta\; \in \ML_{g,n} | \;\ell_{\eta}(X) \leq 1, [\eta] \in U\}.        
\end{equation}       
 
For $g=n=1$ and $\Gamma =\Mod_{1,1}$,  Erlandsson and Souto \cite{ES} show that for 
$X \in \te_{1,1}$ as $T \rightarrow \infty,$ the sequence $\mu^{X}_{T,\gamma, \Gamma}$ will converge to a multiple of $\nu_{X}$ for an arbitrary closed curve $\gamma$.

As Souto pointed to us, combining the results in \cite{ES} and Theorem \ref{morege} implies that for any $X \in \te_{g,n}$
 \begin{equation}\label{ES-general}
 \mu^{X}_{T,\gamma, \Gamma} \rightarrow \frac{n_{\gamma}}{B(X)} \; \nu_{X}
 \end{equation}
 as $T \rightarrow \infty.$ See Corollary $4.4$ in \cite{ES}.

\noindent
{\bf 3.}                                                      
One can reformulate $(\ref{ES-general})$ as the distribution of orbits of $\Mod_{g,n} \cdot \gamma$ on $T^1(X).$ 
Recall that every $\lambda \in \ML_{g,n}$ determines a geodesic flow invariant                  
measure $\mu_{\lambda}$ on $T^{1}(X).$
Then we get a measure $\mu_{X}$ on $T^{1}(X)$ defined by the ergodic
decomposition
$$\mu_{X}=\int \limits_{\mathbb{P}\ML_{g,n}} \frac{\mu_{\lambda}}{\ell_{\lambda}(X)} \; d\lambda.$$
In general,  the limit of measures  supported on the set of $\beta'$s in  $\Mod_{g,n} \cdot \gamma$ (as in $(**)$) exists and it is equal to a multiple of $\mu_{X}$.
 In case $\gamma$ does not have any self intersections this equidistribution result is a corollary of \cite{M:ANN}. 
 \\

\noindent
{\bf 4.} More generally, given any {\it filling} geodesic currents $Z_1$ and $Z_2$, we can consider the growth of the 

\begin{equation}\label{limit:topological}
\frac{\#\{Z \in \Mod_{g}  \cdot Z_1 \;|\; i(Z, Z_2 ) \leq L\}}{L^{6g-6}}
\end{equation}
as $L \rightarrow \infty.$ One can easily prove this is a sequence 
is bounded away from zero. 

We note that if $Y$ is filling, then any limit point of the set $[\Gamma \cdot Y] \subset \mathbb{P}(\mathcal{G}(S_{g}))$ is inside $\mathbb{P}(\ML_{g}).$
As a result of work of Souto and Erlandsson \cite{ES} (see (\ref{ES-general}))
 the limit in $(\ref{limit:topological})$ exists when $Z_1$ and $Z_2$ are both filling closed curves on 
$S_{g}.$  
 This is of special interest  since many geometric structures on $S_{g,n}$ could be represented as geodesic currents. See 
\cite{Otal}, \cite{DRL} and \cite{LM}. \\

\noindent
{\bf 5.} We emphasize that the method used in our paper only works for hyperbolic surfaces.  Our method does not
 work for studying the statistics of ratios of lengths of a random pants decomposition on a negatively curved surface. It would be interesting to generalize Theorem \ref{distribution:pants} to this setting. \\

\noindent
{\bf 6.}
 Some of the results in this paper are analogous to the ones in \cite{ABEM}. However, in our setting there is no natural flow playing the rule of the geodesic flow. The length balls ${\bf B}_{\gamma}(L)$ are not the metric balls (we only use some basic properties of the Thurston lipshitz distance in Lemma \ref{local:main} and \ch \ref{from-eq-to-count}).
 The ergodic properties of the Teichm\"ulcer geodesic flow are used only indirectly: the ergodicity of the earthquake flow is a consequence of the ergodcity of the horocycle flow on the moduli space of meromorphic quadratic differentials with simple poles. The ergodicity of the horocycle flow was obtained using the mixing of the geodesic flow due to Masur and Veech. See \cite{masur:I}, \cite{Masur:mapping}, \cite{veechteich}, \cite{M:EE} and references within. \\

\noindent
{\bf 7.} The results of this paper also hold in general for hyperbolic surfaces with geodesic boundary components. See \ch \ref{G-B}.
  The methods used in this paper could be used to obtain the distribution of
$$(\ell_{\g \cdot \gamma_1}(X), \ell_{\g \cdot \gamma_2}(X))_{\g \in \Mod_{g,n}}, $$ for arbitrary closed curves  
$\gamma_1, \gamma_2 \in \pi_1(S_{g,n}),$ and $X \in \te_{g,n}.$ 
 \\

\noindent
{\bf 8.} Let $F: \te_{g,n} \rightarrow {\Bbb R}_+$  be an arbitrary proper asymptotically piecewise linear map and $X \in \te_{g,n}$
The method in this paper can be used to study the growth of
$$\# \{ \g \cdot X\; |\; \g \in \Mod_{g,n}\;, F(\g \cdot X) \leq L\}$$
as $L \rightarrow \infty.$

\noindent
{\bf 9.} It seems to be hard to make our method effective and obtain error terms for $(\ref{main:theo})$. Even knowing how fast the {\it horosphere measures} ( see \ch \ref{limiting}) get equiduistributed on $\mathcal{M}_{g,n}$ does not seem to be enough.  \\

\subsection{Remark on the case with extra symmetry}
We have to be careful about the statement of the results when $\gamma$ or $X$ have extra symmetries.
Given a finite subgroup index $\Gamma$ of $\Mod_{g,n}$, we consider the following sets:
\begin{enumerate}
\item $\mathcal{A}_1 (X,\gamma, L)=\{\g \cdot \gamma \; | \;\g \in \Gamma, \ell_{\g \cdot \gamma} (X) \leq L\}$ (used to define $s_{X}(L, \gamma)$),  
\item $ \mathcal{A}_{2} (X,\gamma, L)= \{\g^{-1}\cdot  X\; |\; \ell_{\gamma} (g^{-1}\cdot X) \leq L, \g \in \Gamma\}, $ (related to the lattice counting problem $\Gamma \cdot X \cap  {\bf B}_{\gamma}(L)$), and
\item $\mathcal{A}_3 (X,\gamma, L)= \{\g \in \Gamma \; | \; \ell_{\g \cdot \gamma} (X) \leq L\}$ (considered in $(\ref{countl}$)).
\end{enumerate}
It is important to note that in general these sets are not of equal size. The problem arises if $X$ or $\gamma$ have 
non-trivial symmetry groups.
\begin{itemize}
\item 
For a closed curve $\gamma$, we define the symmetry group of $\gamma$ by $$\Sym(\gamma)= \Stab(\gamma)\bigcap_{i(\alpha, \gamma)=0} \Stab(\alpha),$$
where $\Stab(\beta)= \{\g \in \Mod_{g,n}\; | \; \g \cdot \beta=\beta \}.$ If $\gamma$ is filling, then $\Sym(\gamma)= \Stab(\gamma).$
\item
Given  $X \in \te_{g,n}$ we define
 $\operatorname{Aut}(X)=\{\g \in \Mod_{g,n}\; |\; \g \cdot X=X\}.$
 
 \end{itemize}

 It is known that both $\Sym(\gamma)$ and $\Aut(X)$
 are finite subgroups of $\Mod_{g,n}.$ \\

Also, one can easily check that:

$$ |\Sym(\gamma) \cap \Gamma| \cdot |\mathcal{A}_1 (X,\gamma, L)| =  |\mathcal{A}_3 (X,\gamma, L)|,$$ and 

$$ |\Aut(X) \cap \Gamma| \cdot |\mathcal{A}_2 (X,\gamma, L)| =  |\mathcal{A}_3 (X,\gamma, L)|.$$\\

\noindent
{\bf Notes and References.}
For more on properties of length functions see \cite{Wolpert:L}, \cite{BBFS} and references within.
For interesting conjectures and results on word length growth of self-intersecting curves see \cite{C},  \cite{CL}, and \cite{S1}.

The key tool in proving the results on earthquake flow is the  work of Thurston and Bonahon on horocycle foliations and shear coordinates
for the Teichm\"uller space and the space of measured foliations \cite{Thurston:minimal}, \cite{Bonahon:shear}. This construction is a natural generalization 
of the Fenchel-Nielsen coordinates for the Teichm\"uller space. For the dynamics of the horocycle flow see \cite{Masur:mapping} See also \cite {Re}.

 It is known that the intersection $i(\gamma,-)$ defines a piecewise linear function on $\ML_{g}.$ For more on the Dehn-Thurston coordinate change under elementary moves see \cite{P}, \cite{D-Thurston}. Compare with Proposition \ref{simple-case} and \ch \ref{HT:moves}. \\

\noindent
{\bf Acknowledgments}. I am grateful to Alex Eskin, Peter Sarnak, Jenya Sapir and Juan Souto for many illuminating and encouraging discussions related to the results of this paper. I would also like to thank Tian Yang and Alex Wright for their comments.  This work was partially supported by NSF and Simons grants.

\end{section}

\begin{section}{Background}\label{Background}
In this section, we present some familiar concepts on the basic properties of the moduli space, the space of measured geodesic 
laminations, and the earthquake flow. We will give an overview of the results we will use later in the paper. For more details on hyperbolic geometry
and Teichm\"uller theory see \cite{IT} and \cite{B:book} and references within. 

\subsection{Hyperbolic length function and trace identities}\label{trace:id}

Recall that every hyperbolic surface $X$ gives rise to a representation
$$\rho_X: \pi_{1}(S_{g,n}) \rightarrow PSL_{2}({\Bbb R})$$
$$ \gamma \rightarrow A_{\gamma}$$
such that 
$$|\tr(A_{\gamma})|=\cosh(\ell_{\gamma}(X)/2).$$

Generally, the conjugacy class of a matrix in $\operatorname{SL}(2, {\Bbb C})$ is determined by its trace. Also, a generic pair $(A,B)$ is determined up to conjugacy by the triple of traces $(\tr(A),\tr(B), \tr(AB))$ and we have  $$\tr([A,B])= \tr(A)^2+\tr(B)^2+ \tr(AB)^2 - \tr(A) \tr(B) \tr(AB).$$
Moreover, one can calculate the trace of any word in $A$ and $B$ using the following identity:
\begin{equation}\label{TI}
\tr(AB)+\tr(AB^{-1})= \tr(A) \times \tr(B).
\end{equation}
By a theorem of Fricke \cite{F1}, \cite{FK}, the moduli space of equivalence classes of $SL(2, {\Bbb C})$-representations naturally identifies with affine space ${\Bbb C}^3$.

 \subsection{Teichm\"uller Space.} 

A point in the {\it Teichm\"uller space} $\te(S)$
is a complete hyperbolic surface $X$ equipped with a diffeomorphism $f: S
\rightarrow X$. The map $f$ provides a {\it marking} on $X$ by $S$. Two
marked surfaces $f:\; S \rightarrow X$ and $g:\; S \rightarrow Y$
define the same point in $\te(S)$ if and only if $f\circ g^{-1} : Y \rightarrow X$
is isotopic to a conformal map. When $\partial S$ is nonempty, consider 
hyperbolic Riemann surfaces homeomorphic to $S$ with geodesic
boundary components of fixed length. Let $A=\partial S $ and $\eLL=(\eLL_{\alpha})_{\alpha
  \in A} \in {\Bbb R}_{+}^{|A|}$. A point $ X \in \te(S,\eLL)$ is a marked hyperbolic surface   
with geodesic boundary components such
that for each boundary component $\beta \in \partial S$, we have 
 $$\ell_{\beta}(X)=\eLL_{\beta}.$$ 
 Let $S_{g,n}$ be an oriented connected surface of genus $g$ with 
$n$ boundary components $(\CC_{1},\ldots, \CC_{n})$. Then 
$$\te(S_{g,n})(\eLL_{1},\ldots,\eLL_{n})=\te(S_{g,n},\eLL_{1},\ldots,\eLL_{n}),$$
 denote the Teichm\"uller space of hyperbolic structures on $S_{g,n}$
with geodesic boundary components of length $\eLL_{1},\ldots,\eLL_{n}.$
By convention, a geodesic of length zero is a cusp and we have  
$$\te(S_{g,n})=\te(S_{g,n})(0,\ldots,0).$$
Let $\operatorname{Mod}_{g,n}$
denote  the mapping class group of $S_{g,n}$, or the
group of isotopy classes of orientation preserving self
homeomorphisms of $S$ leaving each boundary component set wise fixed. The mapping
class 
group 
$\operatorname{Mod}_{g,n}=\operatorname{Mod}(S_{g,n})$ acts 
on $\te(S_{g,n})$ by changing the marking. The quotient space
$$\mathcal{M}_{g,n}=\mathcal{M}(S_{g,n})=
\te(S_{g,n})/ \operatorname{Mod}_{g,n}$$
is the moduli space of Riemann surfaces homeomorphic to $S_{g,n}$ 
with $n$ punctures. Given $\Gamma \subset \Mod_{g,n},$ define
$$\mathcal{M}_{g,n}[\Gamma]: = \te(S_{g,n})/ \Gamma.$$
Then $\mathcal{M}_{g,n}[\Gamma]$ is a finite cover of $\mathcal{M}_{g,n}.$

\subsection{The Weil-Petersson symplectic form}\label{G-S}
Recall that a {\it symplectic structure} on a manifold $M$ 
is a non-degenerate closed 2-form 
$\omega \in \Omega^{2}(M)$. The $n$-fold wedge product 
$$\omega \wedge \cdots \wedge \omega$$
never vanishes and defines a volume form on $M$.\\
By work of Goldman \cite{G1},
the space $\te_{g,n}(\eLL_{1},\ldots,\eLL_{n})$ carries a natural symplectic
form invariant under the action of the mapping class group.
This symplectic form is called
{\it Weil-Petersson symplectic form,} and denoted by $w$ or $w_{wp}$.
Note that when $S$ is disconnected, we have 
$$\vol_{wp}(\mathcal{M}(S,\eLL))=\prod\limits_{i=1}^{k} \vol_{wp}(\mathcal{M}(S_{i},\eLL_{A_{i}})).$$

\subsection{Thurston distance function}\label{Th:dis} 
Given $X, Y \in \te_{g,n}$, define 
$$d_{1}(X,Y)=\sup_{\lambda \in \ML_{g,n}} (\operatorname{log}(\frac{\ell_{\lambda}(X)}{\ell_{\lambda}(Y)})).$$
For more details see \cite{Thurston:minimal}. Since $d_{1}(X,Y)$ is not symmetric, 
we will work with
$$d_{Th}(X,Y)=\max\{d_{1}(X,Y),d_{1}(Y,X)\}.$$
Now let ${\bf B}_{X}(\epsilon) \subset \te_{g,n}$ be the ball of radius $\epsilon$ around $X\in \te_{g,n}$,
$${\bf B}_X(\epsilon)=\{Y|\; d_{Th}(X,Y)\leq \epsilon/2 \}.$$
Note that as $\epsilon \rightarrow 0,$ $\e^{\epsilon} \sim 1+\epsilon.$ Also, if 
$\epsilon$ is small then $1-2 \epsilon \leq \e^{-\epsilon}$, and $\e^{\epsilon} \leq 1+2 \epsilon.$
Hence if $Y \in {\bf B}_X(\epsilon),$ then for any $\lambda \in \ML_{g,n}$ we have
$$\ell_{\lambda}(X) (1-\epsilon) \leq  \ell_{\lambda}(Y) \leq \ell_{\lambda}(X) (1+\epsilon). $$

Note that since $\dim(\te_{g,n})=6g-6+2n$,  $\dim(\ell_{\gamma}^{-1}(L))=6g-7+2n$. One can estimate the volume of the ball of radius $\epsilon,$ as $\epsilon \rightarrow 0$ \cite{M:ABC}:   

\begin{prop}\label{bbbb}
For any $X \in K$, we have
\begin{itemize}
\item {\bf a)}: For small $\epsilon>0$,  $\vol_{wp}({\bf B}_X(\epsilon)) \asymp \epsilon^{6g-6+2n}.$ 
\item {\bf b)}: For any $\gamma \in \ML_{g,n}({\Bbb Z}),$
$$\vol_{wp}( {\bf B}_X(\epsilon) \cap \ell_{\gamma}^{-1}(L)) =O(\frac{\epsilon^{6g-7+2n}}{L}),$$
where the volume of ${\bf B}_X(\epsilon) \cap \ell_{\gamma}^{-1}(L) \subset H_{\gamma}(L)$ is with respect to the measure $\mu_{\gamma}^{L}.$
\end{itemize}
\end{prop}

\subsection {Marked pants decomposition and twisting information} 
Let $\mathcal{P}$ be a pants decomposition, decomposing $S_{g,n}$ to pairs of pants $\Sigma_1,\ldots, \Sigma_{2g-2+n}$. We will consider a more general form of this statement by considering {\it marked} pants decomposition $\widetilde{\mathcal{P}}$: for each $\alpha_{i}= \Sigma_{k} \cap \Sigma_{l}$ in  we have a pair $\eta_{1}, \eta_{2}$ of simple arcs joining $\alpha_{i}$ to two other curves in $\Sigma_{k}$ and $\Sigma_{l}$. This gives us a way of defining $0 \leq \tau_{\alpha_{i}} (X)\leq \ell_{\alpha_{i}}(X)$, {\it twist} around $\alpha_{i}$ on $X.$ In case $\alpha_{i}$ is in only one pair of pants, we can choose $\eta_1=\eta_2$ to be a simple arc joining $\alpha_i$ to itself.
Given $\alpha_{i}$, let $\Sigma_1$ and $\Sigma_2$ be the two pairs of pants containing $\alpha_i$: choose an arc $\delta_{i}$ intersection $\alpha_i$ at one point transversely joining two points on the boundary of $\Sigma_1 \cup \Sigma_2$. This naturally defines a closed curve $\beta_{i}$ disjoint from all other closed curves in $\mathcal{P}$. If $\Sigma_1=\Sigma_2$, $i(\beta_i, \alpha_i)=1$, and if  $\Sigma_1\not = \Sigma_2$, $i(\beta_i, \alpha_i)=2$. 
%SEE PICTURE.

\subsection{The Fenchel-Nielsen coordinates} \label{FNB}
A {\it pants decomposition} 
of $S$ is a set of disjoint simple closed curves which decompose 
the surface into pairs of pants. Fix a system of 
pants decomposition of $S_{g,n}$, 
$\mathcal{P}=\{\alpha_{i}\}_{i=1}^{k}$, where
$k=3g-3+n$. For a marked hyperbolic surface 
$X \in \te_{g,n} $, the {\it Fenchel-Nielsen coordinates}  
associated with 
$\mathcal{P}$, $\{\ell_{\alpha_{1}}(X),\ldots,\ell_{\alpha_{k}}(X),
\tau_{\alpha_{1}}(X),\ldots,
\tau_{\alpha_{k}}(X)\}$, 
 consists of the 
set of lengths of 
all geodesics used in the decomposition and the set of the 
{\it twisting} parameters used to glue the pieces.

We have an isomorphism 
\begin{equation}\label{FN}
 \operatorname{FN}_{\widetilde{\mathcal{P}}} : \te_{g,n} \cong {\Bbb R}_{+}^{3g-3+n} \times {\Bbb R}^{3g-3+n} 
 \end{equation}
by the map
$$X \rightarrow (\ell_{\alpha_{i}}(X),\tau_{\alpha_{i}}(X)).$$
By work of Wolpert, over Teichm\"uller space the
Weil-Petersson symplectic structure has a simple form in Fenchel-Nielsen
coordinates \cite{Wolpert:FN}. 

\begin{theo}[Wolpert]\label{wolp}
The Weil-Petersson symplectic form is given by
$$\omega_{wp} =\sum\limits_{i=1}^{k} d\ell_{\alpha_{i}}\wedge
d\tau_{\alpha_{i}}.$$
\end{theo}

\subsection{Integrating geometric functions over moduli spaces.}\label{integrate:geometric}
Here, we discuss a method for integrating certain geometric functions over $\mathcal{M}_{g,n}$ developed in \cite{M:In}.
Let $\mathcal{S}_{g,n}$ denote the set of homotopy classes of non-trivial, non-peripheral, 
simple closed curves on the surface $S_{g,n}$.
Let $\Lambda=(\gamma_{1},\ldots,\gamma_{k})$, where $\gamma_i$'s are distinct and disjoint elements of $\mathcal{S}_{g,n}.$
To each
$\Lambda$, we associate the set 
$$\mathcal{O}_{\Lambda}=\{ (\g\cdot \gamma_{1},\ldots, \g \cdot \gamma_{k}) \; | \; \g \in \Mod_{g,n}\}.$$ 
 Given a function $F:{\mathbb R}_{+}^{k} \rightarrow {\mathbb R}_{+}$, 
define $$F^{\Lambda}: \mathcal{M}_{g,n} \rightarrow {\mathbb R}$$
by 
\begin{equation}\label{function}
F^{\Lambda}(X)=\sum\limits_{(\alpha_{1},\ldots, \alpha_{k}) \in \mathcal{O}_{\Lambda}} 
F(\ell_{\alpha_{1}}(X), \ldots, \ell_{\alpha_{k}}(X)).
\end{equation}

Let $S_{g,n}(\Lambda)$ be the result of cutting the surface $S_{g,n}$ 
along $\gamma_{1},\ldots, \gamma_k$; that is $S_{g,n}(\Lambda) \cong S_{g,n}-U_{\Lambda}$, where $U_{\Lambda}$ is 
an open neighborhood of $\gamma_{1}\cup\ldots \cup \gamma_{k}$ homeomorphic to $\cup_{i=1}^{k} \gamma_{i} \times (0,1)$.
Thus $S_{g,n}(\Lambda)$ is a possibly disconnected surface 
with $n+2k$ boundary components; each $\gamma_{i}$ gives rise to two boundary components $\gamma_{i}^{1}$ and $\gamma_{i}^{2}$ of  $S_{g,n}(\Lambda)$.
Given ${\bf x}=(x_{1},\ldots, x_{k})$ with $x_{i}\geq 0,$ we consider the moduli space
$\mathcal{M}(S_{g,n}(\Lambda), \ell_{\Lambda}={\bf x})$ of hyperbolic Riemann surfaces homeomorphic to $S_{g,n}(\Lambda)$ such that for $1 \leq i \leq k$, $\ell_{\gamma_{i}^{1}}=x_{i}$ and $\ell_{\gamma_{i}^{2}}=x_{i}$. Then we have (\cite{M:In}):
\begin{theo}\label{integrate}
For any $\Lambda=(\gamma_{1},\ldots,\gamma_{k})$, the integral of $F^{\Lambda}$ over $\mathcal{M}_{g,n}$ with 
respect to the Weil-Petersson volume form is given 
by 
$$\int \limits_{\mathcal{M}_{g,n}}\!\! 
F^{\Lambda}(X)\,
dX=2^{-M(\gamma)} \;\int\limits_{{\bf x}
 \in {\mathbb R}_{+}^{k}} F(x_{1},\ldots,x_{k})\; V_{g,n}(\Lambda,{\bf x})\; {\bf
x}\cdot d {\bf x},
$$
where 
$ {\bf x}\cdot d{\bf x}=x_{1} \cdots
x_{k} \cdot dx_{1}\wedge \cdots \wedge dx_{k}$,
and 
$$M(\gamma)=| \{i | \; \gamma_{i} \mbox{ \; \mbox{separates} off a one-handle from \;} S_{g,n} \}|.$$
\end{theo}
Here given ${\bf x}=(x_{1},\ldots,x_{k}) \in {\mathbb R}_{+}^{k}$, $V_{g,n}(\Lambda,{\bf x})$ is defined by 
 $$V_{g,n}(\Lambda,{\bf x})= \vol(\mathcal{M}(S_{g,n}(\Lambda), \ell_{\Lambda}={\bf x}).$$
That is, $\ell_{\gamma_{1}^{1}}=x_{1}, \ell_{\gamma_{1}^{2}}=x_{1},\ldots, \ell_{\gamma_{k}^{1}}=x_{k},\ell_{\gamma_{k}^{2}}=x_{k} .$
Also,
$$V_{g,n}(\Lambda,{\bf x})=\prod \limits_{i=1}^{s}
V_{g_{i},n_{i}}(\ell_{A_{i}}),$$
where 
\begin{equation}\label{disj}
S_{g,n}(\Lambda)=\bigcup_{i=1}^{s} S_{i}\;,
\end{equation}
$S_{i} \cong S_{g_{i},n_{i}},$ and 
$\; A_{i}= \partial
S_{i}.$ \\
\noindent
{\bf Remark.} Given a multi curve $\gamma=\sum_{i=1}^{k} c_{i} \gamma_{i}$, the symmetry group of $\gamma$, 
$\operatorname{Sym}(\gamma)$, is defined by 
$$\operatorname{Sym}(\gamma)= \operatorname{Stab}(\gamma)/\cap_{i=1}^{k} \operatorname{Stab}(\gamma_{i}). $$
When $F$ is a symmetric function, we can define 
 $$F_{\gamma}: \mathcal{M}_{g,n} \rightarrow {\mathbb R}$$
$$F_{\gamma}(X)= \sum_{\sum_{i=1}^k c_{i}\alpha_{i} \in \Mod_{g,n} \cdot \gamma} F(c_{1}\ell_{\alpha_{1}}(X),\ldots, c_{k} \ell_{\alpha_{k}}(X)).$$
Then it is easy to check that 
\begin{equation}\label{Gamma}
F^{\Lambda}(X)= \operatorname{Sym}(\gamma) \cdot F_{\gamma}(X),
\end{equation}
where $\Lambda= (c_{1}\gamma_{1},\ldots, c_{k}\gamma_{k}).$

The following result will allow us to calculate the integral in Theorem $\ref{integrate}$ \cite{M:JAMS}, \cite{M:In}:
\begin{theo}\label{volm}
The volume
$V_{g,n}(L_{1},\ldots,L_{n})=\vol_{wp}(\mathcal{M}_{g,n}(L))$ is a polynomial in
$L_{1},\ldots,L_{n}$; namely we have:

$$V_{g,n}(L)=\sum\limits_{\alpha \atop |\alpha| \leq 3g-3+n}
C_{\alpha}\cdot L^{\; 2 \alpha},$$ where
 $C_{\alpha}>0 $ lies in $\pi^{6g-6+2n-|2\, \alpha|} \cdot {\Bbb Q}$.
\end{theo}

 \subsection{Measured laminations}\label{sec:ML}

 The space of measured laminations $\ML_{g,n}$ is the completion of the space of weighted 
 simple closed curves on a hyperbolic surface of genus $g$ with $n$ boundary components. A geodesic measured lamination 
 $\lambda$ consists of a closed subset of  $X$ foliated by complete simple 
 geodesics and a measure on every arc $k$ transverse to $\lambda$.
   For understanding measured laminations, it is helpful to lift them to the universal cover 
 $D$ of $X.$
 A directed geodesic is
determined by a pair of points $(x_{1}, x_{2}) \in  S ^{\infty} \times S^{\infty} - \Delta$, where 
$\Delta$ is the diagonal $\{(x, x)\}$.
A geodesic without direction is a point on $J = (S^{\infty}\times S^{\infty}-\Delta) / {\Bbb Z}_{2}$, where ${\Bbb Z}_{2}$ acts by  interchanging coordinates. Given a measured  geodesic  lamination  $\lambda$,  the  preimage  of  
of  its  underlying  geodesic  lamination  $A\subset  D$ is decomposed  as a union of  
geodesics of  $D$. Unless $\lambda$ is a union of simple closed curves, the complementary
regions are a union of polygons. Therefore geodesic laminations on two homeomorphic hyperbolic surfaces can be compared by passing to the circle at infinity, and the notion of a measured lamination
only depends on the topology of the surface $X.$

The weak topology on measures induce the measure topology on the space of measured
laminations; in other words, this topology is induced by the weak topology on the space 
of measured on a given arc which is transverse to each lamination from an open subset of
$\ML_{g,n}.$

For two simple closed curves $\gamma_{1}, \gamma_{2}$ the intersection number $i(\gamma_{1},\gamma_{2})$ is the minimum number of points in which representatives of $\gamma_{1}$ and $\gamma_{2}$ must intersect. 
 In fact, the intersection pairing extends to a continuous map $$i : \ML_{g,n} \times \ML_{g,n} \rightarrow {\Bbb R}_{\geq 0}.$$  

The  space of measured laminations on $S_{g,n}$ is a piecewise integral linear manifold homeomorphic to 
${\Bbb R}^{6g-6+2n}$ and ${\Bbb R}_{+}$ acts naturally on $\ML_{g,n}.$ However, $\ML_{g,n}$ does not have a natural  differentiable structure. 

The integral points in $\ML_{g,n}$ are in one to one correspondence with integral multi curves 
on $S_{g,n}.$ As a result $\ML_{g,n}$ carries a mapping class group invariant volume form 
$\vol_{Th}$ satisfying the following properties:
\begin{itemize}
\item For any open subset 
$U \subset \ML_{g,n}$ we have  
$$\vol_{Th}(t \cdot U)= t^{6g-6+2n} \vol_{Th}(U).$$

\item For any bounded set $U \subset \ML_{g,n}$, 
$\#|\{  \alpha \in U\;\mbox{ is an integral multi curve\;} \}| < \infty.$
Moreover, if $U$ is convex (in a train track chart) then as $L \rightarrow \infty$
\begin{equation}\label{limit}
\frac{\#|\{ \alpha \;| \;\alpha \in L \cdot U\;\mbox{ is an integral multi curve\;} \}|}{L^{6g-6+2n}} \rightarrow \vol_{{Th}}(U).
\end{equation}
\end{itemize}

See \cite{Thurston:book:GTTM} and \cite{Harer:Penner:Book} for more on the space of measured laminations.
 \subsection{Length function and right earthquake flow}\label{twist:define}
For any simple closed geodesic $\alpha$ on $X\in \te_{g,n}(\eLL)$ and $t \in {\Bbb R}$,
we can deform the hyperbolic structure as follows. We cut the surface
along $\alpha$, twist $\alpha$ distance $t$ to the right, and reglue back. Let us denote the 
new surface by $\operatorname{tw}^{t}_{\alpha}(X)$. As $t$ varies, the resulting continuous path 
in Teichm\"uller space is the Fenchel-Nielsen deformation of $X$ along
$\alpha$. For $t=\ell_{\alpha}(X)$, we have
\begin{equation}\label{D:T}
\tw^{t}_{\alpha}(X)=h_{\alpha}(X),
\end{equation}
where $h_{\alpha} \in \operatorname{Mod}(S_{g,n})$ is the right {\it Dehn
  twist} about $\alpha$.
  We remark that the notion of right and left only depend on the orientation of the underlying surface.\\
  
 By Wolpert's result (Theorem $\ref {wolp}$), the vector field generated by twisting around
$\alpha$
is symplectically dual to the exact one form $d \ell_{\alpha}$.
Recall that both the length function and the twisting deformation
are extended by homogeneity and continuity on $\ML_{g,n}$ \cite{Kerckhoff:eq:analytic}, \cite{Thurston:book:GTTM}, \cite{Thurston:earthquakes}.
More precisely, there is a continuous map   
$$\mathcal{L}: \ML_{g,n} \times \te  \rightarrow {\Bbb R}_{+},$$ defined by 
$$\mathcal{L}(\lambda,X)=\ell_{\lambda}(X)$$
where $\ell_{\lambda}(X)$ is the length of the lamination $\lambda$ on $X$. 

 Similarly, the right {\it earthquake} deformation is the extension of the right twist deformation for a general lamination such that 
 $$\operatorname{tw}_{r\cdot \eta}^{t}(X)=\operatorname{tw}_{\eta}^{r\cdot t}(X). $$ In other words, $\operatorname{tw}^{t}_{\lambda}$ is the limit of time $t$ twist deformation of 
 any sequence $\{r_{i} \gamma_{i} \}$ converging to $\lambda$ in $\ML_{g,n}.$

Any smooth function $H$ on a symplectic manifold $(M, \omega)$ gives rise to a vector field $X_{H}$ satisfying
$$\omega(X_{H},Y)=dH(Y).$$
It is easy to check that the {\it Hamiltonian} flow of the function $H$ generated by the vector field 
$X_{H}$ preserves both $H$ and $\omega$.
See \cite{Kerckhoff:eq:analytic} for more details.
\begin{theo}[Wolpert]\label{ham}
The earthquake flow along $\lambda$, $\operatorname{tw}^{t}_{\lambda}$ is the  the {\it Hamiltonian} flow of the length function $\ell_{\lambda}.$  Therefore, for $\lambda \in \ML_{g,n}$
the flow $\operatorname{tw}_{\lambda}$ is volume preserving.
 \end{theo}
 Let ${\bf E}_{\lambda}(X) \in T_{X} \te_{g,n}$ denote the tangent vector to the right earthquake path along $\lambda$
 at $X \in \te_{g,n}.$ By Theorem $\ref{ham}$, $\omega({\bf E}_{\alpha}, \cdot)= d\ell_{\alpha}.$ Also, the following result holds \cite{Kerckhoff:eq:analytic}: 
 \begin{theo}
 Let $X \in \te_{g,n}$. Then the  map ${\bf E}: \ML_{g,n} \rightarrow T\te_{g,n}(X)$ defined by
 $$\lambda \rightarrow {\bf E}_{\lambda}(X)$$
 is a homeomorphism.
 \end{theo}
 In other words, every tangent vector in $\te_{g,n}$ is tangent to a unique right earthquake path in the Teichm\"uller space.

Also,  for any $\lambda, \eta \in \ML_{g,n}$ 
 \begin{equation}\label{rcosine}
 | \omega({\bf E}_{\lambda}(X),{\bf E}_{\eta}(X))| < i(\lambda,\eta).
 \end{equation}

 \begin{theo}[Kerckhoff]\label{ker}
 For any $X \in \te_{g,n}$, the function $f(t)=\ell_{\eta}(\operatorname{tw}^{t}_{\lambda}(X))$ is convex. Moreover, unless $i(\lambda,\eta)=0$, $f$ is strictly convex and $\lim\limits_{t \rightarrow \pm \infty} f(t)=\infty.$
 \end{theo}
 
 \subsection{Bundle of measured geodesic laminations.}
Let $\mathcal{P}\te_{g,n}= \mathcal{M L}_{g,n} \times \te_{g,n}$ be the
  bundle of geodesic measured laminations over $\te_{g,n}$.
Let $\mathcal{P}^{1}\te_{g,n}$, the unit sub-bundle for the norm 
 $$\parallel(\lambda,X)\parallel= \ell_{\lambda}(X),$$
  and finally,  $\mathcal{P}^{1}\mathcal{M}_{g,n}[\Gamma]=\mathcal{P}^{1}\te_{g,n}/\Gamma$. 
Then Thurston's earthquake flow on $\mathcal{P}\te_{g,n}$
  is defined at time $t$ by 
$$\operatorname{tw}^{t}(X,\lambda)=(\operatorname{tw}_{\lambda}^{t}(X),\lambda).$$

 \subsection{Ergodic properties of the earthquake flow}\label{EPEF}
 First, we sketch the construction of an invariant measure $\nu_{g,n}[\Gamma]$ for the earthquake flow on $\mathcal{P}^{1}\mathcal{M}_{g,n}[\Gamma]$ in the Lebesgue measure class. 
Consider the Thurston volume form $\vol_{Th}$ on $\ML_{g,n}$, and 
 the Weil-Peterson volume form $\mu_{wp}$ on the Teichm\"uller space.
 It is known that:
  \begin{itemize}
 \item  Fixing $\lambda$ and $t \in {\Bbb R}$, by Theorem $\ref{wolp}$ the map
 $$\operatorname{tw}^{t}_{\lambda}: \te_{g,n} \rightarrow \te_{g,n} $$
  is the Hamiltonian flow of the length function of $\lambda$. Therefore, this map 
 preserves the Weil-Petersson volume.
 \item The Thurston volume form is mapping class group invariant  \cite{Harer:Penner:Book}.
\end{itemize}
 Therefore, $\mu_{Th} \times \mu_{wp}$ gives rise to a mapping class group
 invariant measure on $\ML_{g,n} \times \te_{g,n}$ which is invariant under the earthquake
 flow. On the other  hand, the length function
 $$\mathcal{L}: \ML_{g,n} \times \te_{g,n} \rightarrow {\Bbb R}_{+},$$
 defined by $\mathcal{L}(\lambda,X)=\ell_{\lambda}(X)$
 is invariant under the earthquake flow and the action of the mapping class 
 group. Now we define $\nu_{g,n}[\Gamma]$ to be the measure induced by $\mu_{Th} \times \mu_{wp}$
 on $\mathcal{P}^{1}\mathcal{M}_{g,n}[\Gamma]$.
 In other words,
 any point $X \in \te_{g,n}$ defines a measure $\nu_{X}$ on $\mathcal{P ML}_{g,n}$
such that for $A \subset \mathcal{PML}_{g,n}$, we have 
$$\nu_{X}(A)=\vol_{Th}(\{\lambda\;| \;[\lambda] \in A,\; \ell_{\lambda}(X) \leq
  1\}),$$
with respect to the Thurston volume form on $\ML_{g,n}.$ It is easy to check that $\{\nu_{X}\}_{X \in \te_{g,n}}$
is mapping class group invariant; namely,
for any $h \in \Mod_{g,n},$ $\nu_{h\cdot X}(h \cdot A)=\nu_{X}(A).$\\
The Weil-Petersson
 volume form on $\mathcal{M}_{g,n}[\Gamma]$ and the family $\{\nu_{X}\}_{X \in \te_{g,n}}$ of
 measures on $\mathcal{PML}_{g,n}$ combine to give the invariant
 measure $\nu_{g,n}[\Gamma].$
The measure $\nu_{g,n}[\Gamma]$ is a finite invariant measure for the earthquake flow on
$\mathcal{P}^{1}\mathcal{M}_{g,n}$ \cite{M:EE}.
In fact, by the definition, the measure $\nu_{g,n}$ projects to the volume
form given by $B(X) \cdot \mu_{wp}$ on $\mathcal{M}_{g,n}[\Gamma]$, where 
$$B(X)=\nu_{X}(\mathbb{P}\ML_{g,n})=\vol_{Th}(\{ \lambda\; |\; \ell_{\lambda}(X) \leq 1\}).$$ Thus we have
$$\nu_{g,n}[\Gamma](\mathcal{P}^{1}\mathcal{M}_{g,n}[\Gamma])=[\Mod_{g,n}:\Gamma] \times \int_{\mathcal{M}_{g,n}} B(X) \mu_{wp} < \infty .$$ 

In \cite{M:EE}, we establish a relationship between the earthquake flow and the Teichm\"uller horocycle
flow on the moduli space $\mathcal{Q}\mathcal{M}_{g,n}$ of 
meromorphic quadratic differentials with simple poles at $n$ poles, and obtain the following result:
\begin{theo}\label{er}
The earthquake flow on $\mathcal{P}^{1}\mathcal{M}_{g,n}[\Gamma]$ is ergodic 
with respect to the Lebesgue measure class.
\end{theo}

\noindent
{\bf Remark.} The main reason this result holds is that the Lebesgue measure is ergodic for the horocycle flow in the finite cover of $\mathcal{Q}^{1}\mathcal{M}_{g,n}.$
This is a corollary of the proof of this result due to Masur \cite{Masur:mapping}.

\subsection{Non Divergence of the Earthquake flow}\label{non-div}
Given $\epsilon >0,$ consider the set $K_{\epsilon}$ defined by 
 \begin{equation}\label{defmum}
 K_{\epsilon}=\{ X \;|\; \mbox{ for every simple closed geodesic}\; \gamma,  \ell_{\gamma}(X)> \epsilon\}.
  \end{equation} 
 By Mumford's criterion, $K \subset \mathcal{M}_{g,n}$ is compact if and only if
 for some $\epsilon,$ we have $K \subset K_{\epsilon}.$ 
  
Here we recall results of Minsky and Weiss \cite{Minsky:Wei} on non-divergence of the
earthquake flow on the bundle $\mathcal{P}^{1}\mathcal{M}_{g,n}$ of geodesic 
laminations. This result parallels important result in the homogeneous setting
that there are no divergent orbits for unipotent flows \cite{Dani:minimal}, \cite{Dani:orbits}.

Following \cite{Minsky:Wei}, given $y \in \te_{g,n},$ and $\lambda \in \ML_{g,n}$, define $\epsilon_{\gamma}(y,\lambda)$ and $t_{0}(y,\lambda)$, by 
$$\epsilon_{\gamma}(y, \lambda)=\operatorname{Min}_{t \in {\Bbb R}} \ell_{\gamma}(\operatorname{tw}_{\lambda}^{t}(y)),$$
$$\epsilon_{\gamma}(y,\lambda)=  \ell_{\gamma}( \operatorname{tw}_{\lambda}^{t_{0}(y,\lambda)}(y)).$$

Also, for any $\epsilon>0$ and compact subset $K$, there is a (larger) compact subset $K'$ such that for every $p \in K$ and every $T>0$ we have 
$$\operatorname{Avg}_{T,p}(K') \geq (1-\epsilon).$$ See Theorem ${\bf E}1$ of \cite{Minsky:Wei}.

Note that by Theorem \ref{ker}, when $i(\gamma,\lambda) \not =0$, $t_{0}$ is uniquely determined by the second equation. Also, let $J(\rho)=\{t\; |\; \ell_{\gamma}(\tw^{t}_{\gamma}(y)) \leq \rho \}.$
Note that $J(\rho)$ also depends on $y$, $\gamma$ and $\lambda.$
As in Lemma $5.2$ of \cite{Minsky:Wei}, one can control the growth of the length of $\gamma$ when $t$ is close to $t_{0}$ as follows.
 There are constants $\rho$ and $C$ such that for any $y \in \te_{g,n}$,
and $t \in J(\rho)$, we have 
$$ i(\gamma,\lambda)\; |t-t_{0}| -C \epsilon_{\gamma} \leq  \ell_{\gamma}(t).$$
As a result the following theorem holds \cite{Minsky:Wei}: 
\begin{theo}\label{WM}
For any $\epsilon,$ there exists a compact set $K \subset \mathcal{M}_{g,n}$ such that 
for any $(\lambda, X) \in \mathcal{P} \mathcal{M}_{g,n}$ exactly one of the following holds:
\begin{enumerate}
\item $$\liminf\limits_{T \rightarrow \infty} \frac{|\{t \in [0,T]\; | \;\operatorname{tw}^{t}_{\lambda}(X) \in K \} |}{T}> 1-\epsilon, $$
\item There is a simple closed curve $\gamma$ such that the function $f(t)=\ell_{\gamma}(\operatorname{tw}^{t}_{\lambda} (X))$ is equal to a constant smaller than $\epsilon.$
\end{enumerate}
\end{theo}

We remark that by Theorem $\ref{ker}$, the second case holds if and only if $i(\lambda,\gamma)=0$
and $\ell_{\gamma}(X) \leq \epsilon$.

Next, we apply this result for a connected multi curve $\gamma \in \ML_{g,n}({\Bbb Z})$. In this case, given $Y \in \te_{g,n}$ the path 
$\operatorname{tw}^{t}_{\gamma}(Y)$ projects to a closed curve in $\mathcal{M}_{g,n}.$
Thus we obtain:
\begin{coro}\label{mwle}
For any $\delta>0$, there exists $\epsilon>0$ such that for any $X \in \te_{g,n}$ and 
$\gamma \in \ML_{g,n}({\Bbb Z})$, exactly one of the followings hold:
\begin{enumerate}
\item There is a simple closed curve $\beta$, such that $i(\beta,\gamma)=0$, and 
$\ell_{\beta}(X)< \epsilon$
\item  $$\frac{|\{t \in [0,\ell_{\gamma}(X)]\; |\; \operatorname{tw}_{\gamma}^{t}(X) \in K_{\epsilon} \} |}{\ell_{\gamma}(X)}> 1-\delta, $$
where $K_{\epsilon}$ is defined in equation ($\ref{defmum}$).
\end{enumerate}
\end{coro}

\end{section}
\begin{section}{Example of once punctured torus and Markoff triples}\label{C-V}
In this section, we give a brief overview of the counting problem on a one-holed torus.
We recall that a Markoff triple is a solution $(p,q,r)$ of the Markoff equation
$$p^{2}+q^2+r^2= 3pqr $$
Let ${\bf M}_{{\Bbb N}}$ be the set of Markoff triples where $p,q,r \in {\Bbb N}$.
One can use the norm 
$$H(p,q,r)= p+q+r$$ or $$|(p,q,r)|= \max\{p,q,r\}$$
to measure the size of a Markoff triple.
Given a Markoff triple, one can obtain a new Markof triple by 
\begin{equation}\tag{***}
(a,b,c) \rightarrow (a',b,c),\;\;
(a,b,c) \rightarrow (a,b',c),\;\;
 (a,b,c) \rightarrow (a,b,c'),
 \end{equation}
where 
$a'=3bc-a,$ $b'=3ac-b,$ and $c'=3ab-c.$
Markoff showed that any Markoff triple in ${\bf M}_{{\Bbb N}}$ could be obtained from $(1,1,1)$ by applying the moves above. 

The growth of the set
$${\bf M}_{{\Bbb N}}(x)= \{(p,q,r) \in  {\bf M}_{{\Bbb N}} \; | \; |(p,q,r)|  \leq x\}$$
was first obtained by Zagier \cite{Z:M}: there exists $C>0$ such that as $x \rightarrow \infty$
\begin{equation}\label{Z}
\#({\bf M}_{{\Bbb N}}(x))=C \; (\log(x))^2+ O(\log(x) (\log(\log(x)))).
\end{equation}
The related problem of the growth of the set
$${\bf N}(x)= (0, x] \cap \bigcup_{(p,q,r) \in {\bf M}_{{\Bbb N}}} \{p,q,r\}$$ was also studied by McShane and Rivin \cite{MR}. See also \cite{ES}. \\

\noindent
{\bf Simple length spectrum of hyperbolic once punctured tori.}
It is known that the solutions of the Markoff equation are closely related to the geodesic  lengths of simple closed curves on hyperbolic once punctured tori \cite{Cohn:Markoff}. The main idea behind this correspondence is using 
trace identities in $SL_{2}({\Bbb R})$ discussed in \ch $\ref{trace:id}$

As before, let $S_{1,1}$ be a compact oriented surface of genus one with one boundary component. 
Let $(\alpha, \beta, \gamma)$ be a triple of simple closed curves on  $S_{1,1}$ such that $i(\alpha, \beta)=1$ and $\gamma=h_{\alpha}(\beta).$
Then there exists a hyperbolic once punctured torus $X_0 \in \te_{1,1}$ such that 
the map 
$$ \Mod_{1,1}\cdot (\alpha, \beta, \gamma) \rightarrow {\bf M}_{{\Bbb N}} $$
$$(\alpha', \beta', \gamma') \rightarrow (\cosh(\ell_{\alpha'}(X_0)/2), \cosh(\ell_{\beta'}(X_0)/2), \cosh(\ell_{\gamma'}(X_0)/2)),$$
is a one to one correspondence. 
The surface $X_0$ is the unique once-punctured torus which is produced by a the commutator of the modular group. This surface has the longest systole and even has maximal lengths in $\te_{1,1}$ \cite{Sc}.

In fact, Markoff moves correspond to the action of Dehn-twists on the pairs $(\alpha, \beta, \gamma).$ The acton of  $SL(2,{\Bbb Z})$  on ${\bf M}_{{\Bbb N}}$ is so that the 
\begin{displaymath}
 \begin{pmatrix} 1 & 0 \\ 1 & 1 \end{pmatrix}, 
  \begin{pmatrix} 1 & 0 \\ -1 & 1 \end{pmatrix}, 
  \begin{pmatrix} 0 & 1 \\ 0 & -1
  \end{pmatrix} 
\end{displaymath}
correspond to the three Markoff moves sending $(x, y, z)$ to $(x,-z-xy,y)$, $(z,y,-x-yz)$ and $(x,y,-z-xy).$\\
 Following the discussion in Remark {\bf 1} in \ch \ref{N:R}, by applying $(\ref{main:theo})$ in this case, we get:
  \begin{coro} 
  Let $\Gamma_0$ be a finite index subgroup of $SL_{2}({\Bbb Z}).$ Then we have
  \begin{equation}\label{Z:G}
 \# (\Gamma_0 \cdot (1,1,1) \cap {\bf M}_{{\Bbb N}}(x)) \sim  B \; \cdot (\log(x))^2 \times [\Mod_{1,1}: \Gamma_0],
 \end{equation}
 
as $x \rightarrow \infty.$ Here $B$ is a constant independent of $\Gamma_0.$
 \end{coro}

\noindent 
{\bf Character Variety of $S_{1,1}$.}
One can consider the variety $V(\kappa)$ defined by the set of triples $(x,y,z)$ satisfying
\begin{equation}\label{ME}
\kappa = x^2 + y^2 + z^2 - xyz-2
\end{equation}
in ${\Bbb R}^3.$ 
As before, the mapping class group $\Mod_{1,1}$ acts on this set by Markoff moves given by $(***)$.
Note that for equation ($\ref{ME}$) parametrizes the character variety of one holed torus with a boundary component. 
By work of Goldman \cite{G4} this action on each connected component of $V(\kappa)$ defined by $(\ref{ME})$ is either properly continuous or ergodic. 

 For $\kappa <-2 ,$ the variety $V(\kappa)$ parametrizes the moduli spaces of hyperbolic surfaces of genus one with one geodesic boundary component: the representations corresponding to hyperbolic structures give rise to contractible connected components. On the other hand, representations which map a simple closed curve to an elliptic element define an open subsets of on which $\Mod_{1,1}$ acts ergodically. 
See the appendix in \cite{G4} and references within for more details. 

We remark that he main counting result in this paper (Theorem \ref{morege}) holds for the solutions of $(\ref{ME})$ corresponding to hyperbolic surfaces of genus one with one geodesic boundary component.\\

\noindent 
{\bf Remark.}
Similarly, the equation 
$$x^2 +y^2 +z^2 +xyz = Ax+By+Cz+D,$$ where 
$A = ab+cd$, $B = bc+ad$, $C = ac+bd$, and $D = 4-a^2 -b^2 -c^2 -d^2 -abcd$
correspond to the character variety of four holed spheres. See \cite{G6} for more details. 

In fact, one can give a parametrization of other Teichm\"uller spaces in terms of lengths of simple closed curves by degree two equations. 

Recall that by work of Nielsen, the mapping class group is isomorphic to the outer automorphism group $\operatorname{Out}(\pi_{1}(S_{g,n}))$ and acts on the space of equivalence classes of representations $\pi_{1}(S_{g,n} )\rightarrow G$. Even in case $G=\operatorname{SL}(2, {\Bbb C})$ the dynamics of this action is very rich and interesting.

\end{section}

\begin{section} {Piecewise linear behavior of the length of a closed curve on $\te_{g,n}$}\label{PLL}

Let $\widehat{\mathcal{P}}$ be a marked pants decomposition of $S_{g}$. Consider the Fenchel-Nielsen coordinates on $\te_{g,n}$
(discussed in \ch \ref{FNB}). 

The main result of this section is the following statement:

\begin{theo}\label{LAPL}
 Let $\gamma$ be a closed curve on $S_{g,n}.$ The length function 
$$\ell_{\gamma}: \te_{g,n} \rightarrow {\Bbb R}_{+}$$
$$X \rightarrow \ell_{\gamma}(X)$$
is an asymptotically piecewise linear function of rational type with respect to the Fenchel-Nielsen coordinates. 
 
\end{theo}
Roughly speaking, the statement implies that far away from finitely many rational hyperplanes, $\ell_{\gamma}$ {\it asymptotically} behaves like a linear function with rational coefficients. However, we remark that the set of problematic hyperplanes would depend on $\gamma$: we expect this set to grow if we consider $\gamma$ with more self-intersections. \\

\noindent
{\bf Example.}
In the case of once punctured torus, let $\mathcal{P}=\{\alpha\}$ and $\gamma$ be a simple closed curve $i(\gamma,\alpha)=1.$ Then for any $X \in \te_{1,1}$ we have
\begin{equation}\label{1:1}
\cosh(\ell_{\gamma}(X)/2)= \cosh(\tau_{\alpha}(X)/2) \; \sinh(\ell_{\alpha}(X)/2)^{-1} \times \left(\frac{1+\cosh(\ell_{\alpha}(X))}{2}\right)^{1/2}.
\end{equation}
See the discussion in \ch $\ref{ex:ob}.$

The proof of Theorem $\ref{LAPL}$ is straightforward.

First, we will define the notion of 
asymptotically piecewise linear and give a few examples. 

\subsection{Main Definitions.} Let $\mathcal{C}$ be a closed cone in ${\Bbb R}^m$. We say $F: \mathcal{C} \rightarrow {\Bbb R}$ is {\it asymptotically linear} with respect to coordinates $x_{1},\ldots,x_{m}$
iff  there are linear functions $\mathcal{R}_{1},É.,\mathcal{R}_{m'}$ and $\LL: {\Bbb R}^m \rightarrow {\Bbb R}$ such that 
$$F(x_, \ldots, x_{m})- \LL(x_1,\ldots,x_m) \rightarrow c$$ uniformly as 
$\min\{\mathcal{R}_{i}(x_1, \ldots, x_m)\}_{i=1}^{m'} \rightarrow \infty$, where $c \in {\Bbb R}.$

Note that there is no restriction on $F$ {\it close} to the hyperplanes defined by $\mathcal{R}_{1}=0 ,\ldots, \mathcal{R}_{m'}=0.$

We say $F: \mathcal{C} \rightarrow {\Bbb R}$ is {\it asymptotically unbounded} with respect to coordinates $x_{1},\ldots,x_{m}$ iff 
iff  there are linear functions $\mathcal{R}_{1},É.,\mathcal{R}_{m'}$ and $L$ in $x_{1},\ldots, x_{m}$ such that 
$$F(x_, \ldots, x_{m}) \rightarrow \infty $$ uniformly as 
$\min\{\mathcal{R}_{i}(x_1, \ldots, x_m)\}_{i=1}^{m'} \rightarrow \infty$. Note that in general the linear function $L$ could be the zero function on some cones and an asymptotically piecewise linear does not have to be unbounded.  

We say $F$ is {\it asymptotically piecewise linear} iff there are linear functions $\mathcal{W}_{1},\ldots, \mathcal{W}_{k}$ such that for any ${\bf \epsilon}= (\epsilon_1,\ldots, \epsilon_{k}),$ $\epsilon_i=1$ or $-1,$ the restriction of $F$ on each sub cone defined by $\mathcal{C}_{\bf \epsilon}=\{ {\bf x}\; | \; \operatorname{Sign}(\mathcal{W}_{i}({\bf x}))=\epsilon_{i}\}$ is asymptotically linear.
We say  $F$ is of {\it rational type} if $\mathcal{R}_{1},É.,\mathcal{R}_{m'},$ $\mathcal{W}_{1},\ldots, \mathcal{W}_{k}$ and $L$ can be chosen to all have rational coefficients.

We say $F$ is {\it Strongly} asymptotically piecewise linear  if $\log(\sinh(F(x_{1},\ldots, x_{m}))$ is asymptotically piecewise linear. \

\subsection{Examples and basic observations} \label{ex:ob}
Before proving the main result,  we state some simple observations on the properties of asymptotically piecewise linear functions.\\
\noindent 
{\bf 1.} A simple example of an strongly asymptoticly linear function on ${\Bbb R}_{+}=\{x\;|\; x>0\}$ is $F_1(x)=\arccosh (\e^{x})$. This is simply because for $z>0,$ $0< \arccosh(z)-\log(z) <\log(2)$ and 
$\lim_{z\rightarrow \infty} \arccosh(z)-\log(z)= \log(2).$ It is easy to see that the convergence is uniform.  
On the other hand, as $z \rightarrow 0$, $\arccosh(1+z)$ behaves like $\sqrt{z}$: more precisely, we have 
$$\lim_{z\rightarrow 0} \frac{\log(\arccosh(1+z))-\log(z)/2-\log(2)/2}{z}=-1/12.$$
As a result $f(x)=\arccosh(1+\e^{x})$ is strongly asymptotically linear on ${\Bbb R}.$\\
Note that $\log(\frac{\e^x}{\e^y+\e^x})$ is asymptotically piecewise linear, but  not asymptotically unbounded as it is bounded on the cone $\{ (x,y)\; | \; x>y \}.$

 One can check that $F(x,y)= \arccosh(\e^{y}-\e^{x})$ is asymptotically linear on the cone
$\mathcal{C}=\{(x,y)\; |\; 0<x<y\}$: in this case, $L(x,y)=y,$ and $R(x,y)=y-x.$ Similarly  
$$H(x,y)= \arccosh\left(\frac{\e^{y}-(\e^{x}+ \e^{y})^{1/2}}{(1-2 \e^{-x})^{1/2}-(1-3 \e^{-y})^{1/3}}\right)$$ is asymptotically linear on any cone in ${\Bbb R}^{2}$ where it is defined.\\

\noindent 
{\bf 2.} 
 Let ${\Bbb R}^k=\{(x_{1},\ldots, x_{k})|\;  x_{i} \in {\Bbb R} \}$. Define  ${\bf e}_{i}: {\Bbb R}^{k} \rightarrow {\Bbb R}$ by ${\bf e}_{i}({\bf x})=\e^{x_{i}}.$
 Let  $\mathcal{F}_{k}$ to be the smallest family of functions ${\Bbb R}^{k} \rightarrow {\Bbb R}$ containing the functions ${\bf e}_{i}$ for $1 \leq i \leq k$ such that the following holds. For any two $f, g \in \mathcal{F}_k$ and $m\in {\Bbb N}$ 
we have $$\{f+g, f-g, f \times g, \frac{f}{g}, f^{1/m} \} \subset  \mathcal{F}_k.$$
Note that $\mathcal{F}_{k}$ can be naturally embed in $\mathcal{F}_{k+1}.$ Define 
\begin{equation}\label{FF}
\mathcal{F}= \bigcup_{k \geq 1} \mathcal{F}_{k}.
\end{equation}

 One can easily check that for any $P \in \mathcal{F}_{k}$
$$ F: {\Bbb R}^{k} \rightarrow {\Bbb R}$$
 \begin{equation} \label{general-type}
F({\bf x})= \arccosh(P ({\bf x}))
\end{equation}
  is asymptotically piecewise linear on any cone where it is defined. Let $F_{1},\ldots, F_{k}$ be functions in $\mathcal{F}_k$. Define
$$ G ({\bf x}) = F( F_{1}' ({\bf x}),\ldots, F_{k}'({\bf x})),$$
where $F_{i}' ({\bf x})= \arcsinh(F_{i} ({\bf x}))$ or $F_{i}' ({\bf x})= \arccosh(F_{i} ({\bf x}))$. 
Then $G \in \mathcal{F}_{k}; $ this is simply because $$ \e^{\arccosh(x)}= x+ \sqrt{x^2-1}\;\;, \mbox{and} \;\; \e^{\arcsinh(x)}=  x+ \sqrt{x^2+1}.$$ 
Hence, in view of $(\ref{general-type})$, $\arccosh(G({\bf x}))$ is also asymptotically piecewise linear on any cone where it is defined.  \\

More generally, let $\mathcal{A}$ be the family of functions $G: {\Bbb R}^{m} \rightarrow {\Bbb R}^{m}$ such that for each $1\leq i \leq m, $ we have $$\cosh(\pi_{i}(G)) \in \mathcal{F} \; \mbox{or} \; \sinh(\pi_{i}(G)) \in \mathcal{F}.$$ Here $\pi: {\Bbb R}^{m} \rightarrow {\Bbb R} \rightarrow {\Bbb R}$ is the natural projection map to the ith coordinate. 

Then it is easy to check that:
\begin{lemm}\label{composition}
We have 
\begin{itemize}
\item Any function $\mathcal{G} \in \mathcal{A}$ is asymptotically piecewise linear.
\item the coefficients of the linear functions approximating $G \in \mathcal{A}$ are in ${\Bbb Q}.$
\item The composition $G_1\circ G_2$ of any two maps in $\mathcal{A}$ is again in $\mathcal{A}.$
\end{itemize}
\end{lemm}
 \noindent 
{\bf 3.} One can check that if $G_{1} ({\bf x})$ and $G_{2}({\bf x})$ are strongly asymptotically piecewise linear in $\mathcal{C}$ then $G(x)=\log( 1+ \sinh(G_{1}({\bf x}) \times \sinh(G_{1}({\bf x}))$ is also strongly APL. \\

Let $G_{1}, G_{2}, G_{3}$ be strongly APL and let $\varepsilon \in \{+1,-1\}.$
Then 
$$G({\bf x})= \arccosh(2 \cosh(G_{1}({\bf x})) \cdot \cosh(G_{2}({\bf x}))+ \varepsilon \times \cosh(G_3({\bf x}))$$
is also asymptotically piecewise linear on any cone where it is defined.\\

 \subsection{Outline of proof of Theorem \ref{LAPL}.}\label{outline}
One can prove this result using trace identities or explicit calculations of $\ell_{\gamma}$ in terms of Fenchel-Nielsen coordinates \cite{O:explicit}. However, here we prove the statement in 2 steps. Some of the results obtained in the process might be of independent interest.  \\

\noindent
{\bf I.} First we prove the result for the case when $\gamma$ does not have self-intersections (see Proposition \ref{simple-case}). In this case, $\gamma$ can be extended to a pants 
decomposition $\mathcal{P}'$ of $S_{g,n}.$
An {\it elementary move} on a pants decomposition is replacing only one curve $\alpha$ in $\mathcal{P}$ by another curve $\beta$ disjoint from all other curves in $\mathcal{P}$ which has the minimum number of intersections with $\alpha$: if there there are two pants containing $\alpha,$ this number is $2$ and otherwise, it is $1$. 
Following work of Hatcher and Thurston \cite{HT:moves}, given $\mathcal{P}$ and $\mathcal{P}'$, one can obtain $\mathcal{P}'$ from $\mathcal{P}$ by doing finitely many elementary moves. 
One corollary of our proof is that the change of coordinates on $\te_{g,n}$ from $\mathcal{P}$ to a different pants decomposition $\mathcal{P}'$
is asymptotically piecewise linear. 
We use work of Hatcher and Thurston to get this result by only doing the calculations in four holed spheres and one hold tori with geodesic boundary components. \\

\noindent
{\bf II.}
In view of some observations in \ch \ref{ex:ob}, the result for an arbitrary connected closed curve is obtained by induction on the number of self intersections and trace identities (see Lemma $\ref{trace}$).  \\

First we prove:

 \begin{prop}\label{simple-case}
 Let $\mathcal{P}$ and $\mathcal{P}'$ be two marked pants decompositions of $S_{g,n}$. The map 
 $$ {\Bbb R}_{+}^{3g-3+n} \times {\Bbb R}^{3g-3+n}  \rightarrow {\Bbb R}_{+}^{3g-3+n} \times {\Bbb R}^{3g-3+n}$$
 \begin{equation}\tag{*}
  ({\bf l}, {\bf t}) \rightarrow  \operatorname{FN}_{\mathcal{P}'}( \operatorname{FN}_{\mathcal{P}}^{-1}({\bf l}, {\bf t}))
  \end{equation}
 is asymptotically piecewise linear. 

  \end{prop}

\noindent
{\bf Remarks.}\\
\noindent
{\bf 1.} This can be proved by using the calculation in \cite{O:change} (see also \cite{ALPS1}): one can easy check that the length and twisting are both given in terms of functions of type $(\ref{general-type})$. However we give a different proof as the local calculation we do will be useful for us later.  \\
 
\noindent
 {\bf 2.} By results in \cite{ALPS1} the map in $(*)$ is not Lipschitz with respect to natural distance function on $\te_{g,n}$
defined by $\mathcal{P}$ and $\mathcal{P}'.$

 \subsection {Trigonometry of pairs of pants and right-angled hexagons.}\label{Tr-h}

Recall that  for each triple $(x,y,z)$, there is a unique pair of pants with these boundary lengths. 
One can obtain a hyperbolic pair of pants with geodesic boundary components by gluing identical right angled hexagons. 
Following \cite{B:book} (\ch $3$) let $H$ be a right-angled geodesic hexagons with consecutive sides $a, \tilde{c}, b, \tilde{a}, c, \tilde{b}$. Then, we have\\

\noindent
{\bf 1.} If $H$ is convex, we have
\begin{equation}\label{HC1}
  \cosh(c)= \frac{\cosh(\tilde{c})+ \cosh(\tilde{a}) \cosh(\tilde{b})}{\sinh(\tilde{a}) \sinh(\tilde{b})},
\end{equation}
and\;

\noindent
{\bf 2.} If $H$ has intersecting sides $c$ and $\tilde{c}$ we have
\begin{equation}\label{HI}
\cosh(\tilde{c})= \sinh(\tilde{a}) \sinh(\tilde{b}) \cosh(c)+ \cosh(\tilde{a}) \cosh(\tilde{b}).
\end{equation}

See \cite{B:book} for more details.

One can easily check that $\tilde{c} < \tilde{a}+ \tilde{b}+c.$ Note that this inequality is not sharp as $\tilde{c}$ can be much smaller than 
$ \tilde{a}+ \tilde{b}.$ 

However, it is important for us to obtain bounds on $\tilde{c}-c$ in terms of
$\tilde{a}$ and $\tilde{b}.$ Note that if $\tilde{c}$ and $c$ intersect, then 
$\tilde{a} \leq  \tilde{c}$ and $\tilde{b} \leq  \tilde{c}.$ The following elementary lemma will be used in the proof of 
Lemma \ref{local:main}:

\begin{lemm}\label{length:set}
There exists $c_0>0 $ such that for any right-angled hexagon $H$ with consecutive sides $a, \tilde{c}, b, \tilde{a}, c, \tilde{b},$  the followings hold. 
\begin{enumerate}
\item 
 If the sides $c$ and $\tilde{c}$ intersect (as in $(\ref{HI})$), 

$c > 4 \max\{-\log(\tilde{a})- \log(\tilde{b})+c_0, c_0\}$  then we have: 
\begin{equation}\label{basic:estimate1}
\tilde{c}= c+\log(\sinh(\tilde{a}))+\log(\sinh(\tilde{b})) + O(\e^{-c/10}),
\end{equation}

\item If the sides $c$ and $\tilde{c}$ intersect, $\frac{1}{c_0}> \max\{\tilde{a}, \tilde{b}\} ,$ and $\tilde{c} > 4 \max\{-\log(\tilde{a})- \log(\tilde{b})+c_0, c_0\},$
\begin{equation}\label{basic:estimate2}
\tilde{c}= c+\log(\sinh(\tilde{a}))+\log(\sinh(\tilde{b})) + O(\e^{-\tilde{c}/10}),
\end{equation}
\item If $H$ is convex (as in $(\ref{HC1})$),
$c> 2(\tilde{a}+ \tilde{b}+c_0)$, and $\min\{\tilde{a}, \tilde{b}\} > c_0$ then we have
\begin{equation}\label{basic:estimate3}
\tilde{c}= c+\log(\sinh(\tilde{a}))+\log(\sinh(\tilde{b})) + O(\e^{-\tilde{c}/2 }).
\end{equation}
\end{enumerate}
\end{lemm}

Note that $\log(\sinh(x)) \sim x$ as $x \rightarrow \infty$ and $\log(\sinh(x))\sim \log(x) $
as $x \rightarrow 0.$
 For the proof of Lemma \ref{length:set} see \ch \ref{App:1}.\\

We consider the function $F_1: {\Bbb R}^3 \rightarrow {\Bbb R}$ given by
\begin{equation}\label{def:F1}
F_1(x,y,z)= \arccosh(\frac{\cosh(z)+ \cosh(x) \cosh(y)}{\sinh(x) \sinh(y)}).
\end{equation}
Note that by $(\ref{HC1})$,  $F_{1}(x,y,z)$  the length of the shortest path between $\alpha$ and $\beta$ in a pair of pants with boundary lengths 
$\ell_{\alpha}=x,$ $\ell_{\beta}=y$ and $\ell_{\gamma}=z.$

In view of \ch \ref{ex:ob} $({\bf 5})$ this functions is strongly asymptotically PL in any cone inside $\{(x,y,z)|\; 0<x, 0<y, 0<z\}$ where they are defined.  We will also analyze the asymptotic properties of $F_1$ in $\ch \ref{App:1}$. 

Given $\epsilon>0
$, define  ${\bf B}_{\epsilon}(x,y,z) \subset {\Bbb R}_{+}^3$ to be the set of points $(x',y',z')$ such that 
$$\frac{1}{1+\epsilon} < \frac{x}{x'} < 1+\epsilon, \; \frac{1}{1+\epsilon} < \frac{y}{y'} < 1+\epsilon\; \; \mbox {and} \; \frac{1}{1+\epsilon} < \frac{z}{z'} < 1+\epsilon.$$

We will use the following statement in \ch \ref{local:calculation} to prove Lemma \ref{local:main}
which is essential in \ch \ref{from-eq-to-count}:

\begin{lemm}\label{basic:F1}
Let $\epsilon< 1.$ Then there exists $L_0>0$ such that  for any 
$(x',y',z') \in {\bf B}_{\epsilon}(x,y,z)$ with $\min\{x,y,z\} >L_0$ 
we have 
$$ |  \log(\sinh(F_{1}(x,y,z)))-  \log(\sinh(F_{1}(x',y',z')))| \leq 5 \epsilon \times \max\{x,y,z\}.$$
\end{lemm}
 These statements are straightforward and elementary. For completeness we will include the proof of Lemma $\ref{basic:F1}$ and Lemma $\ref{length:set}$ in \ch \ref{App:1}.

\subsection{Length functions in the case of one hold tori and four holed spheres.}\label{HT:moves}
As elementary moves discussed in \ch \ref{outline} (see \cite{HT:moves}) are supported on one-holed torus or four-holed spheres, here we will sketch the local calculation in these two cases.\\

Let $${\bf T}_{0,4} (\ell_1,\ell_2,\ell_3,\ell_4)=\{(\ell,\tau)|\; 0 \leq \ell,\;  \tau \in {\Bbb R}\}.$$
 Let $Y=Y(\ell_1,\ell_2,\ell_3,\ell_4,\ell,\tau)$ denote the a hyperbolic sphere with $4$ geodesic boundary components of lengths $\ell_1,\ldots, \ell_4$ obtained by gluing two pairs of pants with boundary lengths $\ell_1,\ell_2,\ell$ and $\ell_3,\ell_4,\ell$ with gluing $\tau$ along 
 the geodesic $\alpha$ of length $\ell.$  Let $\gamma$ be a closed curve intersecting $\alpha$ at 2 points.

Similarly, let ${\bf T}_{1,1}(\ell_1) =\{(\ell,\tau)|\; 0 \leq \ell, \ell_{1}, \tau \in {\Bbb R}\}$. Let $X=X(\ell_1,\ell,\tau)$ denote the a torus with one geodesic boundary components of lengths $\ell_1$ obtained by gluing two edges of a pair of pants with boundary lengths $\ell_1,\ell,\ell$ with gluing $\tau$ along 
 the geodesic of length $\ell.$ Let $\gamma$ be a closed curve intersecting $\alpha$ at one points.

Then one can easily check:

 \begin{lemm}\label{torus}
The maps 
 $$ L_{0,4}:   {\bf T}_{0,4}(\ell_1,\ell_2,\ell_3,\ell_4) \rightarrow {\Bbb R}$$
 $$(\ell,\tau) \rightarrow \ell_{\alpha}(Y(\ell_1,\ell_2,\ell_3,\ell_4,\ell,\tau))$$ 
and
 $$ \tau_{0,4} : {\bf T}_{0,4} (\ell_1,\ell_2,\ell_3,\ell_4), \rightarrow {\Bbb R}$$
 $$(\ell,\tau) \rightarrow \tau_{\alpha}(Y(\ell_1,\ell_2,\ell_3,\ell_4,\ell,\tau))$$ 
are asymptotically piecewise linear and unbounded.
\end{lemm}

\begin{lemm}\label{sphere}
  Then 
 $$  L_{1,1}:  {\bf T}_{1,1}(\ell_1) \rightarrow {\Bbb R}$$
 $$(\ell,\tau) \rightarrow \ell_{\alpha}(X(\ell_1,\ell, \tau))$$ 
and
 $$ \tau_{1,1} : {\bf T}_{1,1} (\ell_1) \rightarrow {\Bbb R}$$
 $$(\ell, \tau) \rightarrow \tau_{\alpha}(X(\ell_1,\ell,\tau))$$ 
are asymptotically piecewise linear and unbounded.
\end{lemm}

One can prove both lemmas by direct calculation from \cite{O:change} (see also \cite{ALPS1}). \\

\noindent
{\it Sketch of proof of Lemma $\ref{torus}$ and Lemma $\ref{sphere}.$} 
As in \ch $3.3$ in \cite{B:book}, one can use ($\ref{HI}$) to calculate the length $d_{\tau}$ : the arc perpendicular to $\alpha_{1}$ and $\alpha _{3}$ which is homotope to the path $a_{1}, $ followed by inversely parametrized $a_{2}$ after moving around $\alpha$ for length $\tau.$ The length of $a_{1}$ and $a_2$ can be calculated using $(\ref{HC1})$
The length of the closed curve $\delta$ can be calculated using ($\ref{HC1}$).
The geodesic length of $\delta$ can be calculated using
$$ \cosh(\ell_{\delta}/2)= \sinh(x/2) \sinh(y/2) \cosh(d_{\tau}) - \cosh(x/2) \cosh(y/2).  $$
See the proof of Proposition $3.3.12$ in Buser's book \cite{B:book}. The exact calculation is not important; here it is enough to note that following the discussion in \ch $\ref{ex:ob} ({\bf 2}),$ we have

\begin{equation}\label{subset}
 \{\cosh(L_{0,4}), \cosh(\tau_{0,4}), \cosh(L_{1,1}),  \cosh(\tau_{0,4})\} \subset \mathcal{F}.
 \end{equation}
Then the result follows from Lemma \ref{composition}.
We also need to show that $\tau_{1,1}$ and $\tau_{0,4}$ are asymptoticlly unbounded on any cone. This is a corollary of Wolpert's result (Theorem \ref{wolp}) for ${\bf T}_{1,1}$ and ${\bf T}_{0,4}$. In the case of four holed torus, this result states that 
$${\Bbb R}^2 \rightarrow {\Bbb R}^{2}$$
$$(\ell, \tau) \rightarrow (L_{0,4}(\ell_1,\ell_2,\ell_3,\ell_4, \ell, \tau), \tau_{0,4} ( (\ell_1,\ell_2,\ell_3,\ell_4,, \ell, \tau))$$
is volume preserving.  

\hfill $\Box$

\subsection{Proof of Theorem \ref{simple-case}}
By Hatcher-Thurston result, there exists a sequence $\mathcal{P}_0, \mathcal{P}_1, \ldots \mathcal{P}_k$ of pants decompositions such that $\mathcal{P}=\mathcal{P}_0$ 
and $\mathcal{P}'=\mathcal{P}_k$ and each $\mathcal{P}_i$ is obtained from $\mathcal{P}_{i-1}$ by an elementary move supported on a one-holed torus or a four-holed sphere. 

In view of Lemma \ref{composition}, Lemma $\ref{torus}$ and Lemma $\ref{sphere}$ imply that the map $(*)$ belongs to the family in $\mathcal{A}$. By Lemma \ref{composition} 
the map $(*)$ is asymptotically PL. 

\subsection{Proof of Theorem \ref{LAPL}.}
In order to prove Theorem \ref{LAPL}, it is enough to consider resolutions of $\gamma$ at a self intersection $p.$ Note that for each $p$, $\gamma$ has two resolutions: one of them is connected and the other one is disconnected. See \cite{G2}.

We can then apply the following lemma (see \cite{R:Y} and references within):

\begin{lemm} \label{trace}
Let $\gamma$ be a closed curve on $S_{g,n}$ with $i(\gamma,\gamma)>0$. There are closed curves $\gamma_{1},\gamma_2, \gamma_{3}$
such that $i(\gamma_j,\gamma_j) < i( \gamma, \gamma)$ for $1\leq j \leq 3$ and for any $X \in \te_{g,n}$ we have:
$$\ell_{\gamma}(X)= 2 \arccosh(2 \cosh(\ell_{\gamma_{1}}(X)/2) \cdot \cosh(\ell_{\gamma_2} (X)/2)+ \varepsilon\; \cosh(\ell_{\gamma_3}(X)/2)),$$
where 
$\varepsilon \in \{1,-1\}.$
\end{lemm}

One can obtain $\gamma_{1}, \gamma_2$ and $\gamma_3$ by resolving one self intersection of $\gamma$ and using the trace formula for matrices in 
$\operatorname{SL}(2, {\Bbb R}).$ However we have to be careful about signs.
If $\alpha$ has only one self intersection, $\gamma_1$, $\gamma_2$ and $\gamma_3$ form a pair of pants containing $\gamma$ (see also \cite{R2}). 
Then $\varepsilon $ is equal to $-1$ if $\gamma$ contains a geodesic triangle at $p$ (see \cite{R:Y} for more details). 

 \subsection{Remark on non-filling closed curves} \label{non-filling}
 Given $\gamma \in \pi_{1}(S_{g,n}),$ let $$\mathcal{S}(\gamma)= \{\alpha,\; \mbox{simple closed curve}\; | \; \operatorname{Stab}(\gamma)\cdot \alpha= \alpha, i(\alpha, \gamma)\not =0\}.$$
 We also consider subsurface (with boundary) ${\bf S}(\gamma)$ such that $\gamma$ is filling in ${\bf S} (\gamma).$ The subsurface ${\bf S}(\gamma)$ is the subsurface we get by considering a small neighborhood of $\gamma.$ 
 Then $\mathcal{S}(\gamma)$ consists of all the curves inside of ${\bf S}(\gamma).$
 
 If $\gamma$ is filling curve then every simple closed curve on $S_{g,n}$ is in $\mathcal{S}(\gamma).$ Also if $i(\gamma,\gamma) <2,$ we have $\mathcal{S}(\gamma) =\emptyset.$
Also, $\operatorname{Stab}(\gamma) \cdot \beta=\beta$ for every $\beta \in \mathcal{S}(\gamma).$

Note that up to finite index, the group $\operatorname{Stab}(\gamma)$ is generates by twisting around the curves not intersecting $\gamma.$ Let 
${\bf S}(\gamma)= \cup_{i \in {\bf I}} \Mod_{g_{i}, n_{i}} $ and $\partial({\bf S}(\gamma))= \cup _{j \in {\bf J}} {\bf T}_{\alpha_j}.$ Then
 $\operatorname{Stab}(\gamma)$ is virtually the same as 

$$\prod_{i \in {\bf I}} \Mod_{g_{i}, n_{i}} \times \prod_{j \in {\bf J}} {\bf T}_{\alpha_j}, $$ where ${\bf T}_{\alpha_j}$ is the group generated by Dehn twists around 
$\alpha_{j} \in \partial({\bf S}(\gamma)).$

Let 
$$\Gamma_{\gamma}= \Stab({\bf S}(\gamma)).$$

In general, $\ell_{\gamma}$ is a well-defined on $\te_{g,n}/ \Gamma_{\gamma},$ and $\{Z \in \te_{g,n}/ \Gamma_{\gamma}, \ell_{\gamma}(Z) \leq L\}$ has finite Weil-Petersson volume. 

To see this, we can choose $\mathcal{P}$ a pants decomposition of $S_{g,n}$ such that $\partial({\bf S}(\gamma)) \subset \mathcal{P}.$ Let 
$\mathcal{P}_{\gamma}= \mathcal{P} \cap \; \Stab(\gamma) \subset \mathcal{P}.$
Note that $\ell_{\gamma}(X)$ is independent of $\ell_{\alpha}$ or twisting parameter around for $\alpha$ which does not intersect the surface ${\bf S}(\gamma).$ 
So $\ell_{\gamma}$ only depends on the lengths and twists of curves in or on the boundary of ${\bf S}(\gamma).$ 
It is important for us that Theorem \ref{LAPL} also holds for non-filling closed curves and even on $\te_{g,n}/\Gamma_{\gamma},$ $\ell_{\gamma}$ is an asymptotically piecewise linear function in terms of $\{\ell_{\alpha}, \tau_{\alpha}\}_{\alpha \in \mathcal{P}_{\gamma}}.$  

\subsection{Applications of local calculations}\label{local:calculation}
The calculation in the previous subsection would also imply the following lemmas which will be used in \ch \ref{from-eq-to-count}. The main goal is to understand the shape of ball of radius $\epsilon$ with respect to the Thurston distance function (defined in \ch \ref{Th:dis}). In general, this can be complicated. However, we prove the bounds we need by analyzing the case of ${\bf T}_{1,1}$ and ${\bf T}_{0,4}.$

\begin{lemm}\label{local:main}

Let $M>1$ and $ 0<\epsilon< \frac{1}{2 \cdot M^3} .$ There exists $L_{0}>0$ depending on $\epsilon$ and $M$ such that the following statement holds.  Let $\alpha$ be a simple closed curve disjoint from all the components $\{\alpha_{i}\}$ of a four-holed sphere (resp. one holed torus). 
Assume $Y \in {\bf T}_{0,4}$ (resp. ${\bf T}_{0,1}$) is such that for any boundary curve $\alpha_i,$ we have
$$ \frac{1}{M} < \frac{\ell_{\alpha}(Y)}{\ell_{\alpha_i}(Y)}< M,$$
and $L=\max\{\ell_{\alpha_{i}}(Y)\}_{i} > L_0.$ Then we have :\\

\noindent
{\bf 1.} There exists $\beta$ with $i(\beta, \alpha)=2$  (resp. $=1$ for ${\bf T}_{0,1}$) such that 
$$\ell_{\beta}(Y) \leq 4L\; \mbox{and}\; \; \ell_{h_{\alpha}^{k}(\beta)}(Y) \geq (k-2) \ell_{\alpha}(Y),$$ for $k \geq M.$

\noindent
{\bf 2.}  For any $Y'$ with $d_{Th} (Y, Y')< \epsilon,$ we have 

 \begin{equation}\label{twist:change}
 |\tau_{\alpha}(Y)- \tau_{\alpha}(Y') | \leq \; C \cdot M^2 \cdot \epsilon \times L,
 \end{equation}
where $C$ is a constant independent of $M$, $\epsilon, $ $Y,$ and $Y'.$

\end{lemm}
Here fix a notion of twisting on $Y$ and measure $\tau_{\alpha}(Y')$ the same way. \\

We prove Lemma \ref{local:main} in \ch \ref{App:1} using Lemma \ref{length:set} and Lemma \ref{basic:F1}.

We will also need the following Lemma in \ch \ref{proof:morege}: 

\begin{lemm}\label{local:bound}
Let $\alpha$ be a simple closed curve disjoint from all the components $\{\alpha_{i}\}$ of a four-holed sphere (resp. one holed torus). 
Let  $\beta$ be a simple closed curve on $Y$ with $i(\beta, \alpha)=2$  (resp. $=1$ for ${\bf T}_{0,1}$) used for defining $\tau_{\alpha}(Y).$
Assume $Y \in {\bf T}_{0,4}$ (resp. ${\bf T}_{0,1}$) is such that for any boundary curve $\alpha_i,$ we have
$\ell_{\alpha_i}(Y) > L_0,$ and $\ell_{\alpha}(Y) >L_0.$
such that 
$$\ell_{\beta}(Y) \leq L\;$$

Then we have 
$$ \ell_{\beta}(Y)+ L  \geq |\tau_{\alpha}(Y)| \cdot \ell_{\alpha}(Y),$$
where $L= \max \{\ell_{\alpha_{i}} (Y) \}_{i} \cup \ell_{\alpha}(Y).$ 
\end{lemm}

\noindent
{\bf Proof.}
By $(\ref{HI})$
$d \geq \tau -2 L $ and by $(\ref{HC1})$ $\ell_{\beta} > d/2.$

\hfill $\Box$

As in the introduction, 
$\mathcal{B}_{\mathcal{P}} (M) \subset \te_{g,n}$ is the set of points such that 
 for all   $\alpha, \beta \in \mathcal{P}$ we have $\frac{\ell_{\alpha}(X)}{\ell_{\beta}(X)} < M.$ Then in view of Lemma $\ref{local:main}$, we have:  

\begin{coro}\label{coro:est}
Let $\mathcal{P}$ be a pants decomposition on $S_{g,n},$ and $\epsilon, M>1$ with $\epsilon < \frac{1}{2M^3}$. There exists $L_0>0$ such that the following holds. 
If $X, X_{1} \in \mathcal{B}_{\mathcal{P}}(M) \subset \te_{g,n}$ such that 
$d_{Th}(X,X_{1}) <\epsilon$ and $L_{\mathcal{P}}(X) > 2 M\cdot L_0$ then 
$$ |\tau_{\alpha_i}(X)- \tau_{\alpha_i}(X_{1}) | \leq C \; M^2  \epsilon \times L_{\mathcal{P}(X)}, $$
where $L_{\mathcal{P}}(X)=\max\{\ell_{\alpha_i}(X)\}_{i=1}^{3g-3+n},$ 
and $C$ is a constant independent of $M$, $\epsilon, $ $X,$ and $X_1.$
\end{coro}

This Corollary will be used later in \ch \ref{from-eq-to-count}. The important case for us is when $M$ is large but fixed and $\epsilon \ll 1/M.$

\begin{lemm}\label{twist:bound}
 Let $\gamma$ be a closed curve on $S_{g,n}.$ Then for any closed curve $\beta \subset {\bf S}(\gamma)$ and 
 $ Z \in \te_{g,n}$ we have
$$\ell_{\beta}(Z) \leq c_{\gamma,\beta} \; \ell_{\gamma}(Z).$$
Moreover,  let $\mathcal{P}=\{\alpha_1,\ldots, \alpha_{3g-3+n}\}$ be a pants decomposition on $S_{g,n}$  such that 
$\partial({\bf S})(\gamma) \subset \mathcal{P}$.
Then for any $\alpha \in \mathcal{P}$ with $i(\alpha, \gamma) \not =0$ we have 
$$ |\tau_{\alpha}(Z)| \leq c'_{\gamma,\mathcal{P}} \; \ell_{\gamma}(Z), $$
where the constant $c_{\gamma,\beta}$ depends on $\gamma$ and  $\beta$ and $c'_{\gamma,\mathcal{P}}$ only depends on $\gamma$ and $\mathcal{P}.$
\end{lemm}

Note that by the second part of this lemma, $h_{\alpha}^{m}(X)\in {\bf B}_{\gamma}(L)$ for large $m$ only if $\ell_{\alpha
}(X)$ is much less than $L.$

We remark that the first part of this proposition can be done using Lemma $5.1$ of \cite{Bas}:
\begin{lemm}\label{bound:Bas}
Let $\gamma$ be a filling closed curve on $S_{g,n}$ and $Z \in \te_{g,n}.$ Then for any simple closed curve $\beta$ we have
$$ \ell_{\beta}(Z) \leq \frac{1}{2} \ell_{\gamma}(Z) \times i(\gamma, \beta).$$
\end{lemm}
\noindent
{\bf Remark.}
In case $n=0$, one can consider the function $$\mathcal{H}_{\omega}(Y_1,Y_2 )= \frac{i(Y_1, Y_2)}{i(\omega,Y_1) \times i(\omega,Y_2)}$$ on $\mathbb{P}(\mathcal{G}(S_{g})) \times \mathbb{P}(\mathcal{G}(S_{g})).$ 
Then 
$$\mathcal{H}_{\omega}(Y_1, Y_2) \leq c(\omega).$$ 
Note that by work of Thurston and Bonahon \cite{Bon} the intersection pairing is continuous and the space of projectivized geodesic currents is compact. On the other hand, $i(\omega,Y)$ is never zero. Hence,  we have $$\mathcal{H}_{\omega}(Y_1, Y_2) \leq c(\omega).$$ 

We need the following weaker observation which holds for an arbitrary closed curves on $S_{g,n}$:
\begin{lemm}\label{bound:Bas:general}
Let $\gamma$ be a closed curve on $S_{g,n}$ and $\beta\in {\bf S}(\gamma).$ There exists $c_{\beta,\gamma}$
such that 
$$ \ell_{\beta}(X) \leq c_{\beta,\gamma} \times \ell_{\gamma}(X)$$
for $X \in \te_{g,n}.$
\end{lemm}

\noindent
{\bf Proof of Lemma \ref{twist:bound}.}
The first part of the lemma is a consequence of Lemma \ref{bound:Bas:general}. In order to prove the second part, we consider the curve $\beta$ defined by the marked pants decomposition: this curve is transverse to $\alpha$ and disjoint from all other closed curves in $\mathcal{P}$. Let $L$ is the maximum  of lengths of boundary components (for ${\bf T}_{0,1}$ or ${\bf T}_{0,4}$).
 
  In view of Lemma $\ref{length:set}$ and Lemma \ref{local:bound} if $ \tau_{\alpha} > 5 L,$ then the length of the corresponding curve $\beta$ is bigger than $\tau_{\alpha}$. Now the result follows from the estimate obtained in the first part for $\beta.$
 \hfill 
 $\Box$
 
 \end{section}
 
 \begin{section}{General Estimates on the counting function}\label{estimates}
In this section, we obtain basic estimates on the distributions of lengths of closed curves in a random long pants decompositions on a fixed hyperbolic surface $X$. 
 Some of these estimates might be of general interests for other problems. We remark that the results do not use the ergodic properties of the earthquake flow/ horocycle flow. We use 
 the results discussed in \ch $\ref{integrate:geometric}$ instead.

 \begin{lemm}\label{basic:estimate:1}
 Let $\mathcal{P}=\{ \alpha_1,\ldots, \alpha_{3g-3+n}\}$ be a pants decomposition on $S_{g,n}$ and $X \in \te_{g,n}.$ There exists $C(X)$ such 
 that for any ${\bf m}= (m_{1},\ldots, m_{3g-3})$ 
 we have
   $$ |\{ \g \cdot \mathcal{P}\; | \ell_{\g \cdot \mathcal{P}}(X) \leq L,  \ell_{\alpha_{1}}\leq \frac{L}{m_1}, \ldots \ell_{\alpha_{3g-3+n}}\leq \frac{L}{m_{3g-3+n}} \}| \leq C_{X} \frac{L^{6g-6+2n}}{m_1^{2} \ldots m_{3g-3+n}^2}. $$

 for $L> L({\bf m},X).$
 \end{lemm}

 \noindent
 {\bf Proof.} 
 Here we obtain a bound for 
 $$ F(X)= |\{ \g \cdot \mathcal{P}\; | \; \ell_{\g \cdot \alpha_{1}}\leq L_1, \ldots \ell_{\g \cdot \alpha_{3g-3+n}}\leq L_{3g-3+n} \}|.$$
 
In view of Theorem \ref{integrate}, we have 
 $$ \int_{\mathcal{M}_{g,n}} F(X) dX= L_{1}^2 \cdots L_{3g-3+n}^2.$$
 Now it is enough to find an embedded ball ${\bf B}_X (\epsilon_X)$ of radius $\epsilon_{X}$ around $X$ (with respect to the Thurston distance function). Note that  in this ball 
 $$ (1-\epsilon_X) < \frac{\ell_{\eta} (X)}{\ell_{\eta} (Y)} < (1+\epsilon_X).$$ Therefore, 
 $$ (1-\epsilon_X) <  \frac{F(X)}{F(Y)} < (1+\epsilon_X) .$$
 As a result we get 
 $$ |\{ \g \cdot \mathcal{P}\; | \; \ell_{\g \cdot \alpha_{1}}\leq L_1, \ldots \ell_{\g \cdot \alpha_{3g-3+n}}\leq L_{3g-3+n} \}| \;\leq C_X L_{1}^2 \ldots L_{3g-3+n}^2.$$
 One can check that the estimate works for $$C_X= \prod_{\alpha} \frac{1}{\ell_{\alpha}(X)^2},$$ where the product is over the set of closed geodesics of length $\leq \epsilon_0$ on $X.$
 
\hfill $\Box$

Let $\mathcal{P}=\{ \alpha_1,\ldots, \alpha_{3g-3+n}\}$ be a pants decomposition on $S_{g,n}.$ Let $\mathcal{P}[1]$ be the set of simple closed curves $\eta$ such that 
$\eta$ only intersects one simple closed curve $\alpha_i$ in $\mathcal{P},$ and $i(\eta, \alpha_i) \leq 2.$ Note that every simple closed curve $\eta \in \mathcal{P}[1]$ is in a once holed torus or a four holed sphere bounded by some of the curves in $\mathcal{P}.$ However, $ |\mathcal{P}[1]| = \infty.$ We prove:

 \begin{lemm}\label{basic:estimate:2}
 Let $\mathcal{P}=\{ \alpha_1,\ldots, \alpha_{3g-3+n}\}$ be a pants decomposition on $S_{g,n},$ and $X \in \te_{g,n}.$ There exists $C(X)$ such that 
  
 $$ |\{ \g \cdot \mathcal{P}\; | \; \ell_{\g \cdot \mathcal{P}}(X) \leq L, \exists \alpha \in \mathcal{P}[1], \ell_{\g \cdot \alpha}(X) \leq \frac{L}{m}\}| \leq C_{X} \frac{L^{6g-6+2n}}{m}. $$
 for $L> L(W,X).$
 \end{lemm}

Given a connected simple closed curve on $S_{g,n}$,
let $h_{\alpha} \in \Mod_{g,n} $ denote the right Dehn twist around $\alpha.$
Then as a simple corollary of Lemma \ref{bound:Bas:general}, we have:

\begin{lemm}\label{ball:stay}
Let $Y \in {\bf B}_{\gamma}(L)$ and $\alpha_1,\ldots, \alpha_k$ be simple closed curves in $\mathcal{S}(\gamma)$ such that 
$i(\alpha_i, \gamma) \not=0$. If 
$h_{\alpha_1} ^{r_1} \cdots h_{\alpha_k} ^{r_k} (Y) \in {\bf B}_{\gamma} (L),$ then for $1 \leq i \leq k$
$$ |r_i| \leq C \frac{L}{\ell_{\alpha_i}(Y)}, $$
where $C=C(\alpha_1,...\alpha_k, \gamma)$ is a  constant independent of $Y$ and $L$.
\end{lemm}

We also have:
  \begin{lemm}
Let $\alpha \in \mathcal{S}(\gamma)$ be a simple closed curve on $S_{g,n}$ and $X \in \te_{g,n}.$ There exists $C(X)$ such that 
 
  $$ \{ \g \cdot \gamma\; | \ell_{\g \cdot \gamma}(X) \leq L, \ell_{\g \cdot \alpha} (X) \leq \frac{L}{m}\}| \leq C_{X} \frac{L^{6g-6+2n}}{m}. $$
 \end{lemm}

We will have the following simple consequence of Lemma \ref{ball:stay}:

\begin{prop}
Let $\mathcal{A}= \{ \mathcal{P}_{k}= \{ \alpha^i_{1},\ldots, \alpha^i_{3g-3+n} \}$ be a sequence of pants decompositions on $S_{g,n}$ such that 
$$ \frac{|\{ \mathcal{P}_k \in \mathcal{A}\; | \ell_{\alpha^{i}_k}(X) \leq L_i\} |}{L_1^{2}\cdots L_{3g-3+n}^{2} } \rightarrow 0, $$
as $L \rightarrow \infty.$
Then 
 $$ \frac{\{ \g \cdot \gamma\; | \ell_{\g \cdot \gamma}(X) \leq L, \g \cdot \mathcal{P} \in \mathcal{A} |}{L^{6g-6+2n}} \rightarrow 0$$
 as $g \rightarrow \infty.$
\end{prop}

Fix a small $\epsilon>0$ such that any two simple closed curves of length $\leq \epsilon$ on a hyperbolic surface are disjoint. Let
$$G(X)= \prod_{\ell_{\alpha}(X) \leq \epsilon} \frac{1}{\ell_{\alpha}(X)} .$$ 
One can easily check that $G$ defines an integrable function on $\mathcal{M}_{g,n}$ with respect to the Weil-Petersson volume form \cite{M:ANN}.
Then we have the following bound on $s_{X}(L,\gamma)$ on $\te_{g,n}:$

\begin{lemm}\label{upper:bound:int}
Given $\gamma$ on $S_{g,n}$ there exists $C_{\gamma}$ such that for $X \in \te_{g,n}$ we have
$$\frac{s_{X}(L,\gamma)}{L^{6g-6+2n}} \leq C_{\gamma} \times G(X).$$
\end{lemm}

\noindent
{\bf Remarks.}
Lemma \ref{ball:stay} and Lemma \ref{upper:bound:int} are 
analogous to Lemma $5.8$ and Theorem $5.4$ in \cite{ABEM}. Even though the setting is not the same, the ideas involved in the proof are similar. 

The bounds for the counting function $s_{X}(L,\gamma)$ were first obtained by Sapir \cite{S1}. Here, the dependence in $X$ is important for us. It is essential in \cite {proof:cone:count} to bound $\frac{s_{X}(L,\gamma)}{L^{6g-6+2n}}$ by an integrable function over $\mathcal{M}_{g,n}.$\\

\noindent
{\bf Proof of Lemma \ref{upper:bound:int}.}  To simplify the notation, here we assume that $\gamma$ is filling. Let 
$\mathcal{P}=\{\alpha_1, \ldots, \alpha_{3g-3+n}\}$ be a pants decomposition of $S_{g,n}.$

There exists $c_{1}= c_{1}(\gamma,\mathcal{P})$ such that we can find $X_0 \in \Mod_{g,n} \cdot X$ 
with 
$\ell_{\gamma}(X_0) \leq c_1, $ and for each $1 \leq i \leq 3g-3+n$ 
$\ell_{\alpha_i}(X_0) \leq c_1.$

Let
$\g \in \Mod_{g,n}$  be such that $\ell_{\gamma}(\g^{-1} \cdot X_0) \leq L.$ Then 
for $1\leq i \leq 3g-3+n,$ we have  $\ell_{\alpha_i}( \g^{-1} \cdot X_0)) \leq c\; L.$ 
There is a natural map $\g \cdot \gamma \rightarrow \g \cdot \mathcal{P}$. However this map is not one-to-one. 
The kernel of this map is generated by elements of the form 
$ h_{\alpha_1}^{r_{1}} \cdots h_{\alpha_{k}}^{r_{k}}$. On the other hand
The number of the elements of $h_{\alpha_1}^{r_{1}} \cdots h_{\alpha_{k}}^{r_{k}} (Z)$ in 
${\bf B}_{\gamma}(L)$ can be estimates by Lemma $\ref{ball:stay}.$

Let $\beta_{i}= \g \cdot \alpha_i$ and $k_{i}=1+ [\tau_{\alpha_i}(\g \cdot X_0)]-m_{i},$
where $$m_{i}= \min\{ [\tau_{\alpha_i}(\g' \cdot X_0)]\; | \; \g' \cdot X_0 \in B_{\gamma}(L), \; \; 
\forall i \; \; \g' \cdot \alpha_i= \g \cdot \alpha_i \}.$$
By Lemma $\ref{ball:stay},$ $\beta_{\g}= \sum_{i=1}^{3g-3+n} k_i \cdot \beta_{i} $ is an integral multi-curve of length $\leq c' L$ on $X.$ Moreover, by the definition $\beta_{\g}$ uniquely determines $\g.$
Let $b_{X}(L)$ denote the number of integral multi-curves of length $\leq L$ on $X.$
In view of Proposition 3.6 in \cite{M:ANN}, for $L\geq 1$ we have 

$$\frac{s_{X}(L,\gamma)}{L^{6g-6+2n}} \leq c'_{\gamma} \frac{b_{X}(L)}{L^{6g-6+2n}} \leq C_{\gamma} \; G(X).$$

\hfill $\Box$

 \end{section}

\begin{section}{Limiting distribution of projections of cones in $\mathcal{M}_{g,n}$}\label{limiting}

Given a pants decomposition $\mathcal{P}$, we consider
$$H_{\mathcal{P}}^{L}=\{ X \in \te_{g,n}\; | \sum_{i=1}^{3g-3+n} \ell_{\alpha_{i}}(X)=L\} \subset \te_{g,n}.$$
Then, the Weil-Petersson volume form 
induces a volume form on $H_{\mathcal{P}}^{L}.$ On the other hand, for any finite index subgroup $\Gamma \subset \Mod_{g,n}$ the horosphere 
$H_{\mathcal{P}}^{L}$ projects in to a closed subset of $\mathcal{M}_{g,n}[\Gamma]= \te_{g,n}/
\Gamma$.
In this section, we will study the asymptotic distribution of subsets of $H_{\mathcal{P}}^{L}[\Gamma]$ in $\mathcal{M}_{g,n}[\Gamma]$.

We recall that $\mu_{\mathcal{P}}^{L}[\Gamma]$
defines a finite measure on $\mathcal{M}_{g,n}[\Gamma]$ of total mass $M_{\mathcal{P}}^{\mathcal{C}, L} < \infty$. The method used in \cite{M:AB} imply that:
\begin{equation}\label{eq:MAB}
\mu_{\mathcal{P}}^{L}[\Gamma]/L^{6g-7+2n} \rightarrow [\Mod_{g,n}:\Gamma] \times \frac{B(X)}{b_{g,n}} \mu_{wp},
\end{equation}

as $L\rightarrow \infty.$ Here, as before
\begin{equation}\label{B}
B(X)= \vol_{Th}( \{ \lambda \in \ML_{g,n}, \ell_{\lambda}(X) \leq 1\})
\end{equation}
is a continuous proper map on $\mathcal{M}_{g,n},$ and 
$b_{g,n}=\int\limits_{\mathcal{M}_{g,n}} B(X) d\mu_{wp}(X).$

In this paper,  we need a slight generalization of $ (\ref{eq:MAB})$. As before, we consider the Fenchel-Nielsen coordinates on $\te_{g,n}$ defined by $\mathcal{P}.$
For ${\bf m}=(m_1,\ldots, m_{3g-3+n}) \in {\Bbb Z}^{3g-3+n}$, define  
\begin{equation}\label{cone}
\mathcal{C}_{\mathcal{P}}^{\bf m} =\{ X \; | \; X\in \te_{g,n}, \; \forall 1\leq i \leq 3g-3+n,\; \;m_{i} \cdot \ell_{\alpha_{i}}(X) \leq \tau_{\alpha_{i}}(X) \leq (m_{i}+1)\cdot \ell_{\alpha_{i}}(X)\}.
\end{equation}

Note that the group generated by Dehn twists around $\alpha_{1},\ldots, \alpha_{3g-3+n}$ is a finite index subgroup of $\operatorname{Stab}(\mathcal{P})\subset \Mod_{g,n}.$  
Also, each $\mathcal{C}_{\mathcal{P}}^{{\bf m}}$ is a fundamental domain for the action of this group.
We fix a cone $\mathcal{C} \subset \mathcal{C}_{\mathcal{P}}^{{\bf m}}$ with respect to Fenchel-Nielsen coordinates. In this case the map $H_{\mathcal{P}}^{L} \cap \mathcal{C} \rightarrow \mathcal{M}_{g,n}[\Gamma]$ is generically one-to-one. 

Let $ \mu_{\mathcal{P}}^{\mathcal{C},L}[\Gamma]$ be the restriction of $\mu_{\mathcal{P}}^{L}[\Gamma]$ to a Fenchel-Nielsen convex cone $\mathcal{C}.$

\begin{theo}\label{distribution:main}
For any convex cone $\mathcal{C}\subset \mathcal{C}_{\mathcal{P}}^{{\bf m}}$ we have:
$$\frac{\mu_{\mathcal{P}}^{\mathcal{C},L}[\Gamma]}{L^{6g-7+2n}} \rightarrow [\Mod_{g,n}: \Gamma] \;V_{\mathcal{C}}(\mathcal{L}_0) \times  \frac{B(X)}{b_{g,n}} \; \mu_{wp},$$
as $L\rightarrow \infty.$ 
\end{theo}
Here $\mathcal{L}_0: \mathcal{C} \rightarrow {\Bbb R}_+$ is defined by $\Li_0(X)= \sum_{i=1}^{3g-3} \ell_{\alpha_{i}}(X).$\\

\subsection{Ideas of the proof of (\ref{eq:MAB}).}\label{I:ABC}
We briefly recall the arguments in \cite{M:AB} which is basically Theorem \ref{distribution:main} for $\mathcal{C}= \mathcal{C}_{\mathcal{P}}^{{\bf m}}.$
A crucial role in the argument is played by the dynamics of the earthquake flow on the bundle 
$\pi: \mathcal{P}^{1}\mathcal{M}_{g,n}[\Gamma] \rightarrow \mathcal{M}_{g,n}[\Gamma]$ of measured geodesic laminations of length one.  Given $(\lambda, X) \in \ML_{g,n} \times \te_{g,n},$ we define
$\operatorname{tw}^{t}(\lambda,X)=( \lambda, \operatorname{tw}_{\lambda}^{t}(X))$. This flow generalizes the classical Fenchel-Nielsen deformations around simple
closed geodesics.
The space $\mathcal{P}^{1}\mathcal{M}_{g,n}[\Gamma]$ carries a natural finite invariant measure $\nu_{g,n}[\Gamma]$
which projects to the volume form given by $B(X) \cdot \mu_{wp}$ on $\mathcal{M}_{g,n}[\Gamma]$.  Hence $\nu_{g,n}[\Gamma](\mathcal{P}^{1}\mathcal{M}_{g,n}[\Gamma])=b_{g,n} \cdot [\Mod_{g,n}: \Gamma].$ See \ch \ref{EPEF}. We can break the proof of the main result  into four steps as follows.\\

\noindent
{\bf 1.} First, we note that the measure $\mu_{\mathcal{P}}^{L}[\Gamma]$ is the push forward by a measure $\nu_{\mathcal{P}}^{L}[\Gamma]$ on $\mathcal{P}^{1}\mathcal{M}_{g,n}[\Gamma]$ of total mass $M_{\mathcal{P}}^{ L}$ ; in other words, for any open set 
$U \subset \mathcal{M}_{g,n}[\Gamma],$ we have
$$\mu_{\mathcal{P}}^{L} [\Gamma] (U)= \nu_{\mathcal{P}}^{ L}[\Gamma]( \pi^{-1}(U)).$$
Moreover, it is easy from the definition that any weak limit of the sequence of measures $\nu_{\mathcal{P}}^{ L}[\Gamma]$ is invariant under the earthquake flow.\\

\noindent
{\bf 2.} Then as in Theorem $5.9$ of \cite {M:AB} implies that any weak limit of the sequence $\{\nu_{\mathcal{P}}^{ L}[\Gamma]/M_{\mathcal{P}}^{ L}\}_{L}$ belongs to the Lebesgue measure class. 

To understand the limiting behavior of this sequence, one can estimate the volume of small neighborhoods 
of $\mathcal{P}^{1}\mathcal{M}_{g,n}$ with respect to the measure $\nu_{\mathcal{P}}^{L}/ M_{\mathcal{P}}^{\mathcal{C}, L}.$
More precisely, we fix a compact set $K \subset \mathcal{M}_{g,n}$, and show for any open set $U \subset P\ML_{g,n} \times K$ when $L$ is large 
\begin{equation}\label{bound}
 \nu_{\mathcal{P}}^{\mathcal{C},L}(U)/ M_{\mathcal{P}}^{\mathcal{C}, L} \leq C \nu_{g,n}(U).
 \end{equation}
Here $C$ is a constant depending only on $K.$
To prove this claim first we show a similar statement for $\mu_{\mathcal{P}}^{L}$'s. Given $ X \in K$, and 
small $\epsilon >0 $, we will define an open set ${\bf B}_X(\epsilon) \subset \te_{g,n}$ satisfying the following properties:
\begin{enumerate}

\item $$C_{2} \epsilon^{6g-6+2n} \leq \mu_{wp}({\bf B}_X(\epsilon)) \leq C_{1} \epsilon^{6g-6+2n}.$$

\item $$ \frac{1}{L^{6g-7+2n}}\mu_{\mathcal{P}}^{L}({\bf B}_X(\epsilon)) \leq C_{3} \epsilon^{6g-6+2n}.$$
\end{enumerate}
 
 Here $C,$ $C_{2}$ and $C_{3}$ are constants depending only on $K$.\\

We prove that when $L$ is large enough
\begin{equation}\label{large} 
\frac{\nu_{\mathcal{P}}^{L}(V)}{M_{\mathcal{P}}^{ L}} \leq C \cdot \nu_{g,n}(V),
\end{equation}
 where $V$ is a subset of the compact part of the moduli space, and $C$ is a constant independent of $V.$\\
The proof of equation ($\ref{large}$) relies on estimates for the number of integral points in the space of measured laminations, and known facts about convexity behavior of geodesic-length functions along the earthquake paths 
\cite{Kerckhoff:Nielsen}, \cite{Wolpert:SFN}. \\

\noindent
{\bf 3.}  Using work of Minsky and Weiss \cite{Minsky:Wei}, we then show that any weak limit 
of $\{\nu_{\mathcal{P}}^{L}[\Gamma]/M_{\mathcal{P}}^{ L}\}_{L}$ is a probability measure; namely, there is no escape of mass to infinity of the  moduli space.\\

\noindent
{\bf 4.} Finally, the ergodicity of the earthquake flow \cite{M:EE} implies that $\nu_{g,n}[\Gamma]/b_{g,n}$ is the only earthquake flow invariant probability measure on $\mathcal{P}^{1}\mathcal{M}_{g,n}[\Gamma]$ in the Lebesgue measure class. 
Therefore, as $L \rightarrow \infty,$
\begin{equation}\label{mm}
\frac{\nu_{\mathcal{P}}^{ L}[\Gamma]}{M_{\mathcal{P}}^{ L}} \rightarrow \frac{\nu_{g,n}[\Gamma]}{b_{g,n}}.
\end{equation}

In order to prove $(\ref{eq:MAB})$ it is enough to note that in view of Theorem \ref{wolp}, we have: 
\begin{lemm}\label{asympt}
Fix a pants decomposition $\mathcal{P}$ on $S_{g,n}$ and a rational cone $\mathcal{C} \subset \te_{g,n}.$
Assume that $\mathcal{C}\subset \mathcal{C}_{\mathcal{P}}^{{\bf m}}$ for some ${\bf m} \in 
{\Bbb Z}^{3g-3+n}.$
Then there exists a rational number $C_{\mathcal{P}, \mathcal{C}} >0$ such that as $L \rightarrow \infty$
$$ M_{\mathcal{P}}^{ L} \sim C_{\mathcal{P} } \cdot L^{\operatorname{dim}(\te_{g,n})-1}.$$
\end{lemm}

\subsection{Almost invariant measures for the earthquake flow.} 

In order to prove Theorem $\ref{distribution:main}$, we consider the set  
$$\mathcal{H}_{\mathcal{P}}^{\mathcal{C},L}= \{ (\frac{\sum_{i=1}^{3g-3+n} \alpha_{i}}{L},X) \;| \;\ell_{\mathcal{P}}(X)=L, X \in \mathcal{C}\} \subset \mathcal{P}^{1}\te_{g,n} .$$
We remark that 
$$\mathcal{H}_{\mathcal{P}}^{\mathcal{C}, L} \subset \mathcal{H}_{\mathcal{P}}^{L}= \{ (\frac{\sum_{i=1}^{3g-3+n} \alpha_{i}}{L},X) \;| \;\ell_{\mathcal{P}}(X)=L\}.$$

Then the measure $\mu_{\mathcal{P}}^{\mathcal{C}, L}$ on $H_{\mathcal{P}}^{\mathcal{C},L}$, induces a measure $\nu_{\mathcal{P}}^{\mathcal{C}, L}$ supported on $\mathcal{H}_{\mathcal{P}}^{\mathcal{C},L}$. The measure $\nu_{\mathcal{P}}^{\mathcal{C}, L}$ is the restriction of $\nu_{\mathcal{P}}^{L}$ to $\mathcal{H}_{\mathcal{P}}^{\mathcal{C},L}$.

By the definition, $\mu_{\mathcal{P}}^{\mathcal{C}, L}[\Gamma]$ is the push forward of the measure $\nu_{\mathcal{P}}^{\mathcal{C}, L}$ under the projection map 
$$\pi_{\Gamma}: \mathcal{P}^{1}\mathcal{M}_{g,n}[\Gamma] \rightarrow \mathcal{M}_{g,n}[\Gamma].$$
In other words, for any open set $U \subset \mathcal{M}_{g,n},$
\begin{equation}\label{push}
\mu_{\mathcal{P}}^{\mathcal{C}, L}[\Gamma](U)= \nu_{\mathcal{P}}^{L}[\Gamma]( \pi^{-1}(U)).
\end{equation}
Let  
$$ M_{\mathcal{P}}^{\mathcal{C}, L}= \nu_{\mathcal{P}}^{L}(\mathcal{P}^{1}\mathcal{M}_{g,n})=\mu_{\mathcal{P}}^{L}(\mathcal{M}_{g,n}).$$

Note that $\nu_{\mathcal{P}}^{\mathcal{C}, L}[\Gamma]$ is  not necessarily invariant  under the earthquake flow. However, we have:

\begin{lemm}\label{inv}
Let $\mathcal{P}$ be a pants decomposition on $S_{g,n}$. Let $\mathcal{C} \subset \mathcal{C}_{\mathcal{P}}^{{\bf m}}$ be a convex cone in $\te_{g,n}$ with respect to the Fenchel-Nielsen coordinates. Then any weak limit of the sequence $\{\nu_{\mathcal{P}}^{\mathcal{C}, L}[\Gamma]/M_{\mathcal{P}}^{\mathcal{C}, L}\}_{L}$ is invariant under the earthquake flow.
\end{lemm}
\noindent
{\bf Proof.} 
Note that for any $X_{1} \in \mathcal{H}_{\mathcal{P}}(1),$ we have:
$${\bf t} \cdot \tw^{r} ( X_1, \widetilde{\alpha})=  \tw^{{\bf t}^2 \cdot r} ({\bf t} \cdot X_1, \widetilde{\alpha}/{\bf t}). $$
where $\widetilde{\alpha}= \sum_{i=1}^{3g-3+n} \alpha_{i}.$

Fix $r>0.$ It is easy to show that as 
$$ \frac{\nu_{\mathcal{P}}^{\mathcal{C}, L} (\{ v\; | \; v \in  \mathcal{H}_{\mathcal{P}}^{\mathcal{C}}(L), \tw^{r}(v) \not \in  \mathcal{H}_{\mathcal{P}}^{\mathcal{C}}(L)\})}{M_{\mathcal{P}}^{\mathcal{C}, L}} \rightarrow 0$$
as $L \rightarrow \infty.$ This implies that any weak limit of the sequence is invariant under $\tw^{r}.$

\hfill $\Box$
 
 We will show that any weak limit of the sequence 
$\{\ \nu_{\mathcal{P}}^{\mathcal{C}, L}/ M_{\mathcal{P}}^{\mathcal{C}, L}\}_{L}$ is in the Lebesgue measure class.
As in \ch \ref{I:ABC} since the space $\mathcal{M}_{g,n}$ is not compact, the weak limit of a sequence of probability measures need not be a probability measure; 
 it is possible that for a sequence 
 $\mu_{m} \rightarrow \mu$ of 
 probability measures on $\mathcal{M}_{g,n}$, $\mu(\mathcal{M}_{g,n})< 1.$ In other words, we need to check that given $\delta >0,$ there is a compact subset $\mathcal{K}^{\delta} \subset \mathcal{P}^{1} \mathcal{M}_{g,n}$ such that 
$$\liminf_{L \rightarrow \infty} \frac{\nu_{\mathcal{P}}^{\mathcal{C}, L}(\mathcal{K}^{\delta})}{M_{\mathcal{P}}^{\mathcal{C}, L}} > 1-\delta.$$

 The results obtained in \cite{Minsky:Wei} allows one to estimate  the amount of time earthquake paths spend in the thin part of the moduli space. See \ch \ref{non-div}.

\begin{lemm}\label{nomasin}
Any weak limit of the sequence $\{\nu_{\mathcal{P}}^{\mathcal{C}, L}[\Gamma]/M_{\mathcal{P}}^{\mathcal{C}, L}\}_{L}$ is a probability measure in the Lebesgue measure class.
 \end{lemm}
 \noindent
 {\bf Proof.}
 It is enough to observe that there exists a constant depending only on $\mathcal{P}$ and $\mathcal{C}$ such that 
 $$  \nu_{\mathcal{P}}^{\mathcal{C}, L}[\Gamma]/M_{\mathcal{P}}^{\mathcal{C},L} (U) \leq C_{\mathcal{C}} \times \nu_{\mathcal{P}}^{L}/M_{\mathcal{P}}^{\mathcal{C}, L} (U)$$
 for any open set $U$. Now the lemma follows from the same statements for $\nu_{\mathcal{P}}^{L}/M_{\mathcal{P}}^{\mathcal{C}, L} $. See  (\ref{mm}) and also Theorem 5.9 and Theorem 5.11 in \cite{M:AB}. 
 
 \hfill $\Box$

\noindent
{\bf Proof of Theorem \ref{distribution:main}.}
  Consider a weak limit $\tilde{\nu}$ of the sequence $\nu_{\mathcal{P}}^{\mathcal{C},L}/
  M_{\mathcal{P}}^{\mathcal{C}, L}$; namely,
 $$ \nu_{\mathcal{P}}^{\mathcal{C}, L_{i}}/ M_{\mathcal{P}}^{\mathcal{C}, L_i}\rightarrow \tilde{\nu},$$
as $L_{i} \rightarrow \infty.$
 Then by Lemma $\ref{nomasin}$, $\tilde{\nu}$ is a probability measure in the Lebesgue measure
 class. On the other hand, by Lemma $\ref{inv}$ the limit measure is also invariant under the earthquake flow. Hence, by the ergodicity of the earthquake flow (Theorem $\ref{er}$), we have
 $$\tilde{\nu}=\frac{\nu_{g,n}}{b_{g,n}}.$$ Now the result follows from Lemma $\ref{asympt}.$
\hfill $\Box$

\subsection{Remark}\label{remark:NF}

We will need a slight generalization of Theorem \ref{distribution:main} later in \ch $\ref{proof:morege}$ for dealing with the case that $\gamma$ is not filling.

 Let $\mathcal{P}$ be a pants decomposition of $S_{g,n}$ and let ${\bf S}$ be a subsurface of $S_{g,n}$ such that $\partial ({\bf S}) \subset \mathcal{P}.$
Let $\mathcal{P}_{{\bf S}}= \mathcal{P} \cap ({\bf S} \cup \partial ({\bf S}))$ and $\Gamma_{{\bf S}}=
\Stab({\bf S}).$ Then the space $\te_{g,n}/ \Gamma_{{\bf S}}$ has a parametrization by the length and twist around curves in $\mathcal{P}_{{\bf S}}.$ 

The Weil-Petersson volume form 
induces a measure on 
$$H_{\mathcal{P}_{{\bf S}}}^{L}=\{ X \; |\;  \sum_{\alpha \in \mathcal{P}_{{\bf S}}} \ell_{\alpha}(X)=L \} \subset \te_{g,n}/ \Gamma_{{\bf S}}.$$

One can similarly define $\mathcal{C}_{\mathcal{P}_{{\bf S}}}^{\bf m}$ and $\mu_{\mathcal{P}_{{\bf S}}}^{\mathcal{C}_{{\bf S}},L}.$

Let $\mathcal{C}_{{\bf S}} \subset \te_{g,n}/ \Gamma_{{\bf S}}$ be a convex cone in $\mathcal{C}_{\mathcal{P}_{{\bf S}}}^{\bf m}.$ 

Define $\mathcal{L}_0: \mathcal{C}_{{\bf S}} \rightarrow {\Bbb R}_+$ is  by $\Li_0(X)= \sum_{\alpha \in \mathcal{P}_{{\bf S}}}  \ell_{\alpha}(X).$
Following Theorem \ref{volm} and  Theorem $\ref{integrate}$ we have
$$ M_{\mathcal{P}_{{\bf S}}}^{L} \sim C_{\mathcal{P}_{{\bf S}}} \cdot L^{\operatorname{dim}(\te_{g,n})-1}.$$

Then the proof of  Theorem \ref{distribution:main} implies that:
$$\frac{\mu_{\mathcal{P}_{{\bf S}}}^{\mathcal{C}_{{\bf S}},L}}{L^{6g-7+2n}} \rightarrow \;V_{\mathcal{C}_{{\bf S}}}(\mathcal{L}_0) \times  \frac{B(X)}{b_{g,n}} \; \mu_{wp},$$
as $L\rightarrow \infty.$

 \end{section}

\begin{section}{From equidistribution of horospheres to counting estimates}\label{from-eq-to-count}

Let $\mathcal{P}$ be a pants decomposition of $S_{g,n}$ and  $\mathcal{C} \subset \te_{g,n}$ be a cone in the Fenchel-Nielsen coordinates For finite index subgroup $\Gamma$ of $\Mod_{g,n}.$ Our goal is to understand the asymptotics of the number of elements in $$\Gamma \cdot X \cap \mathcal{C}_{L} \subset \te_{g,n}, 
$$ where $$\mathcal{C}_{L}= \{ X \in \mathcal{C}, \Li_0(X) \leq L\},$$
and $\Li_0(X)= \sum_{i=1}^{3g-3+n} \ell_{\alpha_i}(X).$
In this section, we will use the equidistribution results obtained in \ch $\ref{limiting}$ to prove the following result:
\begin{theo}\label{cone:count} 
Let $\mathcal{C} \subset \mathcal{C}_{\mathcal{P}}^{{\bf m}}$ for some ${\bf m}= (m_{1},\ldots, m_{3g-3+n}) \in {\Bbb Z}^{3g-3+n}$ be a convex cone. Then for any $X \in \te_{g,n}$ we have
\begin{equation}\label{countl2}
\#(\{ \g \in \Gamma\; | \; \g \cdot X \in \mathcal{C}_{L}\}) \sim [\Mod_{g,n}: \Gamma] \times \operatorname{V}_{\mathcal{C}}(\mathcal{L}_0) \cdot B(X) \cdot L^{6g-6+2n}, 
\end{equation}
as $L \rightarrow \infty.$ \end{theo}
As before $V_{\mathcal{C}}=\operatorname{V}_{\mathcal{C}}(\mathcal{L}_0)= \vol_{wp}(\mathcal{C}^{+}_{1}) < \infty.$
Note that an arbitrary cone can be decomposed into finitely many such cones as long as the ratios of $\tau_{i}(X)/\ell_{\alpha_i}(X)$ stays uniformly bounded on the cone. 

We remark that in Theorem $\ref{cone:count},$ there is no difference in the arguments for general case of $\Gamma$ and $\Mod_{g,n}.$ To simplify the notation, we present the argument for $\Gamma=\Mod_{g,n}$ and comment on the general case in \ch $\ref{comment:general}.$

\subsection{Basic Idea of proof of Theorem $\ref{cone:count}$ and the Main Counting Lemma}
  
Let $$\mu_+^{{\bf D}, L} =\frac{\int _{0}^{L} \mu^{{\bf D},t} dt}{L},$$
where $\mu^{{\bf D},t}=\mu_{\mathcal{P}}^{{\bf D},t}[\Gamma]$ as defined in \ch $\ref{limiting}$
Consider the map $\pi: \te_{g,n} \rightarrow \mathcal{M}_{g,n}.$ Then Wolpert's theorem (Theorem \ref{wolp}) implies that:
\begin{lemm}\label{mu:vol}
For any open set $U \subset \mathcal{M}_{g,n},$ we have
$$\mu_+^{{\bf D}, L} (U)= \vol_{wp} (\pi^{-1}(U) \cap {\bf D}_{L}).$$ 
\end{lemm}

Let $X \in U.$ In view of Lemma \ref{mu:vol}, since the Weil-Petersson volume form is invariant under the action of $\Mod_{g,n}$ we expect  
$| \Mod_{g,n} \cdot X \cap {\bf D}_{L}| \times \vol_{wp}(U) $ to be a good approximation of $\mu_+^{{\bf D}, L} (U).$ The problem in this approach is taking care of the edge effect: some parts of 
$\g \cdot U$ could be outside of ${\bf D}_{L}$ and it is possible that $\g \cdot U \cap {\bf D}_{L} \not = \emptyset$ but $\g \cdot X \not \in  {\bf D}_{L}.$

However, if we consider the ball ${\bf B}_{X}(\epsilon)$ of radius $\epsilon$ with respect to the Thurston distance function (\ch \ref{Th:dis}), we can approximate the shape of 
$\g \cdot {\bf B}_{X} (\epsilon)$ using the results of \ch \ref{local:calculation}. The estimates hold as long as $\g \cdot X \in \mathcal{B}_{\mathcal{P}} (M)$ and the lengths of all closed curves in $\mathcal{P}$ are large enough on $\g \cdot X.$

Here, as before, $$\mathcal{B}_{\mathcal{P}} (M)= \{ X \; | \; \forall   1 \leq i,j \leq 3g-3+n\;,  \frac{\ell_{\alpha_{i}}(X)}{\ell_{\alpha_{j}}(X)} < M \; \}  \subset \te_{g,n}.$$

\begin{lemm}\label{cones}
Let $\mathcal{P}$ be a pants decomposition of $S_{g,n}$ and $M>0$. For any convex cone ${\bf D} \subset \mathcal{C}_{\mathcal{P}}^{{\bf m}},$ there exists $L_{0}$ and $\epsilon_0 $ such that  the following holds.
For any $0 < \epsilon< \epsilon_0$  there are cones $ {\bf D}^{+} (M, \epsilon)$ and $ {\bf D}^{-} (M, \epsilon)$ such that: 
\begin{enumerate}
\item $ {\bf D}^{-} (M, \epsilon) \subset {\bf D} \subset {\bf D}^{+} (M, \epsilon),$
\item If  $Z \in {\bf D} \cap \mathcal{B}_{\mathcal{P}} (M)$ such that $L_{\mathcal{P}}(Z)> L_0$ and $d_{Th}(Z,Z')\leq \epsilon$ then $Z' \in {\bf D}^{+} (M, \epsilon),$
\item If $Z \in {\bf D}^{-}(M, \epsilon) \cap \mathcal{B}_{\mathcal{P}} (M)$ such that  $L_{\mathcal{P}}(Z)> L_0$ and $d_{Th}(Z,Z') \leq \epsilon$ then $Z' \in {\bf D},$
and 
\item $$ (1- C \; \epsilon M^2)^{6g-6+2n}< \frac{\vol({\bf D}^{-}_{L}(M, \epsilon))}{\vol({\bf D}_{L})}< 1< \frac{\vol({\bf D}^{+}_{L}(M, \epsilon))}{\vol({\bf D}_{L})} < (1+C \epsilon M^2)^{6g-6+2n}$$
where $C$ is a constant independent of $M$, $\epsilon$ and $L$.
\end{enumerate}

\end{lemm}

\noindent
{\bf Proof of Lemma \ref{cones}}. 
Let $\mathcal{P}= \{\alpha_{1}, \ldots, \alpha_{3g-3+n}\}.$
Let $C$ be the constant in Corollary \ref{coro:est}.
Then we can define
 $$ {\bf D}^{+} (M, \epsilon)= \{ \prod_{i=1}^{3g-3+n} h_{\alpha_i} ^{r_{i}} (Z) \; | \; Z \in {\bf D},\; \; |r_{i}| \leq C M^2 \epsilon \ell_{\alpha_i}(Z) \; \},$$
 and
$$ {\bf D}^{-} (M, \epsilon)= \{ Z\;  | \; \prod_{i=1}^{3g-3+n} h_{\alpha_i} ^{r_{i}} (Z) \in {\bf D} \; ;  \forall |r_{i}| \; \leq C \; M^2 \epsilon\;  \ell_{\alpha_i}(Z) \; \}.$$
 
Here $h_{\alpha}(X)$ is the twisting of $X$ around $\alpha$ as in $(\ref{D:T}).$
 
Parts ${\bf 2}$ and ${\bf 3}$ are direct corollaries of  Corollary \ref{coro:est}. Finally, when $\epsilon$ is small enough (depending on $M$ and ${\bf D}$) Part ${\bf 4}$ is a consequence of Theorem $\ref{wolp}.$

\hfill $\Box$\\

As a result, we get:

\begin{prop}\label{counting:lemma:1}
Let $\mathcal{C}\subset \mathcal{C}_{\mathcal{P}}^{({\bf m})}$ for some ${\bf m}= (m_{1},\ldots, m_{3g-3+n}) \in {\Bbb Z}^{3g-3+n}$. Let $\Gamma$ be a finite index subgroup of $\Mod_{g,n}$ and $X \in \te_{g,n}$. Then for any $M>0$ and small  enough $\epsilon$ (depending on $M$ and $X$, $\Gamma$ and $\mathcal{C}$),  we have: 

\begin{equation}\label{lower}
\mu_{+}^{\mathcal{C}^{-}(M,\epsilon),L} [\Gamma]( {\bf B}_{X} (\epsilon))- \vol_{wp}(\mathcal{C}_{L} -  \mathcal{B}_{\mathcal{P}}(M))  \leq \#( \Gamma \cdot X \;\cap\; \mathcal{C}_{L}) \times \vol_{wp} ({\bf B}_{X} (\epsilon))+ R_1(\epsilon, M),
\end{equation}
and 
\begin{equation}\label{upper}
 \#( \Gamma \cdot X \; \cap\; \mathcal{C}_{L} \cap \mathcal{B}_{\mathcal{P}} (M) ) \times \vol_{wp} ({\bf B}_{X} (\epsilon))- R_2(X, \epsilon,M) \leq \mu_{+} ^{\mathcal{C}^{+}(M, \epsilon),L)}[\Gamma] ({\bf B}_{X} (\epsilon)),
 \end{equation}
where the term $R_1(\epsilon, M)$ and $R_2(X, \epsilon,M)$ are independent of $L.$
\end{prop}

\noindent
{\it Sketch of Proof of Lemma \ref{counting:lemma:1}.}
The result follows from a simple estimate for the Weil-Petersson volume of $ \mathcal{C}_{L} \cap \pi^{-1}({\bf B}_{X}(\epsilon))$ and Lemma \ref{cones}. 
Let $L_0$ is the number given by Lemma \ref{cones}. 
Consider the function
$$R_{1}(\epsilon, M)= \vol_{wp}(\{ Z\; | \; Z \in \mathcal{C}\;, \forall \alpha_i \in \mathcal{P}\; \; \ell_{\alpha_i}(Z) \leq M\cdot L_0\}).$$
Then in view of Lemma \ref{mu:vol} the inequality $(\ref{lower})$ follows from Lemma \ref{simple:bound}, and part $(3)$ of Lemma $\ref{cones}.$ 
Also, define
$$R_{2}(X,\epsilon, M)= \# (\{ Z\; | \; Z \in \Mod_{g,n} \cdot X\; \forall \alpha_i \in \mathcal{P} \; \;\ell_{\alpha_i}(Z) \leq M\cdot L_0\}).$$
Inequality $(\ref{upper})$ follows from Lemma \ref{simple:bound}, and part $(2)$ of Lemma $\ref{cones}.$

\hfill $\Box$

\subsection{Proof of Theorem \ref{cone:count}}\label{proof:cone:count}
To simplify the notation, we present the proof in case we have $\Gamma=\Mod_{g,n},$ and $|\Aut(X)|=1.$ Also for any cone $\mathcal{D},$ let 
$V_{\mathcal{D}}=\operatorname{V}_{\mathcal{D}}(\mathcal{L}_0).$

By Lemma $\ref{simple:bound},$ we have
$$\vol_{wp}(\mathcal{C}_{L} -  \mathcal{B}_{\mathcal{P}}(M))= O_{X}(\frac{L^{6g-6+2n}}{M}).$$
On the other hand, following the general estimates proved in \ch \ref{estimates} (See Lemma \ref{basic:estimate:1} and Lemma \ref{basic:estimate:2}), we have
$$ \#( \Mod_{g,n} \cdot X \cap (\mathcal{C}_{L} -  \mathcal{B}_{\mathcal{P}} (M))= O_{X}(\frac{L^{6g-6+2n}}{M}).$$ 

We can choose a subsequence $\{L_{k}\}_{k}$ with $L_{k} \rightarrow \infty $ as $k \rightarrow \infty$
and 
$$\lim _{k\rightarrow \infty} \frac{\#( \Mod_{g,n} \cdot X \cap \mathcal{C}_{L_k})}{L_k^{6g-6+2n}}= a.$$
Our goal is to prove that $a$ is independent of this subsequence.
Note that by Theorem $\ref{distribution:main}$ for any bounded cone $\mathcal{D}$
\begin{equation}\label{limit:use}
 \lim_{L \rightarrow \infty} \frac{\mu_{+}^{{\bf D},L} ( {\bf B}_{X} (\epsilon))}{L^{6g-6+2n}}= \operatorname{V}_{\bf D} \times Q(X, \epsilon),
 \end{equation}
where 
$$ Q(X, \epsilon)= \int_{{\bf B}_{X}(\epsilon)} B(X) d\mu_{wp}.$$
 Fix $0\ll M.$ Let $\epsilon\ll \frac{1}{M}$. We can use $(\ref{limit:use})$ for $\mathcal{C}^{-}(M,\epsilon)$ and $\mathcal{C}^{+}(M,\epsilon).$ By $(\ref{upper})$ and $(\ref{lower})$ we have
 $$    \frac{Q(X, \epsilon)}{\vol_{wp} ({\bf B}_{X}(\epsilon))} \times  \operatorname{V}_{\mathcal{C}^{-}(M, \epsilon)} -\frac{C_2}{M} - \lim_{k \rightarrow \infty}\frac{R_2(X, M,\epsilon)}{L_{k}^{6g-6+2n}}  \leq $$ $$\leq \lim_{k \rightarrow \infty}\frac{\#( \Mod_{g,n} \cdot X \cap \mathcal{C}_{L_k})}{L_{k}^{6g-6+2n}} \leq$$ 
 $$ \leq  \frac{Q(X, \epsilon)}{\vol_{wp} (\mathcal{B}_{X}(\epsilon))} \times \operatorname{V}_{\mathcal{C}^{+}(M, \epsilon)} + \frac{C_1}{M}+\lim_{k \rightarrow \infty} \frac{R_1(\epsilon,M)}{L_{k}^{6g-6+2n}} $$
 where $C_1$ and $C_2$ are independent of $M$ and $\epsilon$.
 
 Let $\epsilon \rightarrow 0.$ Then by part $4$ of Lemma $\ref{cones}$ we get 
  $$B(X) \times  \operatorname{V}_{\mathcal{C}} -\frac{C_2}{M}- \frac{R(X, M,\epsilon)}{L_{k}^{6g-6+2n}}  \leq \lim_{k \rightarrow \infty}\frac{\#( \Mod_{g,n} \cdot X \cap \mathcal{C}_{L_k})}{L_{k}^{6g-6+2n}} \leq  B(X)  \times \operatorname{V}_{\mathcal{C}} + \frac{C_1}{M}, $$
 
Now the result follows if we let $M \rightarrow \infty.$
\hfill $\Box$\\

\subsection{Comments on the general case} \label{comment:general}
There is no difference in the proof for the general case of $\Gamma$ . The main ingredient of the proof is that the horocycle flow is ergodic on the cover 
 $\mathcal{P}^{1}\mathcal{M}_{g,n}[\Gamma]$ corresponding to any $\Gamma.$

\end{section}

\begin{section}{Volume estimates}\label{vol-est}
 
Let $\gamma$ be a filling closed curve. 
$${\bf B}_{\gamma}(L)=\{ X \in \te_{g,n}, \ell_{\gamma}(X) \leq L\} $$
is compact. It is easy to see that if $\gamma$ is not filling ${\bf B}_{\gamma}(L)$ has infinite volume, as  $\operatorname{Stab}(\gamma) \subset \operatorname{Mod}(S_{g,n})$ is infinite. 
\begin{theo}\label{vol-ball}
Let $\gamma$ be a filling closed curve on $S_{g,n}.$ As $L \rightarrow \infty$
$$\operatorname{Vol}_{wp}({\bf B}_{\gamma}(L))\sim L^{6g-6+2n} v_{\gamma}.$$
More generally, 
$$\operatorname{Vol}_{wp}({\bf B}_{\gamma}(L)/\operatorname{Stab}(\gamma))\sim L^{6g-6+2n} v_{\gamma}$$
as $L \rightarrow \infty.$
Moreover for any connected closed curve $v_{\gamma} \in {\Bbb Q}$
\end{theo} 
In particular, for any $\gamma$ closed cure in $\pi_{1}(S_{g,n})$ we have
$$\int_{\mathcal{M}_{g,n}} s_{X}(L,\gamma) dX \sim L^{6g-6+2n} v_{\gamma},$$
where the volume is taken with respect to the Weil-Petersson volume form.  We need the following lemma: 
\begin{lemm}\label{basic:vol:bound}
Given a closed curve $\gamma$ we have 
$$\operatorname{Vol}_{wp}({\bf B}_{\gamma}(L)/ \operatorname{Stab}(\gamma) ) < C_{\gamma} L^{6g-6+2n}.$$
Here the constant $C_{\gamma}$ is independent of $L$.
\end{lemm}

This is due to Sapir \cite{S1}. Here we sketch the proof for completeness.

\noindent
{\it Sketch of proof of Lemma \ref{basic:vol:bound}.} If $i(\gamma, \gamma) \not =0,$ it is enough to choose a subsurface ${\bf S}(\gamma)$ such that $\gamma$ is filling in ${\bf S} (\gamma).$ For every $Z \in \te_{g,n},$ we can extend a pants decomposition of ${\bf S}(\gamma)$ to a pants decomposition $\mathcal{P}_{Z}$ of $S_{g,n}$ such that for every
$$\ell_{\alpha_i}(Z) < C \ell_{\gamma}(Z),$$
where $C$ is independent of $Z.$
This is because every hyperbolic surface with geodesic boundary components of total length $L$ has a pants decomposition with closed curves of length of order $L$.  
The subsurface ${\bf S}(\gamma)$ is bounded by curves $\mathcal{S}(\gamma)$. Then for any $\alpha$
in the surface ${\bf S}(\gamma)$ or in $\mathcal{S}(\gamma)$ we have $\Stab(\gamma) \subset \Stab(\alpha).$ 
Now the lemma follows from Theorem $\ref{wolp}$ and the proof of Lemma $\ref{twist:bound}.$

\hfill $\Box$\\

We will also need the following straightforward corollary of Theorem $\ref{wolp}$:

\begin{lemm}\label{simple:bound}
 For any convex cone $\mathcal{C}\subset \mathcal{C}_{\mathcal{P}}^{{\bf m}}$ and let $\mathcal{W}$ be a linear function in terms of the Fenchel-Nielsen coordinates, $\{\ell_{i}, \tau_{i}\}.$
 $$\frac{\vol_{wp}( \{ X\; | \; \mathcal{W}(X) \leq \sqrt L\} \cap \mathcal{C}_{L})}{L^{6g-6+2n}} \rightarrow 0$$
 $$\frac{\vol_{wp}( \mathcal{B}^c_{\mathcal{P}} \cap \mathcal{C}_{L})}{L^{6g-6+2n}}= O(\frac{1}{M})$$
 \end{lemm}

 In other words, we have:
 $$1- \frac{C_1}{M} \leq \frac {\vol_{wp} (\mathcal{C}_{L} \cap \mathcal{B}_{\mathcal{P}} (M))}{\vol_{wp} (\mathcal{C}_{L})},$$
where $C_1$ is independent of $M$ and $L$.

\noindent
{\it Sketch of proof of Theorem $\ref{vol-ball}$.} 
Let $F$ be asymptotically linear approximated by $\Li$ in a cone $\mathcal{C}.$ Assume that 
$ \vol( F ({\bf x}) \leq L)/L^{m} $ is bounded. By the definition, we can find $\{\mathcal{R}_{i}\}_{i}$ and $c>0$ such that the following holds: for every $\epsilon>0,$ there exists $M(\epsilon)>0$ such that if for all $i$, $|\mathcal{R}_{i}({\bf x})| > M(\epsilon),$ then $$c -\epsilon<  F({\bf x})- \Li({\bf x}) < c+ \epsilon.$$

Therefore, in view of Lemma $\ref{simple:bound}$ we have
$$ \frac{\vol( F ({\bf x}) \leq L) } {\vol( \LL ({\bf x}) \leq L) } \rightarrow 1,$$
as $L \rightarrow \infty.$

If $\gamma$ is filling, then the result follows since by Theorem $\ref{LAPL},$ $\ell_{\gamma}$ is asymptotically piecewise linear and the linear functions approximating $\ell_{\gamma}$ have rational coefficients and each cone is defined by finitely many rational linear functions.
If $\gamma$ is not filling, as before \ch \ref{non-filling} we choose a pants decomposition $\partial({\bf S}) \subset \mathcal{P}.$ Now in view of Theorem $\ref{volm}$ the same argument works on  
$\te_{g,n}/ \Gamma_{\gamma}.$

\hfill $\Box$

\end{section}

\begin{section}{Proof of Theorem \ref{morege} and Theorem \ref{distribution:pants}.}\label{proofs}
 
 In this section, we prove the main results of this paper. We start with the proof of Theorem \ref{distribution:pants} which is simpler.

 \subsection{Proof of Theorem \ref{distribution:pants}.}
 We consider the Fenchel-Nielsen coordinates corresponding to $\mathcal{P}$ on $\te_{g,n}.$ Note that any $\mathcal{C}_{\mathcal{P}}^{{\bf m}}$ is a fundamental domain for the action of  $ \operatorname{Stab}(\mathcal{P})$ on $\te_{g,n}.$

Hence the following corollary of Theorem \ref{cone:count} is equivalent to Theorem $\ref{distribution:pants}$.
As in the introduction, we have $$\Delta= \{(x_1,\ldots, x_{3g-3+n}) | \; \sum x_{i}=1\} \subset {\Bbb R}_{+}^{3g-3+n}.$$ 
Given $A \subset \Delta,$ let

$$\widehat{A}_{L}= \{(t\; \cdot y_1,\ldots, t \;\cdot  y_{3g-3+n})\; | 0 \leq t \leq L\;,  (y_i)_{i} \in A\; \}.$$

 \begin{coro}\label{eq:theo2}
 Let $X \in \te_{g,n},$  and $\mathcal{P}=\{\alpha_{1},\alpha_{2},\ldots, \alpha_{3g-3+n}\}$ be a pants decomposition of $S_{g,n}.$ Let $A \subset \Delta.$ Then we have:
 
 \begin{equation}
 \#(\{ \g \; | \; (\ell_{\g \cdot \alpha_i}(X))_{i} \in \widehat{A}_{L}\}) \sim [\Mod_{g,n}: \Gamma] \times c \cdot B(X) \cdot L^{6g-6+2n} ,
 \end{equation}
as $L \rightarrow \infty.$
 \end{coro}

  \subsection{Idea of the Proof of Theorem \ref{morege}}\label{proof:main:idea}
 Similarly, we would like to approximate the number of points  
 $$\Mod_{g,n} \cdot X \cap \mathcal{C} \cap {\bf B}_{\gamma}(L) \subset \te_{g,n},$$
 using the results about the asymptotically linear behavior of $\ell_{\gamma}$ on $\te_{g,n}.$

 There are two main technical issues in this approach for provingTheorem $\ref{morege}$:\\

 \noindent
{\bf I.} In general, we do not know much about the linear functions approximating $\ell_{\gamma}$ on $\te_{g,n}.$ We need the counting results in cones for asymptotically piecewise linear functions (See Proposition $\ref{count:PL}$). We will also need the following estimate:
\begin{lemm}\label{cone:count:upper} 
Let $\mathcal{C} \subset \mathcal{C}_{\mathcal{P}}^{{\bf m}}$ for some ${\bf m}= (m_{1},\ldots, m_{3g-3+n}) \in {\Bbb Z}^{3g-3+n}$. Let $\Li$ be a linear function in terms of the Fenchel-Nielsen coordinates. 
Then for any sequence $a_{k} \rightarrow 0,$ we have
\begin{equation}\label{countl2:2}
\frac{\#(\{ \g \cdot X\; | \; \g \in \Mod_{g,n} \;  \Li(\g \cdot X) \leq a_k \cdot L_{k} \} \cap \mathcal{C}_{L_{k}})}{L_{k}^{6g-6+2n}} \rightarrow 0, 
\end{equation}
as $k \rightarrow \infty.$
\end{lemm}

\noindent
{\bf II.} We need to use part $({\bf I})$ for infinitely many cones, as the counting results for cones only holds for cones in $\te_{g,n}$ on which $\ell_{\alpha}/\tau_{\alpha}$ for each 
$\alpha \in \mathcal{P}$ is bounded uniformly on $\mathcal{C}$. The reason for the difference between this case and the case of simple closed curves is that $\operatorname{Stab}(\alpha)$ is large, but 
$\operatorname{Stab}(\gamma)$ is a finite group if $\gamma$ is filling.

\subsection { Part (I):  Counting in one bounded cone.}

 Let ${\bf D}$ be a cone in ${\Bbb R}^{m}=\{ (x_1,\ldots,x_{m}), x_{i} \in {\Bbb R}\}.$ 
We consider the Euclidean volume form $\operatorname{Vol}_{E}$ is defined using the form $d x_{1} \cdots dx_{m}$ in ${\bf D}$. 
  We say a linear function $\mathcal{L}$ {\it is of compact type} in ${\bf D}$ if $V_{\bf D}(\Li)= \vol_{E}(\{{\bf x}\; | {\bf x} \in {\bf D}, \Li({\bf x}) \leq 1\})$ is finite.
  In this case $ \{{\bf x}\; | \; {\bf x} \in {\bf D}, \Li({\bf x}) \leq 1\}$ has compact closure.

Let  $${\bf P}=\{P_1,\ldots, P_{m},\ldots\} \subset {\bf D}$$ be a discrete set of points inside this cone.
We say that ${\bf P}$ becomes equidistributed with respect to the ${\Bbb R}$-action on $({\bf D}, \Li_0)$ iff for any open set 
$U \subset \mathbb{P}({\bf D}),$ we have
 $$| \{ P \in {\bf P}\; | \Li_0(P) \leq T,  [P] \in U\}| \sim c_{{\bf P},\Li_0} \times T^{m} \times V_{\bf D, U}(\Li_0)$$ 
 as $T \rightarrow \infty .$
 Here $V_{\bf D, U}(\Li)= \vol_{E}(\{{\bf x}\; | {\bf x} \in {\bf D}, \Li({\bf x}) \leq 1, [{\bf x}] \in U\})$ and  and $c_{{\bf P}, \Li_0} \in {\Bbb R}_+$ is independent of $U$. We will use the following elementary observation: 
 
 \begin{prop}\label{count:PL}
 In terms of the above notation, assume that the set ${\bf P}$ becomes equidistributed with respect to the ${\Bbb R}$-action on $({\bf D}, \Li_0).$
 Then for any asymptotically  linear function $F$ approximated by a linear function $\Li$ of compact type,  we have:
 $$|{\bf P} \cap \{F({\bf x}) \leq T\}| \sim c_{{\bf P}, \Li}\times  T^{m}$$ as $T \rightarrow \infty $, where
 $$ \frac{c_{{\bf P}, \Li} } {c_{{\bf P}, \Li_0}}= \frac{V_{\bf D}(\Li)}{V_{\bf D}(\Li_0)}. $$
 \end{prop}

\noindent
{\it Sketch of Proof of Proposition $\ref{count:PL}$.}
First, assume $F= \Li$ is a linear function. 
The main idea is to approximate the set $\{{\bf x} \in {\bf D}, \Li({\bf x}) \leq 1\}$ by pieces of the level sets of $\Li_{0};$ this way we can approximate 
$\{{\bf x} \in {\bf D}, \Li({\bf x}) \leq L\}$ by pieces of the form
$\{{\bf x} \in {\bf D}, \Li_0({\bf x}) \leq a_{i} L, [{\bf x}] \in U_{i}\},$ where $\bigcup_{i} U_{i}= \mathbb{P}({\bf D}).$

For the general case of asymptotically piecewise linear functions, we will use the definition of asymptotically piecewise linear.
Exactly as in the proof of Lemma $\ref{cone:count:upper},$ for any linear function $\mathcal{R}$ and any sequence $a_{k} \rightarrow 0,$ we have
\begin{equation}\label{countl2:2:2}
\frac{\#(\{ P \in {\bf P} \; | \;  \mathcal{R}(P) \leq a_k \cdot L_{k} \} \cap {\bf D}_{\Li_0})}{L_{k}^{m}} \rightarrow 0, 
\end{equation}
as $k \rightarrow \infty.$ Now we can use the argument in the proof of Theorem $\ref{vol-ball}$.
\hfill $\Box$\\

\noindent
{\bf Proof of Lemma $\ref{cone:count:upper}$ using Theorem $\ref{cone:count}$.}
Given $\epsilon>0$, one can find a cone $\mathcal{C}_{\epsilon}$ with $\operatorname{V}(\mathcal{C}_{\epsilon}) < \epsilon$ such that for $k$ large enough we have 
$$\{ \g \cdot X\; | \; \g \in \Mod_{g,n} \;  \Li(\g \cdot X) \leq a_k \cdot L_{k} \} \cap \mathcal{C}_{L_{k}} \subset \mathcal{C}_{\epsilon}.$$
Now in view of Proposition \ref{count:PL}, the result follows from using $(\ref{countl2})$ for $\mathcal{C}_{\epsilon},$ and letting $\epsilon \rightarrow 0.$
\hfill $\Box$\\

 \subsection{Proof of Theorem \ref{morege}}\label{proof:morege}

First, to simplify the notation assume that $\gamma$ is filling. Fix $X \in \te_{g,n},$ and let $\mathcal{P}$ be a pants decomposition of $S_{g,n}.$ In view of Proposition $\ref{count:PL}$,  Theorem $\ref{cone:count}$ and Theorem $\ref{LAPL}$ imply that 
for any ${\bf m} \in {\Bbb Z}^{3g-3+n},$ we have 

\begin{equation}\label{count:in:cones}
\#(\{ \g \in \Gamma\; | \; \g \cdot X \in \mathcal{C}_{\mathcal{P}}^{{\bf m}}\; , \; \ell_{\gamma} (\g \cdot X) \leq L\}) \sim [\Mod_{g,n}: \Gamma] \times  c_{\gamma, {\bf m}} \; B(X) \cdot L^{6g-6+2n}, 
 \end{equation}
as $L \rightarrow \infty.$ 

See $(\ref{cone}).$ We remark that we need infinitely many of these cones of the form $\mathcal{C}_{\mathcal{P}}^{{\bf m}}$ to cover $\te_{g,n}.$  However, we have the following:\\
 
\noindent
{\it  Claim.} If $\g \cdot X \in \mathcal{C}_{\mathcal{P}}^{{\bf m}}$ with $\ell_{\gamma}(\g \cdot X) \leq L,$ then 
 $$ \ell_{\alpha_{i}}(\g \cdot X) = O_{X}(\frac{L}{m_i}),$$
 where ${\bf m}=(m_1,\ldots, m_{3g-3+n}) \in {\Bbb Z}^{3g-3+n}.$

Let $\beta_{\Sigma}$ be the simple closed curve used for defining $\tau_{\alpha}$ for each pair of pants $\Sigma$ in the pants decomposition of $\mathcal{P}.$ 
Then the claim follows from Lemma $\ref{local:bound}$ and Lemma $\ref{bound:Bas:general}.$

On the other hand, since
$$\sum_{\bf m \in {\Bbb N}^{3g-3+n}.} \frac{1}{m_{1}^2} \times\ldots \times \frac{1}{m_{3g-3+n}^2} < \infty,$$ 
in view of Lemma $\ref{basic:estimate:1}$ we can use the dominated convergence theorem and 
$(\ref{count:in:cones})$ to prove $(\ref{main:theo}).$ Moreover, we have  
$$n_{\gamma} = \sum_{{\bf m} \in {\Bbb Z}^{3g-3+n}} c_{\gamma, {\bf m}}.$$
In particular, $n_{\gamma}$ is independent of $X$. 

In view of Lemma $\ref{upper:bound:int},$ dominated convergence theorem implies that
$$ \int_{\mathcal{M}_{g,n}} \lim_{L \rightarrow \infty} \frac{s_{X}(L,\gamma)}{L^{6g-6+2n}} dX = 
\lim_{L \rightarrow \infty} \int_{\mathcal{M}_{g,n}}  \frac{s_{X}(L,\gamma)}{L^{6g-6+2n}} dX.$$

 On the other hand, 
 $$ \int_{\mathcal{M}_{g,n}}  \frac{s_{X}(L,\gamma)}{L^{6g-6+2n}} dX= \frac{\operatorname{Vol}_{wp}({\bf B}_{\gamma}(L))} {L^{6g-6+2n}},$$
and 
Theorem $\ref{vol-ball}$ implies that  $n_{\gamma}=v_{\gamma} \in {\Bbb Q}.$ 

If $\gamma$ is not filling, it is enough to choose $\mathcal{P}$ so that ${\bf S}={\bf S}(\gamma) \subset \mathcal{P}$ as in \ch \ref{non-filling}. For $\Gamma_{\gamma}= \Gamma \cap \Stab({\bf S})$, we can use the same argument to count the orbit  $\Gamma \cdot X$ on 
$\te_{g,n}/ \Gamma_{\gamma}.$ The result follows from the discussion in \ch $\ref{remark:NF}.$

\hfill $\Box$

\subsection{Remark on the case of surfaces with geodesic boundary components}\label{G-B}
The same proof works for hyperbolic surfaces with geodesic boundary components.
 We remark that the ergodicity of the earthquake flow is proved in \cite{M:EE} for cusped surfaces however the same result holds for the moduli spaces of surfaces with geodesic boundary components. Note that as in 
\ch \ref{G-S} there is a natural symplectic form on $\mathcal{M}_{g,n}(L).$ 

In fact, we don't need the ergodicity of the earthquake flow in order to prove Theorem \ref{morege}. We do this just makes the proof simpler. Using the technical Lemmas in \ch \ref{local:calculation}, one would only need the ergodicity of $\Mod(S)$ on $\ML(S)$ \cite{Masur:mapping}.  The limit measure in the proof of 
Theorem \ref{distribution:main} is invariant under the {\it horosphere} foliation. For more on the ergodic property of the action of the mapping class group on $\ML(S)$ for any topological surface $S$ with boundary components see \cite{LMi}. 
 \end{section}

\begin{section}{Appendix on Trigonometry of pairs of pants}\label{App:1}

In this section we use some elementary calculations and estimates on the trigonometry of hyperbolic pairs of pants
to prove Lemma \ref{length:set}, Lemma \ref{basic:F1} and Lemma \ref{local:main}. There is no new ideas involved, but we have included the details for completeness. 

We will use the following elementary estimates in this section: 

\begin{equation}\label{star1}
\arccosh(1+R)=
\left\{
	\begin{array}{ll}
		\sqrt{2R}\; (1+ O(R)), & \mbox{if } R <1 \\
	          \log(2)+\log(R)+ O(\frac{1}{R}) & \mbox{if } R> 1\\
                    
	\end{array}
\right.
\end{equation}

\begin{equation}\label{star2}
\log(\sinh(R))=
\left\{
	\begin{array}{ll}
		\log(R)+ O(R^2) & \mbox{if } R <1 \\
	          R-\log(2)+ O(e^{-R}) & \mbox{if } R>  1\\
                    
	\end{array}
\right.
\end{equation}\\
and 
\begin{equation}\label{star3}
\log(1+R)= R+ O(R^2) \;\;\; \mbox{if}\;\;\; R <1 .                
\end{equation}\\

\noindent
{\bf Proof of Lemma \ref{length:set}.} Following $(\ref{HI})$ and $(\ref{HC1})$,  we have 
$$\cosh(\tilde{c})= \sinh(\tilde{a}) \sinh(\tilde{b}) \cosh(c) +\varepsilon \cosh(\tilde{a}) \cosh(\tilde{b}),$$
where $\varepsilon=-1$ iff $H$ is convex and otherwise $\varepsilon=1.$

We claim that $c_0>0$ can be chosen so that in all the cases $1/2< \tilde{c}/c< 2$:

\begin{enumerate}
\item  Note that in part $({\bf 1})$ and $({\bf 2})$, we have
 $$ c< \log(\cosh(c))+ \log(\sinh(\tilde{a}))+ \log(\sinh(\tilde{b})) < \log(\cosh(\tilde{c})) < 2\tilde{c}.$$  
Since $\tilde{c} < \tilde{a}+ \tilde{b}+ c,$ in part $({\bf 2}),$  $c_0$ can be chosen so that we have  $\tilde{c} < 2 c.$
\item In part $({\bf 3})$, we have
$$ \frac{\tilde{c}}{2}<\tilde{c}-\tilde{a}-\tilde{b} \leq c.$$
Moreover, in this case $c_0$ can be chosen so that 
$$ c < 2 \tilde{c}. $$  
\end{enumerate}

 Let  $T=\sinh(\tilde{a}) \sinh(\tilde{b}) \cosh(c) +\varepsilon \cosh(\tilde{a}) \cosh(\tilde{b}).$
 Now, we can use $(\ref{star1})$ and approximate $\tilde{c}=\arccosh(T)$ by $\log(T)$ with an error term of $O(1/T).$
Then $\tilde{c}- c-\log(\sinh(\tilde{a}))- \log(\sinh(\tilde{b}))$ can be approximated by $$\log(1+ \varepsilon e^{-c} \frac{\cosh(\tilde{a}) \cosh(\tilde{b})}{\sinh(\tilde{a}) \sinh(\tilde{b})}). $$
 Now in view of $(\ref{star3})$ the result follows since in all the cases we have
 $$e^{-c} \; \times \frac{\cosh(\tilde{a}) \cosh(\tilde{b})}{\sinh(\tilde{a}) \sinh(\tilde{b})} < e^{-c/2}.$$

\hfill $\Box$

\subsection{Asymptotic behavior of the function $F_{1}$}\label{F1}

Here we analyze the asymptotic behavior of $F_{1}(x,y,z)$ (defined by $(\ref{def:F1})$) on ${\Bbb R}_{+}^{3}$. These calculations will be used in the proof of Lemma $\ref{local:main}.$ We write

\begin{equation}\label{expand:F}
F_{1}(x,y,z)= \arccosh \left(1+ \frac{2 \; \e^{z-x-y}+ 2 \e^{-z-x-y}+ 2 \e^{-2x}+ 2\e^{-2y}}{ (1-\e^{-2x})( 1-\e^{-2y})}\right). 
\end{equation}

We use the expansion of $\arccosh(1+R)$ at $R=0$ and $R= \infty$ as in $(\ref{star1})$ to investigate the asymptotic properties of $F_{1}(x,y,z).$

We are interested in understanding the asymptotics of $\log(\sinh(F_{1}(x,y,z)))$. The behavior of this function is different in different parts of the cone ${\Bbb R}_{+}^3$. \\

Let 
$$\Delta_{1}=\{ (x,y,z)\; | \; x+y\leq z, \; x,y,z\geq 0 \},$$
$$\Delta_{2}=\{ (x,y,z)\; | \;  |x-y|\leq z\leq x+y, \;x,y,z\geq 0 \},$$
and
$$\Delta_{3}=\{ (x,y,z)\; | \;   z\leq |x-y|, \; x,y,z\geq 0 \}.$$

We define the piecewise linear function $E: {\Bbb R}_{+}^{3} \rightarrow {\Bbb R}$ by 

$$
E(x,y,z) =
\left\{
	\begin{array}{ll}
		z-x-y & \mbox{if }  (x,y,z) \in \Delta_1\\
		(z-x-y)/2 & \mbox{if }  (x,y,z) \in \Delta_2 \\
                     -\min\{x,y\} & \mbox{if } (x,y,z) \in \Delta_3
	\end{array}
\right.
$$
Then we have: 

\begin{lemm}\label{E:estimate}
In terms of the above notation, there exist $c_0, c_1 >0$ such that if $\min\{x,y,z\} >c_0$, then we have 
$$|\log(\sinh(F_1(x,y,z)))- E(x,y,z)|\leq c_1.$$
\end{lemm}
\noindent
{\bf Proof.} 
Note that one can easily bound the denominator of $(\ref{expand:F})$ since for $r>\log(2)$ we have
$$1+ \e^{-r} < \frac{1}{1-\e^{-r}} <  1+ 2 \e^{-r}. $$ 
Next, we will find the dominant term in the numerator.\\

\noindent
{\bf 1.} If $x+y \leq z$ then $1< \arccosh(2)\leq F_{1}(x,y,z).$ Then $(\ref{star1})$ implies that 
\begin{equation}\label{F1-1}
 F_{1}(x,y,z)= (z-x-y)+ \log(2)+ O( \frac{1}{\e^{z-x-y}}+ \e^{-\min\{x,y,z\}}+1).
 \end{equation}
Now the bound follows from $(\ref{star2}).$\\

\noindent
{\bf 2.} If $|x-y| < z <x+y,$ then since $-2 \min\{x,y\} < z-x-y <0$, $\e^{z-x-y}$ is again the leading term in $(\ref{expand:F}).$ Following $(\ref{star1})$ we have
$$F_{1}(x,y,z)= 2 \e^{(z-x-y)/2} (1+  O(\e^{-\min\{x,y,z\}}+ \e^{3/2(z-x-y)}+1).$$

Then $(\ref{star3})$ implies that:
\begin{equation}\label{F1-2}
\log(F_{1}(x,y,z))= (z-x-y)/2+ \log(2)+ O(\e^{-\min\{x,y,z\}}+ \e^{3/2(z-x-y)}+1).
\end{equation}

\noindent
{\bf 3.} Finally, if $z < |x-y|$ then $z-x-y < -2\min\{x,y\}$. In this case, $\e^{-2 \min\{x,y\}}$ is the leading term in $(\ref{expand:F})$ . Hence by $(\ref{star1})$ we have 
$$ F_{1}(x,y,z)= \sqrt{2} \e^{-\min\{x,y\}} (1+  O(\e^{-\min\{x,y,z\}}+1)).$$
Then $(\ref{star3})$ implies that:
\begin{equation}\label{F1-3}
\log(F_{1}(x,y,z))= -\min\{x,y\}+ \log(2)/2+  O(\e^{-\min\{x,y,z\}}+1).
\end{equation}
\hfill $\Box$\\

\noindent
{\bf Proof of Lemma \ref{basic:F1}.} The argument is an elementary case by case analysis using 
lemma $\ref{E:estimate}.$ 
We claim that 
\begin{equation}\label{E-1}
|E(x,y,z)-E(x',y',z')| \leq 20 \; \epsilon \times \max\{x,y,z\}.
\end{equation}

Note that:
\begin{itemize} 
\item
If $(x',y',z') \in B_{\epsilon}(x,y,z)$ then
\begin{equation}\label{C-1}
| (x'+y'-z')- (x+y-z)| \leq 3 \; \max\{x,y,z\} \times \epsilon.
\end{equation} 
Also, easy to check that 
\begin{equation}\label{M-1}
1/(1+\epsilon)< \frac{\min\{x,y\}}{\min\{x',y'\}} < 1+\epsilon.
\end{equation}
\item By the definition $E(x,y,z)\geq 0$ if and only if $(x,y,z) \in \Delta_1.$
\end{itemize}

To simplify the notation, let $D=E(x,y,z)-E(x',y',z').$ Then we have:\\

\noindent
{\bf 1.} If $(x,y,z), (x',y',z') \in \Delta_i$ for some $1 \leq i \leq 3,$ then ($\ref{E-1}$) follows from 
$(\ref{C-1})$ and $(\ref{M-1}).$

\noindent
{\bf 2.}
Assume that $(x,y,z) \in \Delta_1,$ and $(x',y',z') \in \Delta_2$ or $\Delta_3.$ In this case, it is enough to bound both $E(x,y,z)$ and $|E(x',y',z')|$ from above. In this case  ($\ref{C-1}$) implies that
$E(x,y,z)+z'-x'+y' \leq 4 \; \max\{x,y,z\} \times \epsilon.$

Now if $(x',y',z') \in \Delta_2$, this is the same as that $E(x,y,z)+2|E(x',y',z')| \leq 4 \; \max\{x,y,z\} \times \epsilon.$ which is enough for us. 

 If $(x',y',z') \in \Delta_3,$ we can assume without loss of generality that $x'>y'.$
 By definition,  $z'< x'-y'< x'.$
 Then get that $0<z-x= (z-x+x'-z')+ z'-x' < |z-z'+x'-x| < 2 \epsilon \max\{x,y,z\} $. The bound follows  since : $0< (z-x-y)+y' < (z-x)+ |(y'-y)|.$\\

\noindent
{\bf 3.}
If $(x,y,z) \in \Delta_2$, $(x',y',z') \in \Delta_3$,  
\begin{itemize}
\item 
If $y<x$, and $y'< x' $ then $2D=(z-x-y)+2y'$. Note that since $z'-x'+y'<0,$ we have
$$(z-x-y)+2y'= (z'-x'+y' )+ (z-z')-(x-x')-(y-y') < |z-z'|+ |x-x'|+ |y-y'|. $$
Moreover, in order to obtain the lower bound on $D$, by $\ref{M-1}$ it is enough to note that we have $$ -2y+2y' <  (z-x-y)+2y' .$$ 

\item 
If $ y<x$, and $x'< y',$ then $2 D=(z-x-y)+2x'$. In this case $z'-y'+x'<0 $, and the upper bound exactly the same as before
$$ 2D= (z'-y'+x')+ (x'-x)+(y'-y)+(z-z') < |z-z'|+|y-y'|+|x-x'|.$$
For the lower bound, we have 
 $$ -2y+2x' < (z-x-y)+2x' .$$ On the other hand, we have  
 $$  \frac{1}{1+\epsilon} y   < \frac{1}{1+\epsilon} x < x' ,$$
 and hence 
 $$ -x' \epsilon< x'- y .$$ 
 \end{itemize}
Now the result follows from $(\ref{E-1})$ and Lemma $\ref{E:estimate}$  if we have $\max \{\frac{c_1}{\epsilon}, c_0\} <L_0,$ where $c_0$ and $c_1$ are given by Lemma $\ref{E:estimate}.$

\hfill $\Box$\\

We will use Lemma \ref{length:set} and Lemma \ref{basic:F1} to prove Lemma \ref{local:main}.\\

\noindent
{\bf Proof of Lemma \ref{local:main}.} Here we prove the statement for ${\bf T}_{0,4};$ the proof for ${\bf T}_{0,1}$ is similar. Assume that $\alpha$ separates the surface into two pairs of pants $\Sigma_1$ with boundary components $\alpha_1, \alpha_4$ and $\alpha$ and $\Sigma_2$  with boundary components $\alpha_3, \alpha_2$ and $\alpha.$ 
 In view of the assumptions $1$, $2$, without loss of generality, we can assume that
 \begin{enumerate}
 \item 
 the length $|a|$ of the shortest geodesic arc $a$ joining $\alpha_1$ to $\alpha$ on $Y$ is $\leq \e^{-L/M}.$ 
 \item the length $|b|$ of the shortest geodesic arc $b$ joining $\alpha_3$ to $\alpha$ on $Y$ is $\leq \e^{-L/M}.$ Note that we also have $\e^{-2L} \leq |a|,$ and 
 $\e^{-2L} \leq |b|.$
 \end{enumerate}
 
Let $d_k$ be an arc joining $\alpha_{1}$ to $\alpha_3$, homotope to the arc going along $a$ then going to around $\alpha,$ $k$ times and then going along $b$ with reverse orientation. 
Following Lemma \ref{length:set} ({\bf 1}), using $(\ref{HI})$, one can show that if 
$k_0 =10\cdot M^2,$  $d=d_{k_0} > 2 M L$. 
Consider the corresponding closed curve $\beta_{k},$ ($ k \geq k_0$) and $\beta= \beta_{k_0}.$ Then $ \beta_{k}= h_{\alpha}^{k-k_0}(\beta).$

Moreover, Lemma \ref{length:set} ({\bf 1}) could be used to relate $\tau_{\alpha}(Y)$ (resp. $\tau_{\alpha}(Y)$) to $d=|d_{k}|$:
\begin{equation}\label{part-one}
 d= \log(\sinh(a))+ \log(\sinh(b))+ T+ \log(2)+ O(\e^{-L/M}),
 \end{equation}
$T= \tau_{\alpha}(Y)+ k_0 \ell_{\alpha}(Y).$

Note that the length of $\beta$ can be approximated in terms of $d$, $ \ell_{\alpha_1}$ and $\ell_{\alpha_3}$ using Lemma \ref{length:set} $({\bf 3}).$
The first part of the Lemma follows from Lemma \ref{length:set}.

On $Y'$ let $|a'|$ and $|b'|$ denote the length of the shortest geodesic arcs joining  joining $\alpha_1$ to $\alpha$ and $\alpha_3$ to $\alpha$.  Also, let $d'=d_{k_0}'$ be the length of the arc corresponding to $d=d_{k_0}.$
Note that 
$$ \frac{1}{1+\epsilon}< \frac{\ell_{\beta_{k}}(Y')}{\ell_{\beta_{k}}(Y)} < 1+\epsilon .$$ 
We consider the three pairs of pants formed by $(\beta, \alpha_3, \alpha_1)$, $(\alpha, \alpha_3, \alpha_2)$ and $(\alpha, \alpha_1,\alpha_4)$. In term of the notation used in \ch \ref{Tr-h}, we have:

$$ (\ell_{\alpha_1}(Y'), \ell_{\alpha_3}(Y'), \ell_{\beta}(Y')) \in B_{\epsilon} (\ell_{\alpha_1}(Y), \ell_{\alpha_3}(Y), \ell_{\beta}(Y)),$$
$$ (\ell_{\alpha_1}(Y'), \ell_{\alpha_4}(Y'), \ell_{\alpha}(Y')) \in B_{\epsilon} (\ell_{\alpha_1}(Y), \ell_{\alpha_4}(Y), \ell_{\alpha}(Y)),$$
and 
$$ (\ell_{\alpha_2}(Y'), \ell_{\alpha_3}(Y'), \ell_{\alpha}(Y')) \in B_{\epsilon} (\ell_{\alpha_2}(Y), \ell_{\alpha_3}(Y), \ell_{\alpha}(Y)).$$

In view Lemma \ref{basic:F1}, we have: 
\begin{equation}\label{part-two}
 | \log(\sinh(d')) -\log(\sinh(d)) | \leq 5 \cdot M \epsilon L, 
 \end{equation}
$$ | \log(\sinh(a')) -\log(\sinh(a)) | \leq 5 \cdot M \epsilon L,$$
and 
$$
  | \log(\sinh(b')) -\log(\sinh(b)) | \leq 5 \cdot M \epsilon L.$$

In view of $(\ref{star2})$ since $\epsilon M <1/M,$ and $d > 2 M L,$ we have $d' > M L.$ 
Similarly, we have $-3L < \log(\sinh(a')) < -\frac{L}{2M}$ and $-3L < \log(\sinh(b')) < -\frac{L}{2M}.$
By Lemma \ref{length:set} $({\bf 2})$, we have

\begin{equation}\label{part-three}
 d'= \log(\sinh(a'))+ \log(\sinh(b'))+ T'+ \log(2)+ O(\e^{-L/M}),
 \end{equation}
where $T'= \tau_{\alpha}(Y')+ k_0 \ell_{\alpha}(Y').$ 

On the other hand $|\ell_{\alpha}(Y)- \ell_{\alpha}(Y)| \leq \epsilon M L.$ Hence in view of $(\ref{star2})$, $(\ref{twist:change})$ follows from Lemma $\ref{basic:F1},$ $(\ref{part-one}),$ 
$(\ref{part-two}),$ and $(\ref{part-three}).$

\hfill $\Box$\\

\end{section}

\bigskip
Department of Mathematics, Stanford University, Stanford CA 94305 USA;\\
mmirzakh@math.stanford.edu

\end{document}